\documentclass[preprint]{imsart}
\RequirePackage{amsthm,amsmath,amsfonts,amssymb}
\usepackage{dsfont}
\usepackage{pdfsync}
\usepackage{graphicx, subcaption}
\usepackage{float}
\RequirePackage[numbers,sort&compress]{natbib}
\usepackage{bbm}
\usepackage[nothm]{yk}
\RequirePackage[colorlinks,citecolor=blue,urlcolor=blue]{hyperref}

\newcommand{\R}{\mathbb{R}}
\newcommand{\ord}{\mathrm{ord}}
\newcommand{\mult}{\mathrm{mult}}

\startlocaldefs
\theoremstyle{plain}

\newtheorem{theorem}{Theorem}[section]
\newtheorem{lemma}[theorem]{Lemma}
\newtheorem{corollary}[theorem]{Corollary}
\theoremstyle{remark}
\newtheorem{definition}[theorem]{Definition}
\newtheorem{assumption}[theorem]{Assumption}

\newtheorem{conjecture}[theorem]{Conjecture}

\endlocaldefs

\begin{document}
\begin{frontmatter}
\title{Minimax Optimal rates of convergence in the shuffled regression, unlinked regression, and deconvolution under {vanishing noise}}
\runtitle{Shuffled and unlinked regression models under vanishing noise}
\begin{aug}
\author[A]{\fnms{Cecile}~\snm{Durot}\ead[label=e1]{cecile.durot@parisnanterre.fr}} and 
\author[B]{\fnms{Debarghya}~\snm{Mukherjee}\ead[label=e2]{mdeb@bu.edu}\orcid{0000-0000-0000-0000}}
\address[A]{MODAL'X, UPL, Univ. Paris Nanterre, CNRS \printead[presep={,\ }]{e1}}

\address[B]{Department of Mathematics and Statistics,
Boston University\printead[presep={,\ }]{e2}}
\end{aug}

\begin{abstract}
Shuffled regression and unlinked regression represent intriguing challenges that have garnered considerable attention in many fields, including but not limited to ecological regression, multi-target tracking problems, image denoising, etc. However, a notable gap exists in the existing literature, particularly in vanishing noise, i.e., how the rate of estimation of the underlying signal scales with the error variance. 
This paper aims to bridge this gap by delving into the monotone function estimation problem under vanishing noise variance, i.e., we allow the error variance to go to $0$ as the number of observations increases. 
Our investigation reveals that, asymptotically, the shuffled regression problem exhibits a comparatively simpler nature than the unlinked regression; if the error variance is smaller than a threshold, then the minimax risk of the shuffled regression is smaller than that of the unlinked regression. On the other hand, the minimax estimation error is of the same order in the two problems if the noise level is larger than that threshold.
Our analysis is quite general in that we do not assume any smoothness of the underlying monotone link function. 
Because these problems are related to deconvolution, we also provide bounds for deconvolution in a similar context. 
Through this exploration, we contribute to understanding the intricate relationships between these statistical problems and shed light on their behaviors when subjected to the nuanced constraint of vanishing noise.
\end{abstract}



\end{frontmatter}


\section{Introduction}
In a standard supervised learning problem, we observe finitely many (say $n$ many) realizations of a paired random variable $(X, Y)$ (where $Y$ is the response variable and $X$ is the covariate) and construct a predictor of $Y$ from $X$ based on those realizations. 
This setup encompasses widely used statistical models, such as standard regression and classification. 
A fundamental assumption for this class of problems is the availability of \emph{paired} observations $\{(X_i, Y_i)\}_{1 \le i \le n}$, and, in particular, the knowledge of which $X$ corresponds to which $Y$.  
However, the information about this pairing is unavailable in many real-world applications. 
To illustrate, consider a concrete example originated in \cite{degroot1971matchmaking} and elaborated in \cite{slawski2022permuted}: suppose we have a collection of photographs of $n$ actors at young ages (denoted as $\cX_n$) and a collection of photos of the same group of actors at older ages (denoted as $\cY_n$). The challenge is to match the young-age photographs to the corresponding old-age pictures of the same actors. In this scenario, if we lack labels such as the actors' names in the images, we face the challenge that each $X_i \in \cX_n$ corresponds to a certain $Y_j \in \cY_n$, yet the specific pairings are unknown. This type of data is called \emph{shuffled data}, and the corresponding problem is termed as \emph{shuffled regression/classification} problem. 
Another common application of shuffled regression, frequently encountered in ecological regression, involves the analysis of voting patterns among groups of individuals. For instance, consider a scenario where two competing political parties, party A and party B, exist within a region of $n$ voters. Each voter is represented by a pair of variables $(X_i, Y_i)$, where $X_i$ is their background information, and $Y_i$ is a binary indicator for their vote (with $Y_i = 1$ denoting a vote for party A and $Y_i = 0$ indicating a vote for party B). We do not observe $\{Y_1, \dots, Y_n\}$, i.e., the specific votes of the individuals, but we have information about the total votes garnered by parties A and B in the region. Consequently, we possess a permuted version of ${Y_1, \dots, Y_n}$ to ensure privacy and anonymity. There is a vast body of literature on a variety of applications and analyses of ecological regression; interested readers may consult \cite{freedman1991ecological, jiang2020ecological, kousser1973ecological, wakefield2003sensitivity, brown1986aggregate, crewe1976another} and references therein. 
In addition to ecological regression, permuted or shuffled data are frequently encountered in various other fields of research. 
For example, the multi-target tracking problem (\cite{poore2006some, durrant2006simultaneous}), aims to match measurements obtained from different environments. 
In computer vision, shuffled regression is used to match two different sets of images of a group of objects, where one set of images may be a noisier/rotated/distorted version of the other (e.g., see \cite{hartley2003multiple, szeliski2022computer} for more details). 
In \cite{ma2021optimal}, the authors studied a permuted monotone matrix model motivated by bacterial growth dynamics based on genome assemblies. 
Another application of the shuffled regression lies in the seriation problem is statistics, as analyzed in \cite{flammarion2019optimal}. The seriation problem has found its application in various fields, including but not limited to sociology (\cite{forsyth1946matrix}), biology (\cite{sokal1963principles}), archaeology (\cite{petrie1899sequences}), and anthropology (\cite{czekanowski1909differentialdiagnose}). 

So far, we have presented examples where we have access to $\{X_1, \dots, X_n\}$ and $\{Y_1, \dots, Y_n\}$, but we do not know which $X_i$ corresponds to which $Y_i$. 
However, in many other applications, $\cX_n$ and $\cY_n$ may be generated from different units. 
Consider the example of finding the relation between wages and housing transactions as presented in \cite{carpentier2016learning}. 
In general, we expect a monotonically increasing relationship between these two variables, i.e., individuals with higher incomes are likely to pay more for housing. 
However, it is almost impossible to find a set of data in which we have information about the income and housing transactions of the same group of people, as different agencies typically collect information regarding wage and housing prices. 
Therefore, we have some realization of $X$ variable (wage), say $\cX_m= \{X_1, \dots, X_{m}\}$ and some other realizations of $Y$ (housing price), say $\cY_n = \{Y_1, \dots, Y_{n}\}$. 
Furthermore, $\cX_m$ and $\cY_n$ can be independent as they can be collected from different groups of people or may contain some overlap, i.e., some individuals are present in both groups. 
This type of data is called \emph{unlinked data}, and the corresponding statistical problem to comprehend the relation between $X$ and $Y$ is called \emph{unlinked regression/classification}. 
Estimation of monotone link function (i.e., when $\bbE[Y \mid X] = m_0(X)$ for some monotone function of $X$) using unlinked data is even less explored in the literature, with notable exceptions are \cite{carpentier2016learning} and \cite{balabdaoui2021unlinked}. In \cite{balabdaoui2021unlinked}, the authors first rigorously established that the monotone link function is identifiable and proposed an estimator of $m_0$ under various smoothness assumptions on the tail of the characteristic function of the error distribution. 

A common thread between shuffled and unlinked regression lies in the absence of links between the pair $(X, Y)$, i.e., we do not know which $X$ corresponds to which $Y$. 
This poses a natural question regarding the informativeness of these unpaired $X$'s in estimating the signal/mean function. 
The key distinction between shuffled and unlinked regression is that in shuffled regression, we possess the set of $X_i$'s corresponding to the set of $Y_i$'s. In contrast, in unlinked regression, $X_i$'s can be entirely independent of $Y_i$'s, i.e., they can be obtained from two separate datasets. 
Consider the following additive noise model with a monotone link function studied in this paper: 
\begin{equation}
\label{eq:main_model_equation}
Y_i = m_0(X_i) + \sigma_n \delta_i \,,
\end{equation}
where $\sigma_n$ represents the scale of the error that may change with $n$.  
At first blush, it appears that shuffled regression provides \emph{more} information than unlinked regression, which is indeed true for the noiseless case where $Y = m_0(X)$; in shuffled regression, we can recover the precise link between $X_i$ and $Y_i$ as $Y_{(i)} = m_0(X_{(i)})$ for all $1 \le i \le n$, where $X_{(j)}$ (resp. $Y_{(j)}$) denotes the $j^{th}$ order statistic of $X_i$'s (resp. $Y_i$'s). Conversely, no such recovery is possible for unlinked regression, even in the absence of noise. However, at the other extreme (i.e., when the noise variance is of the order $1$), whether the shuffled regression model holds an advantage over the unlinked regression model for estimating the underlying mean function is still being determined. 
Hence, it is natural to anticipate that a threshold exists regarding the magnitude of noise: if the standard deviation of the noise falls below that threshold, shuffled data is more informative than unlinked data, resulting in a fast estimation of the signal. 
This observation is the key motivation of our paper: to quantify the precise dependence of the noise variance on the estimation rate of $m_0$ and determine a threshold below which the shuffled regression model exhibits a faster rate for estimating $m_0$ compared to the unlinked regression. 

Estimation of the mean function under the shuffled and the unlinked model shares a resemblance with the classical deconvolution problem. 
The relationship  between the models we are interested in and the deconvolution model has also led us to study the deconvolution model from a minimax point of view under vanishing noise.
In deconvolution, the parameter of interest is the distribution (or its functional) of some random variable $Z$ (say, the signal of some experiment). 
However, we do not observe any realization of $Z$; rather, we observe $n$ independent observations $Y_1, \dots, Y_n$ generated as $Y_i = Z_i + \eps_i$ where $\eps_i$ is independent from $Z_i$. 
In other words, we have access to the realizations of a random variable $Y$, whose distribution is a convolution between the distribution of $Z$ (parameter of interest) and the distribution of error $\eps$  (to relate it to \eqref{eq:main_model_equation}, we may think $Z_i = m_0(X_i)$ and $\eps = \sigma_n \delta$). 
Therefore, one needs to \textit{deconvolute} the 
distribution of $\eps$ from the distribution of $Y$ to get back to the distribution of $Z$. 
Analysis of such deconvolution problems boasts a rich historical background and has found extensive applications across various domains. Examples include its use in medical imaging (\cite{michailovich2007blind, katohrecent}), astronomy (\cite{starck2002deconvolution, jefferies1993restoration, li2011application}), signal processing (\cite{javidi1989image, mendel2012maximum, javidi1989multifunction}), image deblurring (\cite{satish2020comprehensive, yuan2007image, hall2007blind}), and volatility estimation (\cite{comte2004kernel, van2011estimation, miguel2022volatility}).
From a statistical perspective, the deconvolution problem was initially studied in a series of papers, including \cite{fan1991optimal,Fan91AsympNorm,fan1992deconvolution,stefanski1990deconvolving,zhang1990fourier,carroll1988optimal}, where the authors establish the minimax optimal rate of the deconvolution problem depending on the tail of the characteristic function of $\eps$ and construct estimators of the distribution/density of $Z$ that are minimax rate optimal, primarily via a kernel smoothing technique under certain smoothness assumptions on the distribution of $Z$. 
The loss function typically used is the $L_p$ distance between the estimated density or distribution of $Z$ and the true one. 
This line of research established that the rate of convergence is particularly slow if the tail of the characteristic function of the error is very thin (e.g., Gaussian distribution), where the rate decays polynomially in $\log{n}$ ($n$ is the number of observations), e.g., Theorem 1 and 4 of \cite{fan1991optimal}. 
Interested readers may consult \cite{crowder2010deconvolution} for further references. 
More recently, deconvolution with respect to the Wasserstein metric has gained significant attention in the literature. 
One key reason is that Wasserstein distance is a flexible metric for the space of all probability distributions. 
It provides a non-trivial notion of distance between two measures even when they have non-overlapping/disjoint support. Furthermore, the Wasserstein deconvolution problem is related to geometric inference, which aims to recover various geometric properties (e.g., Betti number) of a set upon observing some points from that set (e.g., see \cite{chazal2011geometric, caillerie2013deconvolution}). The minimax optimal rate of Wasserstein deconvolution in the presence of supersmooth error (error whose characteristic function has a thin tail, precisely defined later) in any dimension has been established in \cite{dedecker2013minimax} and for ordinary smooth error (errors whose characteristic function has a thick tail) in one dimension is investigated in \cite{dedecker2015improved}. As before, errors with thin-tailed characteristic functions lead to a very slow rate of convergence for the Wasserstein deconvolution problem.      

Although theoretically, the rate of deconvolution can be very slow, 
general deconvolution estimators seem to work well in practice, even in the presence of noise with thin-tailed characteristic functions. One way to explain this phenomenon, as illustrated in \cite{delaigle2008alternative}, is that the scale of the error is small. In other words, it is intuitive to expect that if the standard deviation of the error is small, then estimation of the signal density should be easier. The rate of convergence of the deconvolution problem involving the error variance was initially investigated in \cite{fan1992deconvolution}. Later, \cite{delaigle2008alternative} established the rate of convergence of a deconvolution kernel density estimator when the noise variance $\sigma_n^2$ goes to $0$ as $n \uparrow \infty$, and the distribution of the signal is smooth. However, the minimax optimal rate of estimation with respect to the variance of the error is not known. In this paper, we bridge this gap. As in \cite{delaigle2008alternative}, we assume that observations take the form $Y_i = Z_i + \sigma_n \delta_i$ where we assume $\sigma_n \downarrow 0$ as $n \uparrow \infty$ (similar to \eqref{eq:main_model_equation} with $m_0(X_i)$ replaced by $Z_i$), but we do not make any smoothness assumption on the signal distribution. Here, we consider the Wasserstein deconvolution problem, i.e., we aim to estimate the distribution of $Z$ and measure the estimation error with respect to the Wasserstein metric. 
%
We summarize our contributions below: 
\\\\
\noindent
{\bf Our contributions: }Our primary contribution in this paper is to present an asymptotic comparison of the minimax risks of the unlinked regression and the shuffled regression problems under a vanishing noise condition. In particular, we study the model \eqref{eq:main_model_equation} under the assumption that the noise standard deviation  $\sigma_n$ goes to $0$ as $n \uparrow \infty$, without any smoothness assumption on the link function $m_0$, and quantify the minimax risk in terms of $(n, \sigma_n)$. We also construct an estimator that attains the minimax risk. Our analysis reveals a phase-transition phenomenon: there exists a threshold of order (close to) $ n^{-1/2}$, such that if $\sigma_n$ is below that threshold, then the minimax risk of estimation of $m_0$ in shuffled regression is smaller than that of the unlinked regression. On the other hand, the minimax risks are of the same order in the two problems when $\sigma_n$ is larger than the threshold. 
Furthermore, we also establish the minimax risk (in terms of $(n, \sigma_n)$) for estimating the signal distribution in the Wasserstein-deconvolution problem when the error standard deviation $\sigma_n \downarrow 0$ as $n \uparrow \infty$ and the characteristic function of $\delta$ has a thin tail. Our theory also indicates that the minimax risk of the unlinked regression and deconvolution has the same order irrespective of the value of $\sigma_n$, i.e., both problems showcase a similar level of difficulty.  
\\\\
{\bf Organization of the paper: }We organize our paper as follows: in Section \ref{sec: deconvolutionSM}, we analyze the deconvolution problem under vanishing noise conditions and establish both minimax upper and lower bounds on the estimation of the distribution of the signal. In Sections \ref{sec:unlinked} and \ref{sec:shuffled}, we present our results for the unlinked and shuffled regression models, respectively. In Section \ref{sec:simulation}, we present some simulation results in the context of a conjecture (Conjecture \ref{conj}). Finally, we end this paper with our conclusion and some future research direction in Section \ref{sec:conclusion}. The proofs of the main results are postponed to Appendix \ref{app:main} whereas the proofs of intermediate results are given in Appendix \ref{app:lemmas}.
\section{Deconvolution under vanishing noise}
\label{sec: deconvolutionSM}

This section presents our main results for the deconvolution problem with vanishing noise. In this setting, we observe $\cY_n = \{Y_1, \dots, Y_n\}$ where $Y_i = Z_i + \sigma_n \delta_i$ for some unobserved $\cZ_n = \{Z_1, \dots, Z_n\} \overset{i.i.d.}{\sim} \mu_Z$ and $\{\delta_1, \dots, \delta_n\} \overset{i.i.d.}{\sim} \mu_\delta$ independent of $\cZ_n$. Here, we assume $E[\delta_i] = 0$, $\var(\delta_i) = 1$, and  $\sigma_n$ is a bounded sequence that is allowed to approach zero as $n \uparrow \infty$. Moreover, as is customary in deconvolution problems for identifiability purposes, $\sigma_n$, and the distribution $\mu_\delta$ are assumed to be known.
We address the problem of Wasserstein deconvolution, i.e., we aim to estimate $\mu_Z$ (the distribution of the signal) based on observed $Y_i$'s and evaluate the performance of an estimator via $W_1$-Wasserstein distance between that estimator and $\mu_Z$. Recall that, given any two distributions $(\mu, \mu')$ on $\mathbb{R}$, the $W_1$ distance between them is defined as: 
$$
W_1(\mu, \mu') = \inf_{\pi \in \Pi(\mu, \mu')} \int \int |x - y| \pi(dx, dy), 
$$
where $\Pi(\mu, \mu')$ is the collection of all possible coupling $\pi$ with the marginals $(\mu, \mu')$.  

Before presenting the main results, we first lay the setup. 
Our first assumption is on the structure of the distribution of $Z$, namely $\mu_Z$, in terms of the finiteness of its first few moments as elaborated below:  
\begin{assumption}[Moments of $\mu_Z$]
\label{assm:dist_Z}
Let $\cD(M, a)$ be the set of all distributions $\mu$ whose $(a+2)^{th}$ moment is bounded by $M$. Assume that $\mu_Z \in \cD(M, a)$ for some $a,M>0$.
\end{assumption} 
The above assumption is quite general in that it only requires a few finite moments for the signal distribution. Unlike several other research papers, we do not need the distribution of $Z$ to be smooth or even have a density function.
For any distribution $\mu$, we denotes its characteristic function by $\mu^*(t)$, defined as: 
$$
\mu^*(t)= \bbE_{X \sim \mu}[\exp{(itX)}] = \int e^{itx} \ d\mu(x) \,.
$$
The rate of estimation of $\mu$ in this deconvolution problem typically depends on the tail of the characteristic function of the error distribution $\mu_\delta$ (e.g., see \cite{fan1991optimal}). 
Roughly speaking, if the tail of $\mu_\delta^*(t)$ is \emph{thin}, then the deconvolution problem is harder than when the tail of  $\mu_\delta^*(t)$ is \emph{thicker}. Based on the thickness of the tail of the characteristic function, the error distribution in the deconvolution problem can be broadly classified into two categories: i) supersmooth error and ii) ordinary smooth error (e.g., see \cite{fan1991optimal,dedecker2013minimax}). In this paper, we work with the supersmooth error, which is defined below: 
\begin{definition}[Supersmooth distribution]
\label{def:supersmooth}
The distribution $\mu$ is said to be supersmooth if its characteristic function $\mu^*$ satisfies the following bound for some positive $\gamma_1,c_1,\beta$, and $\tilde\beta\in\R$: 
\begin{equation*}
\vert \mu^*(t)\vert (1+\vert t\vert)^{-\tilde\beta}\exp(\vert t\vert^\beta/\gamma_1)\leq c_1
 \mbox{ for all }t\in\R.
\end{equation*}
\end{definition}
In other words, a distribution is supersmooth if the modulus of its characteristic function decays exponentially as $t \uparrow \infty$. Henceforth, we assume that the distribution of $\delta$ is supersmooth: 
\begin{assumption}
\label{assm:supersmooth}
The distribution of $\delta$ is supersmooth with some parameter $(c_1, \tilde \beta, \beta, \gamma_1)$ as mentioned in Definition \ref{def:supersmooth}. 
\end{assumption}
 
In the sequel, we denote by $\mu_\epsilon$ the distribution of $\epsilon_i:=\sigma_n\delta_i$ so that the distribution of the sample $\cY_n$ is $(\mu_Z\star\mu_\epsilon)^{\otimes n}$ where $\star$ stands for the convolution operation, and $\otimes n$ denotes the $n$-fold product measure. Under the above assumption on the distribution of $\delta$, we have the following minimax lower bound on the Wasserstein deconvolution problem: 
 \begin{theorem}
 \label{theo: DMlower}
Let $(\sigma_n)$ be a bounded sequence of positive numbers and the distribution of $Z$ satisfies Assumption \ref{assm:dist_Z} for some $a>0$ and $M>0$.
 If the error distribution $\mu_\delta$ satisfies Assumption \ref{assm:supersmooth} and $E\vert \delta\vert^{2+a}<\infty$  then for  sufficiently large $M$, there exists $C>0$  such that for all $n$, 
\begin{equation}
\label{eq: DM16'}
\inf_{\hat \mu_n}\sup_{\mu\in\mathcal D(M,a)}E_{(\mu\star\mu_\epsilon)^{\otimes n}}W_1(\mu,\hat \mu_n)\geq C\left(\sigma_n\left(\log n\right)^{-1/\beta} + n^{-1/2}\right) 
\end{equation}
where  the infimum is over all estimators $\hat \mu_n$ of the distribution of $Z$ based on the observed responses $\cY_n$. 
\end{theorem}
Theorem \ref{theo: DMlower} quantifies a lower bound on the minimax risk of the Wasserstein deconvolution problem in terms of $(n, \sigma_n)$. 
It is evident from this theorem that the lower bound is of the order of $\sigma_n\left(\log n\right)^{-1/\beta}$ if $\sigma_n\geq n^{-1/2}(\log n)^{1/\beta}$, and of the order of $n^{-1/2}$ otherwise.
It is instructive to compare our lower bound with Theorem 2 of \cite{dedecker2013minimax}. In that theorem, the authors obtained the lower bound $(\log{n})^{-p/\beta}$ for $W_p^p$ metric. Setting $p = 1$, their results yield the lower bound $(\log{n})^{-1/\beta}$, which coincides with ours when $\sigma_n = 1$. The key difference between our result and that of \cite{dedecker2013minimax} is the quantification of the precise effect of $\sigma_n$ in the lower bound. Our lower bound shows that the smaller the value of $\sigma_n$, the faster one can estimate $\mu_Z$ based on $\cY_n$. 
However, as the theorem shows, the fastest achievable rate is $n^{-1/2}$. This is intuitive because $n^{-1/2}$ is the price one must pay even when $\sigma_n = 0$, i.e. $Y_i = Z_i$, as $n^{-1/2}$ is the fastest possible rate at which one can estimate $\mu_Z$ upon observing $n$ i.i.d. observations $Z_1, \dots, Z_n$. 

Now that we have a lower bound for the minimax risk of the Wasserstein deconvolution problem, we turn to establish an upper bound. For the upper bound, we need an additional assumption on $\mu^*_\delta$, the characteristic function of $\delta$. Towards that end, define  $r_\delta=1/\mu_\delta^*$, the reciprocal of the characteristic function, and for $\ell=0,1,2$, we denote by $r_\delta^{(\ell)}$  the $\ell$-th derivative of $r_\delta$ (hence, $r_\delta^{(0)}=r_\delta$). We make the following assumption on the function $r_\delta$ and its derivatives: 
\begin{assumption}
\label{assm:error_char_deriv}
The distribution of $\delta$ has a symmetric density around zero with a positive characteristic function, and there exist positive $\gamma_2$, $c_2$, $\beta$, $\tilde\beta$,  such that for every $\ell\in\{0,1,2\}$ and every $t\in\R$,
  \label{assm:DMU}
\begin{equation*}
\vert r_\delta^{(\ell)}(t)\vert \leq c_2 (1+\vert t\vert^{\tilde\beta})\exp(\vert t\vert^\beta/\gamma_2).
\end{equation*}
\end{assumption}
The above assumption is primarily motivated by the assumption of Theorem 4 of \cite{dedecker2013minimax}. One example is the Gaussian distribution, which satisfies the above constraint. Our following theorem presents an upper bound on the minimax risk of the Wasserstein deconvolution problem.

\begin{theorem}
\label{theo: DFMupper}
Let $(\sigma_n)$ be a bounded sequence of positive numbers,
 $a,M,\eta$  be positive numbers where $\eta$ is arbitrarily small. Assume that the distribution of $Z$ satisfies Assumption \ref{assm:dist_Z}, $E\vert \delta\vert^{a+2}\leq M$, and $r_\delta$ is twice continuously differentiable and satisfies Assumption \ref{assm:error_char_deriv}. Then, there exists  a constant $L>0$ such that for any small $\eta > 0$: 
\begin{equation*}
\inf_{\hat \mu_n}\sup_{\mu\in\mathcal D(M,a)}E_{(\mu\star\mu_\epsilon)^{\otimes n}}W_1(\mu,\hat\mu_n)\leq Lv_n
\end{equation*}
where 
\begin{eqnarray}
\label{eq:def_v_n}
v_n= \begin{cases}
\sigma_n\left(\log\left( n\sigma_n^2 \log n\right)\right)^{-1/\beta}&\mbox{ if }\sigma_n\geq n^{-1/2}(\log n)^{\eta}\\
n^{-1/2}(\log\log n)^{-1/\beta}(\log n)^{\eta}&\mbox{ if }\sigma_n\in(n^{-1/2},  n^{-1/2}(\log n)^{\eta})\\
n^{-1/2}&\mbox{ if }\sigma_n\leq n^{-1/2}
\end{cases}
\end{eqnarray}
\end{theorem}

In the above theorem, one may take $\eta > 0$ to be as close to $0$ as possible. 
Combining Theorems \ref{theo: DMlower} and \ref{theo: DFMupper} gives the minimax rates in the deconvolution model with supersmooth noise for an extensive range of noise levels that are allowed to tend to zero as $n\to\infty$, provided that the noise distribution satisfies simultaneously Assumptions \ref{assm:supersmooth} and \ref{assm:error_char_deriv}. 
As mentioned previously, centered Gaussian distribution satisfies both assumptions with $\beta=2$. Let us now compare the upper and lower bounds of the deconvolution problem as obtained in Theorem \ref{theo: DMlower} and Theorem \ref{theo: DFMupper} to elaborate on when they match exactly and when there is a gap: 
\begin{enumerate}
\item If $\sigma_n\geq n^{-1/2+\eta}$ or $\sigma_n \le n^{-1/2}$  then the upper and lower bound on the minimax risk match exactly, which is equal to $\sigma_n\left(\log n \right)^{-1/\beta}$ in the former case and $n^{-1/2}$ in the latter case.  
\item For the intermediate values of $n^{-1/2} \le \sigma_n \le n^{-1/2+\eta}$,  then although the lower and upper bounds do not match completely, they match up to a log factor (that is negligible as compared to the rate close to $n^{-1/2}$) as described below: 
\begin{enumerate}
\item If $\sigma_n\in(n^{-1/2},n^{-1/2}(\log n)^\eta],$ then our upper bound for the minimax risk is $n^{-1/2}(\log\log n)^{-1/\beta}(\log n)^{\eta}$, whereas the  lower bound  is $n^{-1/2}$ (assuming that $\eta<1/\beta$). Therefore the bounds match up to the factor $(\log n)^\eta/(\log \log n)^{1/\beta}$.
\item If $\sigma_n\in(n^{-1/2}(\log n)^\eta,n^{-1/2}(\log n)^{1/\beta}],$ then the upper bound is of the order of $\sigma_n(\log\log n)^{-1/\beta}$ and the lower bound is $n^{-1/2}$ (assuming again that $\eta<1/\beta$), which implies that the rates match up to the maximal factor $(\log n/\log\log n)^{1/\beta}$. 
\item If $\sigma_n\in(n^{-1/2}(\log n)^{1/\beta},n^{-1/2+\eta}],$ then the upper bound is at most $\sigma_n(\log\log n)^{-1/\beta}$ and the lower bound is $\sigma_n(\log n)^{-1/\beta}$, i.e., the upper and lower bound match up to a factor of $(\log n/\log\log n)^{1/\beta}$.
\end{enumerate}
\end{enumerate}
The minimax rates obtained above are summarized in Table \ref{tab: minimaxdeconvolution}. To summarize, we show that our minimax upper bound and lower bound match exactly for the wide ranges of $\sigma_n$, i.e., when $\sigma_n\leq n^{-1/2}$ or $\sigma_n\geq n^{-1/2+\eta}$ for an arbitrarily small $\eta>0$. For the small intermediate region, the bounds match (and are equal to $n^{-1/2}$ ) up to a $\log$ factor.

\begin{table}
\begin{center}
\begin{tabular}{l|l|l}
\hline
range of $\sigma_n$&minimax rate& up to factor\\
\hline
$(0,n^{-1/2}]$&$n^{-1/2}$&1\\
$(n^{-1/2}, n^{-1/2}(\log n)^{\eta}]$&$n^{-1/2}$&$(\log n)^\eta/(\log\log n)^{1/\beta}$\\
$(n^{-1/2}(\log n)^\eta,n^{-1/2}(\log n)^{1/\beta}]$&$n^{-1/2}$&$(\log n/\log\log n)^{1/\beta}$\\
$(n^{-1/2}(\log n)^{1/\beta},n^{-1/2+\eta}]$&$\sigma_n\left(\log n\right)^{-1/\beta}$&$(\log n/\log\log n)^{1/\beta}$\\
$(n^{-1/2+\eta},A]$&$\sigma_n\left(\log n\right)^{-1/\beta}$&1\\
\hline
\end{tabular}
\end{center}
\caption{Minimax rates in the deconvolution problem and the unlinked regression with arbitrarily small $\eta\in(0,1/\beta)$ and $A$ an upper bound for $\sigma_n$. }
\label{tab: minimaxdeconvolution}
\end{table}

%

%
%

\section{Unlinked regression under vanishing noise} 
\label{sec:unlinked}

This section presents a theoretical analysis for estimating the monotone link function in the presence of unlinked data. Recall that, under the setup of the unlinked regression, the observations are assumed to be generated from the model $Y = m_0(X) + \sigma_n \delta$ where $m_0$ is a left-continuous non-decreasing function, the error $\delta \sim \mu_\delta$ for some known distribution function $\mu_\delta$, $\sigma_n$ is a known constant that is typically assumed to be less than or equal to $1$ (may approaches to $0$ as $n \uparrow \infty$ or remains bounded/constant), and $X \sim \mu_X$ for some unknown distribution $\mu_X$. 
As the error scale is absorbed in $\sigma_n$, $\delta$ is assumed to be centered, with variance 1, and independent of $X\in[0,1]$.
We observe $\cX_n= \{X_1, \dots, X_{n}\}$, a set of i.i.d. copies of $X$ and $\cY_n= \{Y_1, \dots, Y_{n}\}$, a set of i.i.d. copies of $Y$. 
Here, the data $\cX_n$ and $\cY_n$ can potentially be collected from different samples (i.e., they can be independent), from the same sample with an unknown link or anything in between (i.e., there is a non-trivial intersection). The two samples could be of different sizes. Still, since we aim to compare the unlinked regression model with the shuffled regression model, in which both samples are necessarily of the same size, we restrict ourselves to the case of samples of the same size. 

Our goal here is to estimate $m_0$ based on  our observations $\cX_n$ and $\cY_n$ and our knowledge about $(\mu_\delta, \sigma_n)$. 
Although we assume that  $X$ is supported on $[0,1]$, our proofs generally apply to any compactly supported $X$. For all positive $M$ and $a$, define $\cM(M, a)$ as the set of non-decreasing functions $m$ on $[0,1]$ that are left-continuous with right-hand limits at every point (LCRR) and satisfy: 
\begin{equation}
\label{eq:def_M}
\int \vert m(x)\vert^{a+2}d\mu_X(x)\leq M \,. 
\end{equation}
Henceforth we assume $m_0 \in \cM(M, a)$ for some $(M, a)$. Therefore, our parameter space is general, incorporating non-smooth and unbounded non-decreasing functions. One can also identify $\cM(M, a)$ with the set of all quantile functions on $\reals$ with certain moment conditions.  The minimax risk of estimating $m_0$ is defined as: 
\begin{equation}
\label{eq: defRU}
\mathcal R_U(M,a,n,\mu):=
\inf_{\widehat m\in\mathcal M}\sup_{m\in\mathcal M(M,a)}\sup_{\pi\in \Pi(m,n)}E_{\pi}\left[\int|\widehat m-m|d\mu \right]
\end{equation}
where $\Pi(m,n)$ is the set of distributions of observations $\cX_n\cup\cY_n$ with $m_0=m$, and $\mu$ can be either the population distribution $\mu_X$ of $X$ (in which case the minimax risk is called  {\it the population minimax risk}) or the empirical distribution $\hat \mu_X$ of the observations $\cX_n$ (in which case the minimax risk is called {\it the empirical minimax risk}).

Here the infimum is taken over all monotone estimators of $m$, i.e. $\widehat m \in \cM$ where $\cM$ is the set of all LCRR functions on $[0,1]$ (more specifically, set of all quantile functions). 
Similar to the preceding section, our primary focus is on quantifying the specific impact of $\sigma_n$ on the minimax risk of estimation, particularly when $\sigma_n$ approaches zero as $n$ increases. 
We start with a lower bound on the minimax risk (as defined in \eqref{eq: defRU}) on the unlinked regression model in terms of $(\sigma_n, n)$, presented in the theorem below: 
\begin{theorem}
\label{theo: DMlowerUL}

Let $(\sigma_n)$ be a bounded sequence of positive numbers, $a>0$ and $M>0$.
Assume that the characteristic function of the error $\delta$ satisfies Assumption \ref{assm:supersmooth},
 that $E\vert 2\delta\vert^{2+a}\leq M$, and that $X$ has a continuous distribution. Then we have: 
$$
\cR_U(M,a,n,\mu) \ge C\left(\sigma_n\left(\log{n}\right)^{-1/\beta} + n^{-1/2}\right) \,,
$$
for all $\sigma_n > 0$, both for $\mu = \mu_X$ or $\mu = \hat \mu_X$. 
\end{theorem}

It is instructive to compare the rate obtained in the above theorem with that of Theorem \ref{theo: DMlower}; in particular, we obtain the same lower bound as in the deconvolution problem. This indicates that the unlinked regression problem is, at most, as hard as the deconvolution problem. If this lower bound is sharp, then the minimax risk of unlinked regression and deconvolution share the same scaling law concerning $(n, \sigma_n)$. That these two problems are related has previously been alluded to in \cite{carpentier2016learning} and \cite{balabdaoui2021unlinked}, however, without rigorous justification. The critical difference between these two problems is the observed $\cX_n$, but as it seems, it does not contribute to the rate of estimation of the signal distribution, i.e., the distribution of $Z:=m_0(X)$. Nevertheless, in the unlinked regression model, we are interested in estimating $m_0$, and the loss function is the $L_1$-distance w.r.t. distribution of $X$, and without $\cX_n$ we cannot identify $m_0$.

In the definition of the minimax risk $\cR_U(M,a,n,\mu)$ (equation \eqref{eq: defRU}), we restricted ourselves to monotone estimators of $m_0$, as that is our parameter space.  
In the following corollary, we show that the infimum could be relaxed to contain all possible estimators (including non-monotone estimators of $m_0$) both for $\mu = \mu_X$ (population distribution of $X$) or $\hat \mu_X$ (empirical distribution of $\cX_n$) in \eqref{eq: defRU}, provided that $X$ has a density that is bounded from above and away from zero:

\begin{corollary}
\label{cor: DMlowerUL}
Under the same setting and assumptions pertaining to Theorem \ref{theo: DMlowerUL}, the same minimax lower bound holds 
even if the infimum is extended to all possible estimators (not necessarily monotone) as long as $X$ has a density bounded from above and bounded away from zero on its support. 
\end{corollary}


Now that we have a lower bound for the minimax risk, we turn to an upper bound. For this task, we assume that $X$ has a density bounded from above and away from zero on its support.

\begin{theorem}
\label{theo: DMUpperUL}
Let $(\sigma_n)$ be a sequence of bounded positive numbers. Assume that the error $\delta$ satisfies Assumption  \ref{assm:error_char_deriv}, that $E\vert \delta\vert^{a+2}\leq M$, and that $X$ has a density that is bounded between $c_X$ and $C_X$ on its support where $c_X>0$ is known. Then,
\begin{eqnarray*}
\mathcal R_U(M,a,n,\mu)\leq Lv_{n}
\end{eqnarray*}
for $\mu \in \{\mu_X, \hat \mu_X\}$, where 
$v_n$ is same as defined in Theorem \ref{theo: DFMupper}. 
\end{theorem}

It is evident from the above theorem that 
we obtain the same upper bound as in the deconvolution problem (Theorem \ref{theo: DFMupper}). 
Since the lower bounds are the same in the two problems, the minimax risk of unlinked regression and deconvolution share the same scaling law concerning $(n, \sigma_n)$. Hence, the minimax rates given in Table \ref{tab: minimaxdeconvolution} for the deconvolution problem also apply to the unlinked regression problem. Consequently, both problems share the same inherent difficulty regarding the rate of estimation of the underlying parameter of interest.

\section{Shuffled regression  under vanishing noise} 
\label{sec:shuffled}

In this Section, we present our theoretical analysis for the shuffled regression model. 
Recall that, in shuffled regression,  the i.i.d. pairs $\{(X_i, Y_i)\}_{i = 1}^n$ are generated from the model $Y_i = m_0(X_i) + \sigma_n \delta_i$ where $\delta_i$ is centered, has variance 1, and is independent of $X_i\in[0,1]$. The distribution of $\delta$ and the standard-error $\sigma_n,$ are assumed to be known.
We only observe the samples $\cX_n= \{X_1, \dots, X_n\}$ and $\cY_n = \{Y_1, \dots, Y_n\}$, where each $Y_i$ is generated by some $X_j \in \cX_n$, but we do not now which $X_j$, as if the $X$ values got shuffled. Another equivalent assumption is that we observe  $\cY_n $ and $\cX_\ord = \{X_{(1)} \le X_{(2)} \le \dots \le X_{(n)}\}$, i.e. the 
the order statistics of $X_i$'s, and this is what we assume in this section. 
In general, the set of covariables $\cX_n$  could be either random or deterministic, but since we aim to compare the shuffled regression model with the unlinked regression model, in which the covariates have to be random, we restrict ourselves here to the case of random covariates.

The monotone link function $m_0$ is our parameter of interest. 
The retrieval of an unknown monotone link function in shuffled regression has recently been explored in \cite{rigollet2019uncoupled} and \cite{slawski2022permuted}. However, the previous papers investigated the \emph{fixed design} model, i.e., $X_1, \dots, X_n$ are assumed to be fixed, whereas as mentioned above, we focus on the \emph{random design} model. Moreover, contrary to earlier literature, we allow for unbounded regression function $m_0$, and we study the effect of a diminishing variance $\sigma_n$ explicitly in the rate of estimation. Estimation of $m_0$ in the presence of random $X$s, all the more so under vanishing error variance, remains largely unknown. We aim to bridge this gap in this subsection. 

Akin to the unlinked regression model, here also the true link function $m_0$ is assumed to be in $\cM(M, a)$ (see equation \eqref{eq:def_M}) for some positive $a$ and $M$. The minimax risk of estimation of $m_0$ is defined as follows: 
\begin{equation*}
\mathcal R_S(M,a,n,\mu):=\inf_{\widehat m\in\cal M}\sup_{m\in\mathcal M(M,a)}E_{m}\left[\int\vert\widehat m-m\vert d\mu_X\right]
\end{equation*}
where $E_{m}$ denotes the expectation with respect to the joint distribution of $\cY_n$ and $\cX_\ord$ when $m_0 = m$,  $\mu$ can be either the population distribution $\mu_X$ of $X$ (in which case the minimax risk is called  {\it the population minimax risk}) or the empirical distribution $\hat \mu_X$ of the observations $\cX_n$ (in which case the minimax risk is called {\it the empirical minimax risk}),
 and infimum is taken over the set $\mathcal M$ of all non-decreasing functions $m$ on $[0,1]$ that are left-continuous with right-hand limits at every point.

We first compute an upper bound on the minimax risk. We obtain a general bound that is valid under arbitrary noise distribution, which is of the order of $\sigma_n$ for the empirical risk, and possibly a slightly larger order for the population risk. This bound is of particular interest when $\sigma_n$ is small. For larger $\sigma_n$ we obtain a more precise bound under supersmooth noise. The bounds are given in the following theorem.

\noindent
\begin{theorem}
\label{thm:shuffled_random_upper}
Let $(\sigma_n)$ be a sequence of bounded positive numbers. Suppose that $X$ has a density bounded below by a known $c_X>0$ on its support and $m_0$ satisfies \eqref{eq:def_M} for some $(a, M) > 0$.
Then  we have the following upper bound on the empirical minimax risk given $E\delta^2<\infty$: 
$$
\mathcal R_S(M,a,n,\hat\mu_X) \le K\sigma_n
$$
where $K > 0$ does not depend on $n$, and  for the population minimax risk we have: 
$$
\mathcal R_S(M,a,n,\mu_X)  \le K\left(\sigma_n + n^{-\frac{a+1}{a+2}} \log{n}\right)\,.
$$
Finally, assume that the error $\delta$ satisfies Assumption  \ref{assm:error_char_deriv}, that $E\vert \delta\vert^{a+2}\leq M$, and that $X$ has a density that is bounded between $c_X$ and $C_X$ on its support where $c_X>0$ is known. Then, if $\sigma_n\geq n^{-1/2}$ and $\eta>0$ is arbitrarily small, we have
$$
\mathcal R_S(M,a,n,\mu) \le K\min\left\{\sigma_n, v_n\right\}.
$$
where $v_n$ is as defined in \eqref{eq:def_v_n}. 
\end{theorem}

It is instructive to compare the upper bounds of Theorem \ref{thm:shuffled_random_upper} with the bounds for the unlinked monotone regression obtained in Section \ref{sec:unlinked}. To do that, first note that the last inequality in Theorem \ref{thm:shuffled_random_upper} concludes that if $\sigma_n\in(n^{-1/2}, n^{-1/2}(\log\log n)^{-1/\beta}(\log n)^\eta]$, then 
$$
\mathcal R_S(M,a,n,\mu) \le K\sigma_n ,,
$$
and iif $\sigma_n\ge n^{-1/2}(\log\log n)^{-1/\beta}(\log n)^\eta$, then 
$$
\mathcal R_S(M,a,n,\mu) \le Kv_n
$$
Hence, the upper bounds of the unlinked and the shuffled regression matches for large $\sigma_n$. As mentioned in the previous section, it is evident from Theorem \ref{theo: DMlowerUL}, that one can never estimate $m_0$ faster than $n^{-1/2}$ in the unlinked regression setup, even if $\sigma_n = 0$. However, for the shuffled regression, the empirical and population minimax risks decay at a rate faster than $n^{-1/2}$ as soon as $\sigma_n \ll n^{-1/2}$. This reveals that estimation of $m_0$ using shuffled data is easier than unlinked data in the presence of a small noise. 

Now we present our results corresponding to the lower bound on the minimax risk of monotone shuffled regression. In Theorem \ref{thm:shuffled_random_upper}, the upper bound for the population risk consists of two terms: $\sigma_n$ and $n^{-(a+1)/(a+2)} \log{n}$. We first focus on the second term, which is free of $\sigma_n$ which indicates that the upper bound cannot be faster than $n^{-(a+1)/(a+2)} \log{n}$ even if $\sigma_n = 0$. The case $\sigma_n = 0$ deserves special attention in the shuffled regression model. 
In the absence of noise, we can recover the unobserved link between $\cX_\ord$ and $\cY_n$ precisely as $Y_{(i)} = m_0(X_{(i)})$ due to the monotone non-decreasing nature of $m_0$. 
Consequently, we know the precise value of $m_0$ at the design points.
This is the reason why the bound that we have for the empirical minimax risk is zero when $\sigma_n=0$.  However, we cannot recover $m_0$ exactly over the whole interval $[0,1]$ even if $\sigma_n = 0$ and the price we pay in the population minimax risk is precisely the second term, i.e., $n^{-(a+1)/(a+2)} \log{n}$ for the population minimax risk. The following theorem establishes that this price is minimax optimal up to a log factor: 
\begin{theorem}
\label{thm:mlb_shuffled_random_noiseless}
Consider the noiseless scenario, i.e. $Y_i = m_0(X_i)$ for some $m_0 \in \cM(M, a)$. Then we have: 
$$
\mathcal R_S(M,a,n,\mu_X) \ge Cn^{- \frac{a+1}{a+2}}(\log{n})^{-\frac{1+\eps}{a+2}} \,,
$$
for any small $\eps > 0$ and some constant $C > 0$. 
\end{theorem}
Theorem \ref{thm:mlb_shuffled_random_noiseless} underscores the importance of the extra factor $n^{-(a+1)/(a+2)}$ of the upper bound in Theorem \ref{thm:shuffled_random_upper}. 
However, the lower bound does not capture the optimal dependence on $\sigma_n$; furthermore, it is not apparent from Theorem \ref{thm:shuffled_random_upper} how the smoothness of the noise (in terms of the decay of its characteristic function) affects the rate of convergence. 
Obtaining a lower bound for the estimation of the monotone link function in a shuffled regression model is a challenging problem. To the best of our knowledge, it is yet to be addressed in full generality. The only known lower bound in this area is Theorem 3 of \cite{rigollet2019uncoupled}, which is i) for the fixed design model and ii) for bounded regression function $m_0$ (with a known bound for the sup-norm of $m_0$)  and iii) under a fixed noise variance. 
Moreover, extending the proof from \cite{rigollet2019uncoupled} to accommodate the random design setup within their framework is not straightforward.
We take a step towards bridging this gap by obtaining a lower bound on the random design model when the errors are generated from a Gaussian distribution with mean $0$ and variance $\sigma_n$, that is allowed to tend to zero as $n\to\infty$: 
\begin{theorem}
\label{thm:mlb_shuffled_random_noisy}
Let $(\sigma_n)$ be a sequence of bounded positive numbers.
Consider the above setup where $X$ has a continuous distribution and $Y = m_0(X) + \sigma_n \delta$, where we observe $\cY_n, \cX_\ord$. Assume $\delta \sim \cN(0, 1)$. Then, there exist positive constants $c, C$ such that
\begin{eqnarray*}
\mathcal R_S(M,a,n,\mu_X)  &\ge& c\sigma_n\left(\frac{\log{n}}{\sqrt{\log{\log{n}}}}\right)^{-1}
\end{eqnarray*}
provided that $\sigma_n\geq Cn^{-1/2}\log{n}(\log{\log{n}})^{-1/2}.$
\end{theorem} 
A few remarks are in order; the lower bound above reveals that for Gaussian noise with sufficiently large $\sigma_n$, the rate of estimation can never be faster than $\sigma_n\sqrt{\log{\log{n}}}/\log{n}$. Since the standard Gaussian distribution satisfies Assumption \ref{assm:supersmooth} with $\beta=2$, this rate is faster than the upper bound obtained for the minimax rate  (see Theorem \ref{thm:shuffled_random_upper}), which is of the order $\sigma_n(\log{n})^{-1/2}$. 
Therefore, we have a gap between the upper and lower bounds on the minimax risk of the shuffled regression problem, and it is not immediate whether the upper bound or the lower bound on the minimax risk is improvable. 
Although we do not have a definite answer, we strongly believe that with a finer analysis, it is possible to improve the lower bound. In fact, we can refine Theorem \ref{thm:mlb_shuffled_random_noiseless}, so that the minimax lower bound matches with the minimax upper bound (which is also the minimax lower bound of the deconvolution problem, see Table \ref{tab: minimaxdeconvolution}) if the following conjecture is true: 
\begin{conjecture}
\label{conj}
Let $N=(n_1,\dots,n_n)$ be a multinomial vector with paramter $n$ and probabilities $(n^{-1},\dots,n^{-1})$. Then for arbitrary $C>0$ one can find $c>0$ such that
\begin{equation}
\label{eq:lemconj}
\lim\inf_{n\to\infty}\bbE\left[\prod_{j\ s.t.\ n_j>0} \left(1-C e^{-c \frac{\log n}{n_j}} \right)^2\right] >0
\end{equation}
Here, both $(C, c)$ are constants independent of $n$. 
\end{conjecture}
\noindent
\begin{theorem}
\label{thm:mlb_shuffled_random_noisy2}
Consider the above setup where $X$ has a continuous distribution and $Y = m_0(X) + \sigma_n \delta$, where we observe $\cY_n, \cX_\ord$. Assume $\delta \sim \cN(0, 1)$. If Conjecture \ref{conj} is true, then there exist positive constants $c,C$ such that
$$
\mathcal R_S(M,a,n,\mu) \ge c\sigma_n\left(\log n\right)^{-1/2}
$$
for $\mu\in\{\mu_X,\hat\mu_X\}$, provided that $\sigma_n \ge Cn^{-1/2}\left(\log n\right)^{1/2}$.
\end{theorem} 
Simulation results support the conjecture, see Section \ref{sec:simulation} for the details. 
Therefore, we believe that for all sufficiently large $\sigma_n$, the optimal rate of estimation for the shuffled monotone regression should coincide with that of the standard deconvolution problem, and all three problems (deconvolution, unlinked regression, and shuffled regression) are asymptotically equivalent in the presence of large noise. Our conclusion is summarized in Table  \ref{tab: minimaxdeconvolutionsuffled}.

\begin{table}
\begin{center}
\begin{tabular}{ll|l|l}
\hline
\multicolumn{2}{l|}{range of $\sigma_n$}&$(0,n^{-1/2}]$ &$(n^{-1/2+\eta},A]$\\
\hline
\multicolumn{2}{l|}{Empirical rate:}&&\\
&Deconvolution&$n^{-1/2}$&$\sigma_n\left(\log n\right)^{-1/2}$\\
&Unlinked&$n^{-1/2}$&$\sigma_n\left(\log n\right)^{-1/2}$\\
&Shuffled&$\leq \sigma_n$&$\sigma_n\left(\log n\right)^{-1/2}$\\
\hline
\multicolumn{2}{l|}{Population rate:}&&\\
&Deconvolution&$n^{-1/2}$&$\sigma_n\left(\log n\right)^{-1/2}$\\
&Unlinked&$n^{-1/2}$&$\sigma_n\left(\log n\right)^{-1/2}$\\
&Shuffled&$\leq \sigma_n  + n^{-\frac{a+1}{a+2}} \log{n}$&$\sigma_n\left(\log n\right)^{-1/2}$\\
\hline
\end{tabular}
\end{center}
\caption{Minimax rates under Gaussian noise with arbitrarily small $\eta\in(0,1/2)$ and $A$ an upper bound for $\sigma_n$. The rates in Shuffled regression are based on Conjecture \ref{conj}.}
\label{tab: minimaxdeconvolutionsuffled}
\end{table}


\section{Simulation for Conjecture \ref{conj}}
\label{sec:simulation}
\begin{figure}
\centering
\begin{subfigure}[h]{0.49\textwidth}
     \centering
     \includegraphics[width=\textwidth]{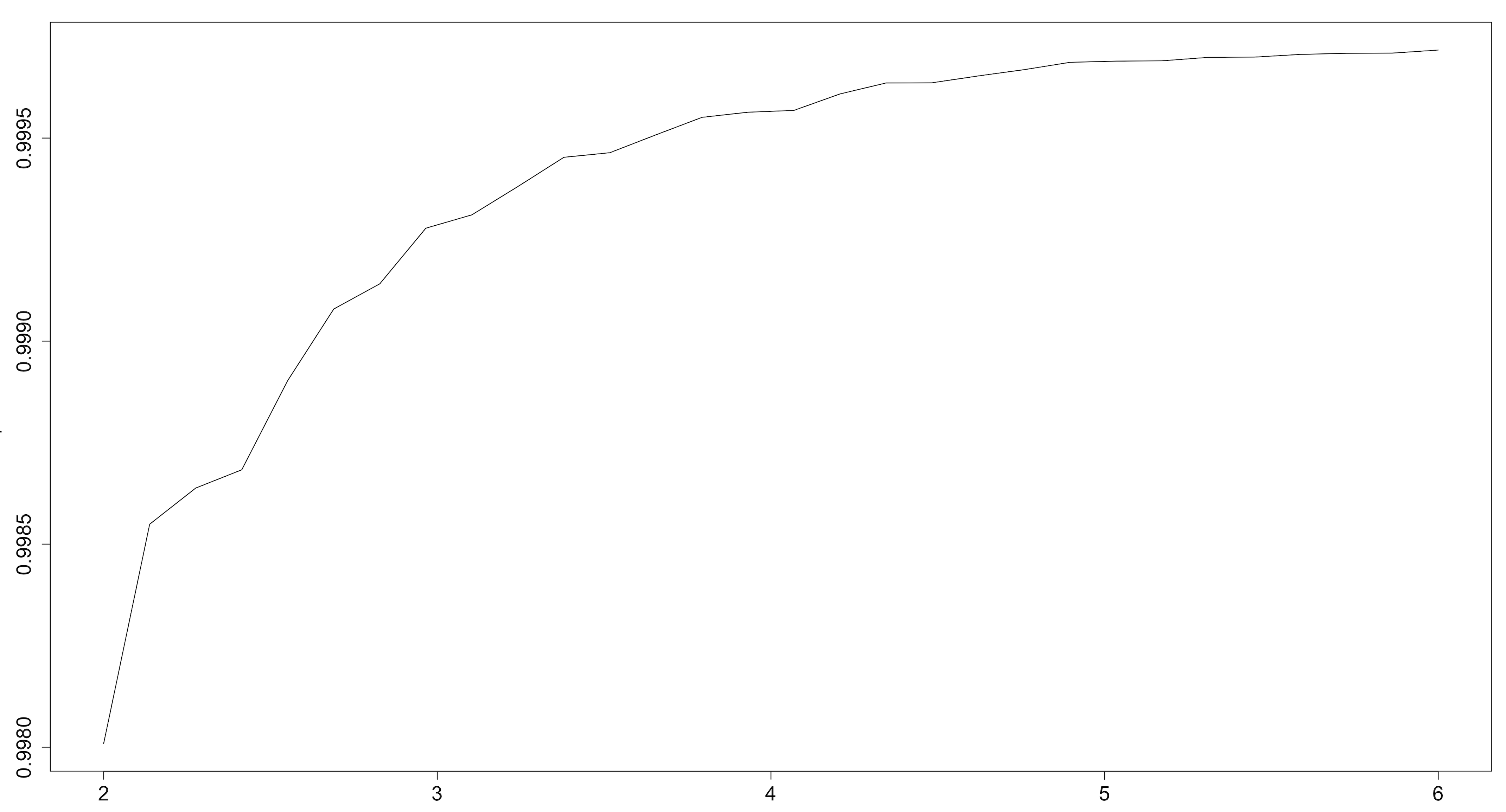}
     \caption{$C = 1, c = 20$}
     \label{fig2:C_500}
 \end{subfigure}
 \hfill
 \begin{subfigure}[h]{0.49\textwidth}
     \centering
     \includegraphics[width=\textwidth]{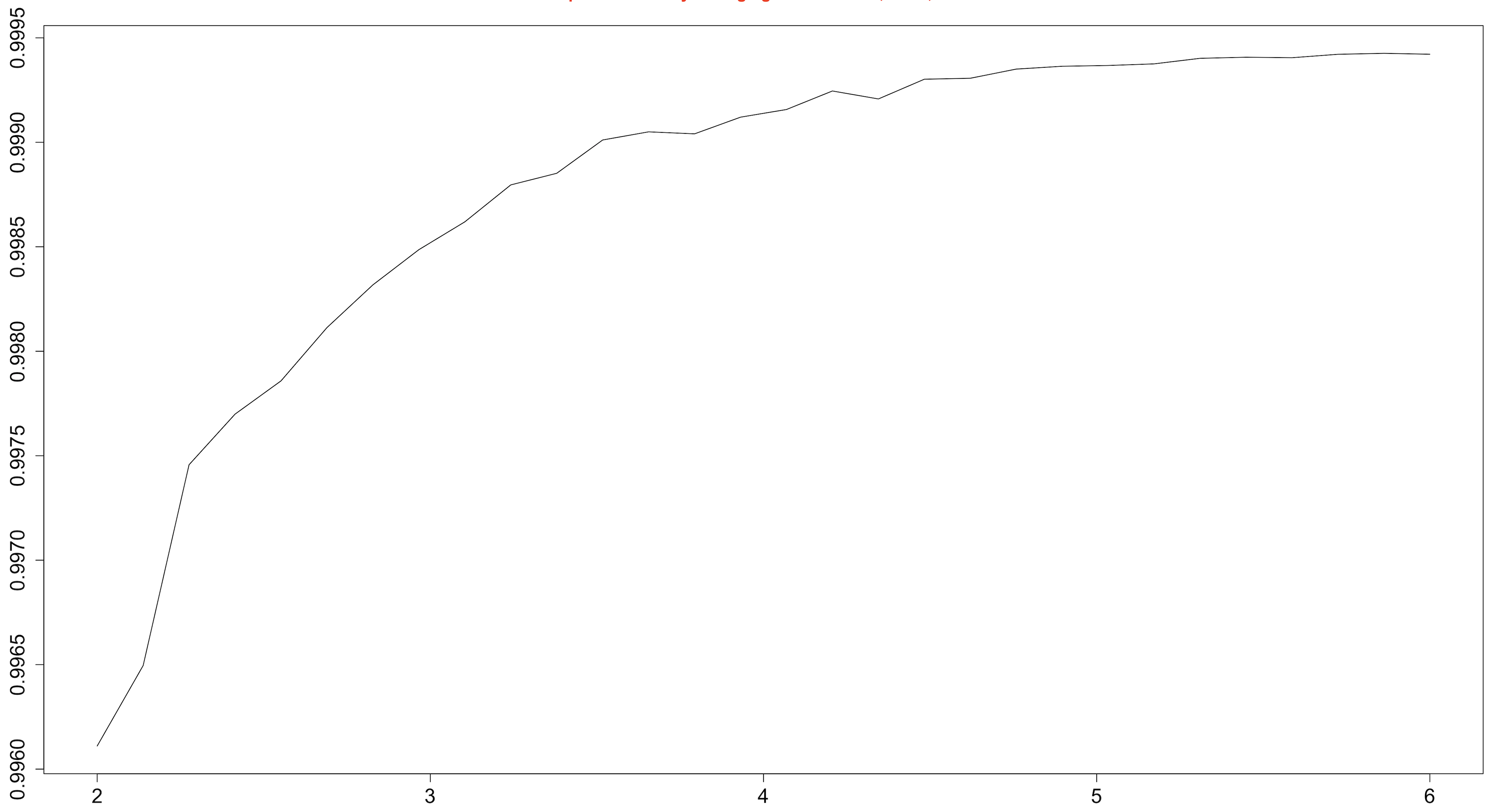}
     \caption{$C = 2, c = 20$}
     \label{fig2:C_1000}
 \end{subfigure}
 \hfill
 \begin{subfigure}[h]{0.49\textwidth}
     \centering
     \includegraphics[width=\textwidth]{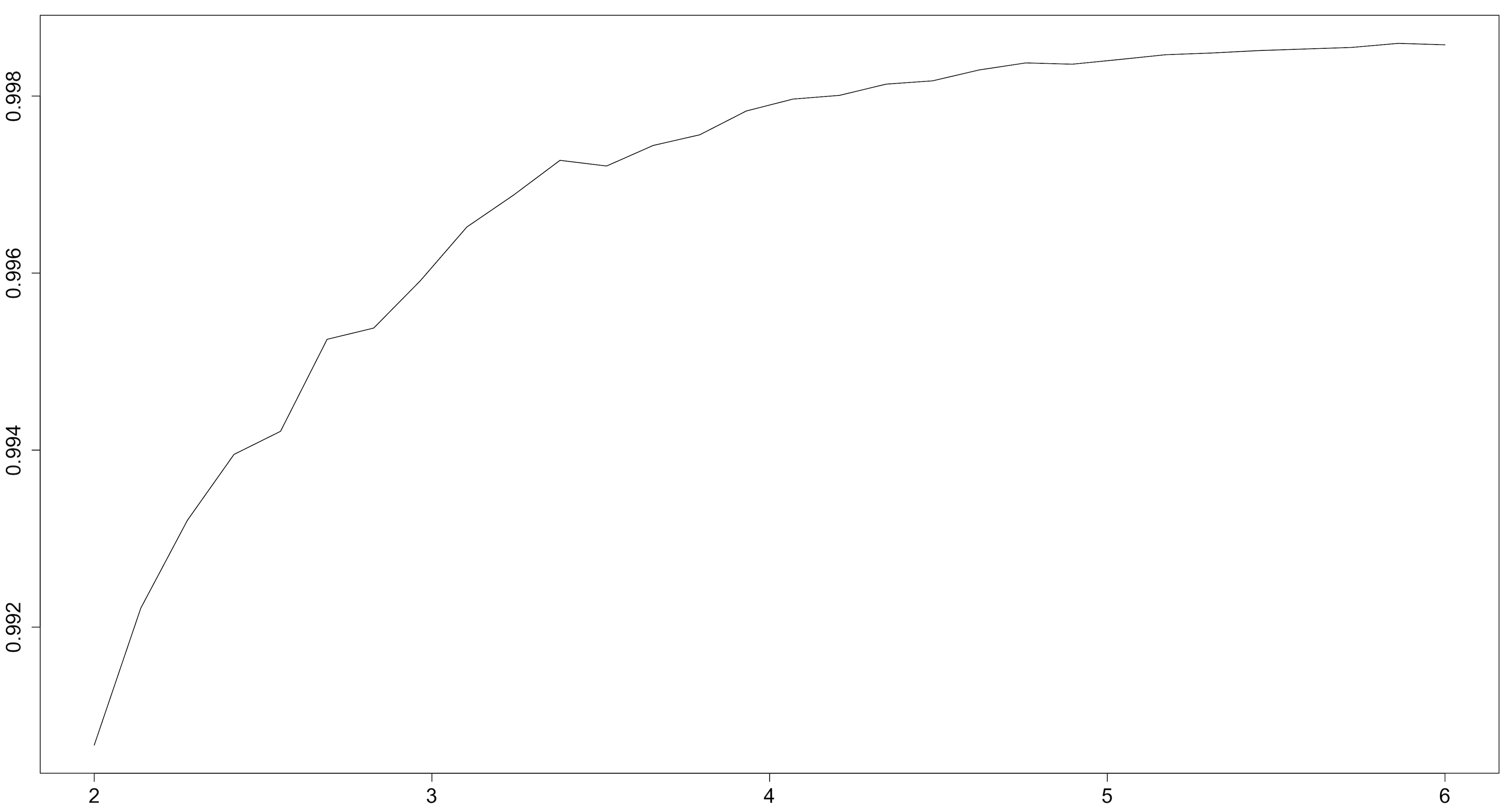}
     \caption{$C = 5, c = 20$}
     \label{fig2:C_2000}
 \end{subfigure}
 \hfill
 \begin{subfigure}[h]{0.49\textwidth}
     \centering
     \includegraphics[width=\textwidth]{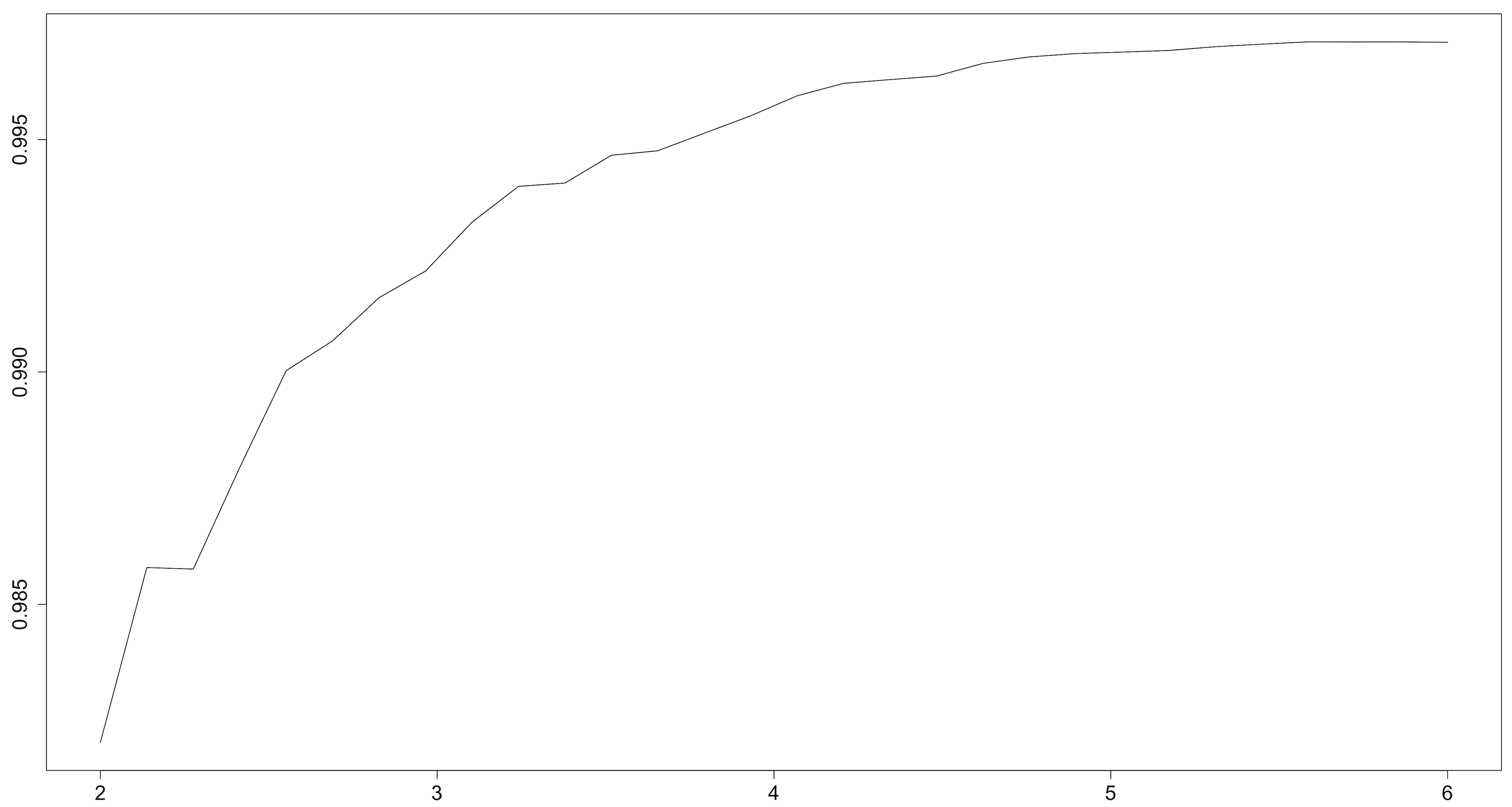}
     \caption{$C = 10, c = 20$}
     \label{fig2:C_5000}
 \end{subfigure}
\begin{subfigure}[h]{0.49\textwidth}
     \centering
     \includegraphics[width=\textwidth]{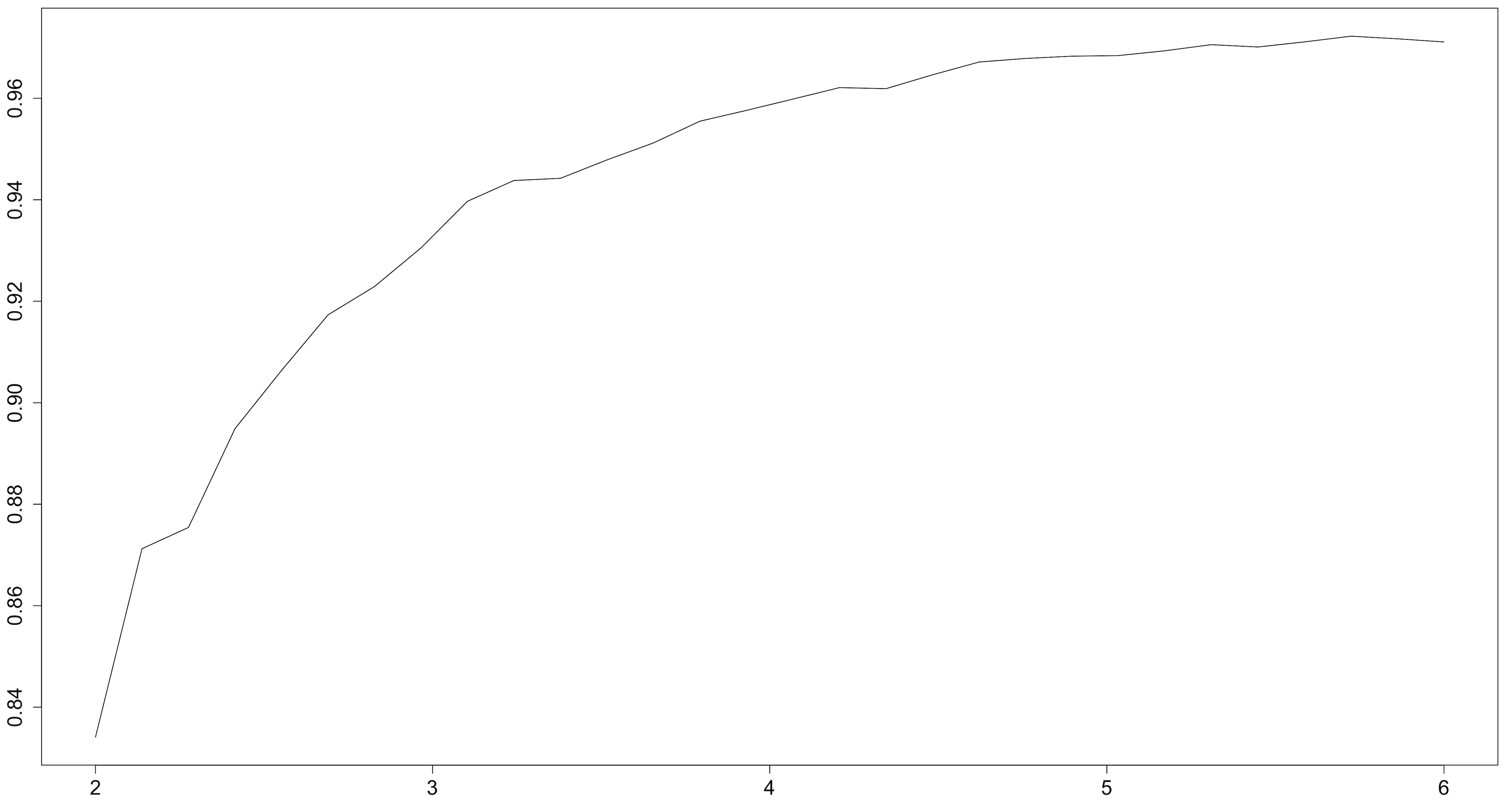}
     \caption{$C = 100, c = 20$}
     \label{fig2:C_500}
 \end{subfigure}
 \hfill
 \begin{subfigure}[h]{0.49\textwidth}
     \centering
     \includegraphics[width=\textwidth]{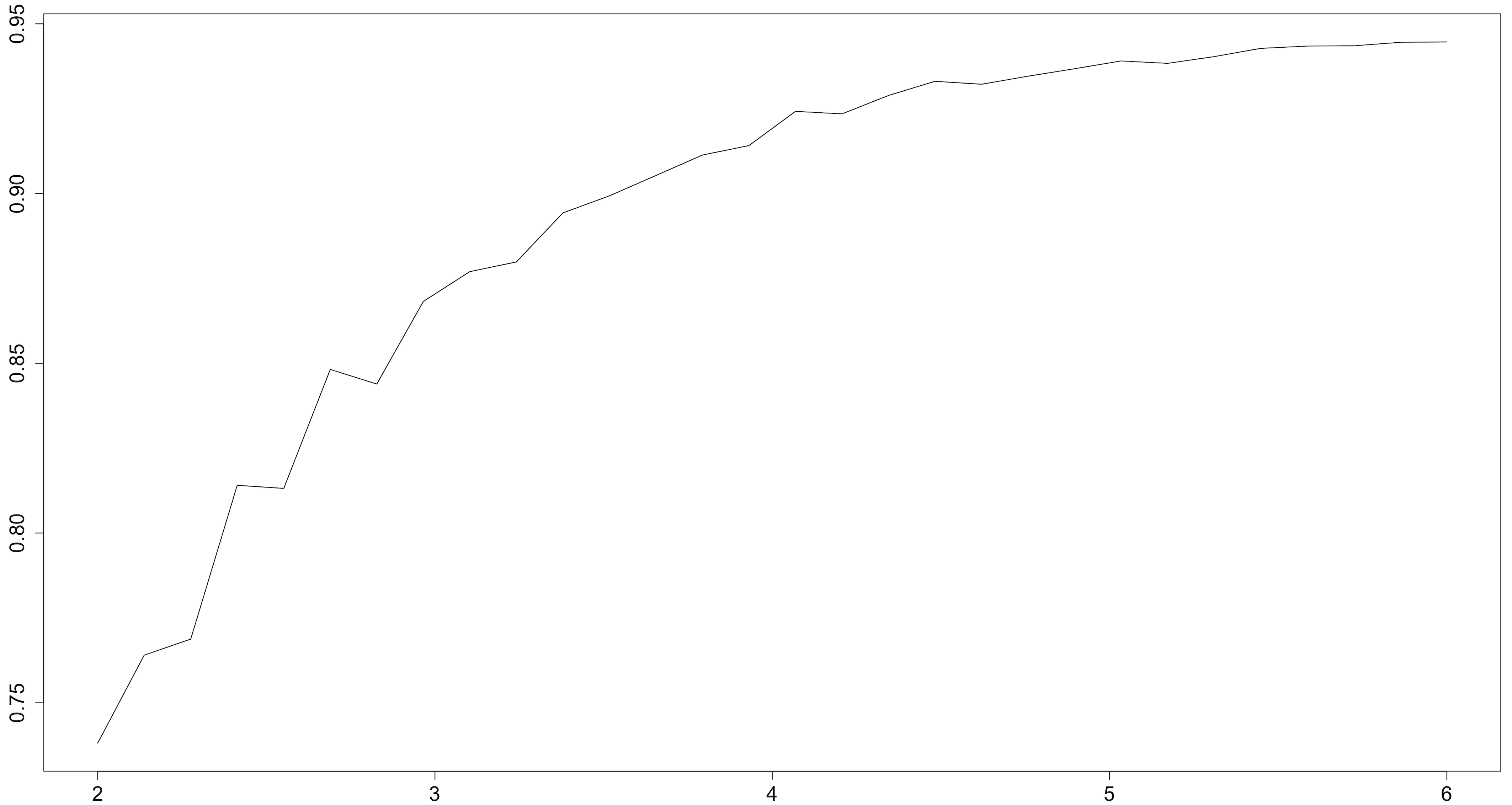}
     \caption{$C = 200, c = 20$}
     \label{fig2:C_1000}
 \end{subfigure}
 \hfill
 \begin{subfigure}[h]{0.49\textwidth}
     \centering
     \includegraphics[width=\textwidth]{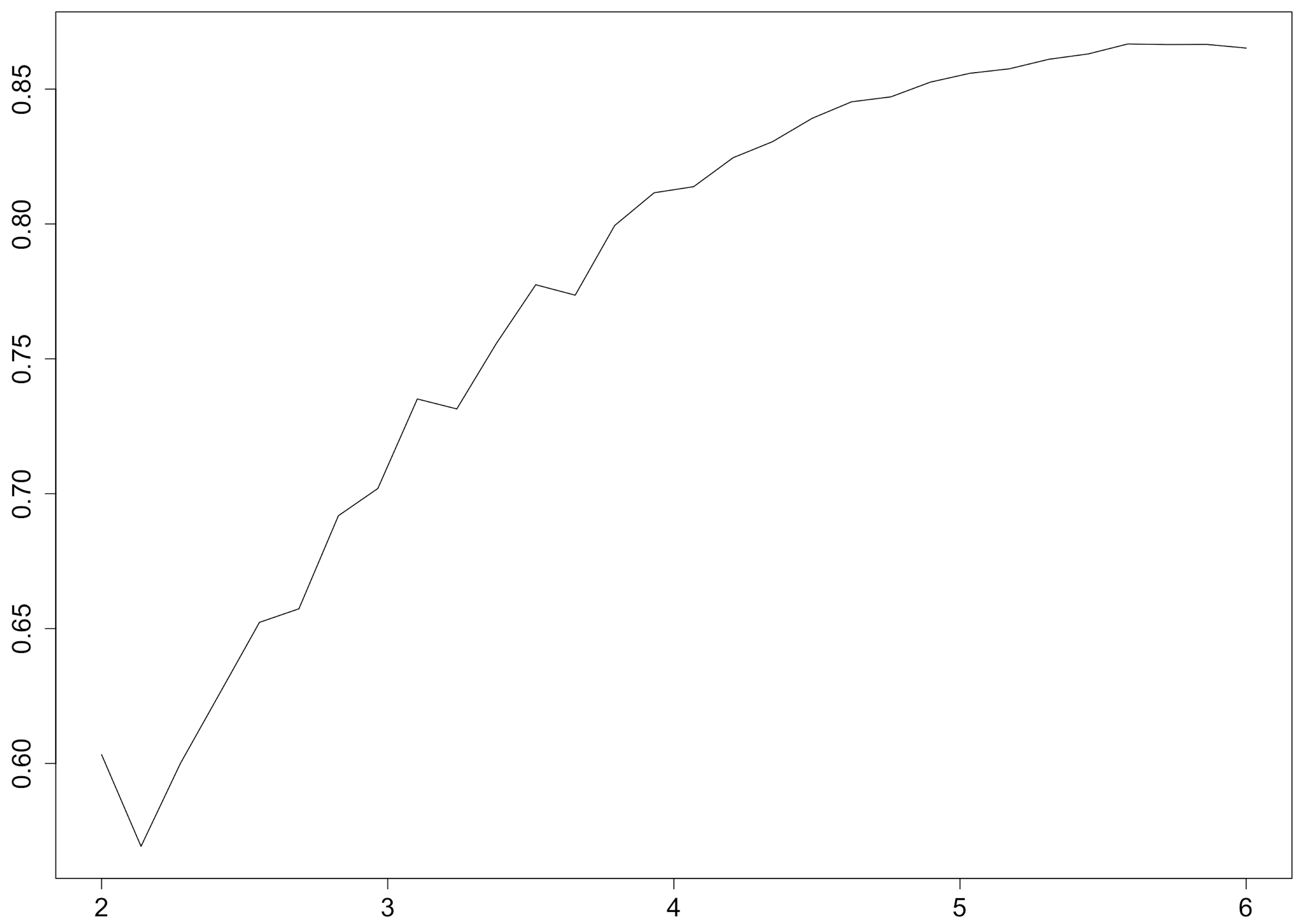}
     \caption{$C = 500, c = 20$}
     \label{fig2:C_2000}
 \end{subfigure}
 \hfill
 \begin{subfigure}[h]{0.49\textwidth}
     \centering
     \includegraphics[width=\textwidth]{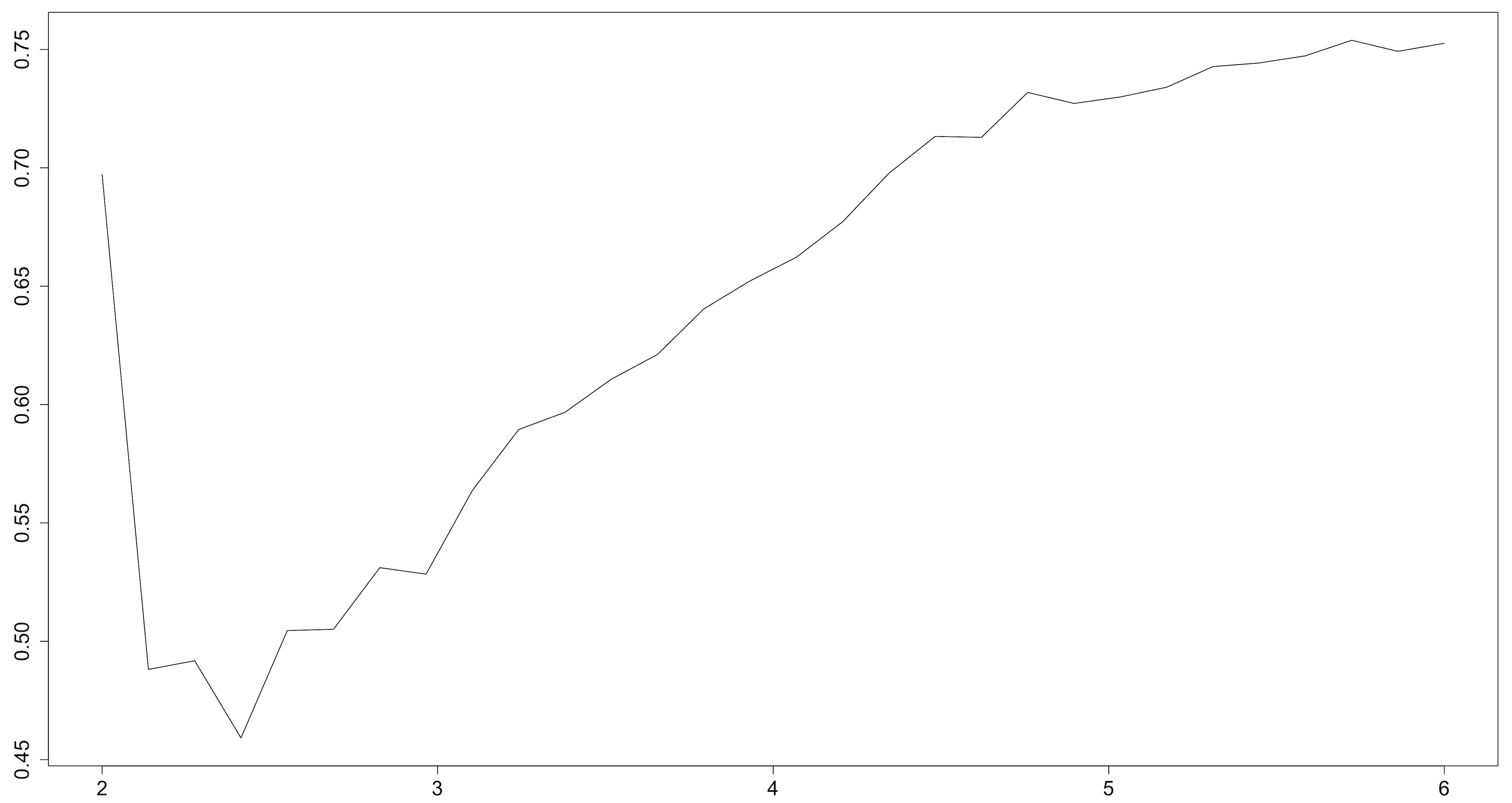}
     \caption{$C = 1000, c = 20$}
     \label{fig2:C_5000}
 \end{subfigure}
    \caption{The eight figures correspond to the approximated value of the expectation in Conjecture \ref{conj} for eight different values of $C$ keeping $c = 20$. The $X$-axis represents the value of $n$ (in $\log$ scale with base 10), and the $Y$-axis denotes the value of expectation approximated by averaging 500 monte-carlo iterations.}
    \label{fig:fixed_c}    
\end{figure}

In this section, we present some simulation results in support of our conjecture (Conjecture \ref{conj}). 
As the power $e$ depends on $\log{n}$ and $n_j$, and both of them are typical of very small order as compared to $n$, it is hard to believe that the right pattern would emerge out of the simulation studies because a potentially huge $n$ is necessary to change the expression in Conjecture \ref{conj} significantly. However, in our small-scale simulation, the results are affirmative. 
 This simulation aims to precisely understand the choice of $c$ depending on $C$. We have run experiments for $8$ different values of $C$, i.e., $C \in \{1, 2, 5, 10, 100, 200, 500, 1000\}$. We show that for each value of $C$, there exists $c$ for which the conjecture seems to be true by our simulation. We especially observe that $c = 20$ works for all these values of $C$. The expectation in Conjecture \ref{conj} is approximated by a Monte Carlo average over 500 iterations. 
Our simulation setup is as follows: 

\begin{enumerate}
\item We vary $n$ between $10^2$ to $10^6$ and we take $30$ different values of $n$ within this region. The values are equally spaced in the $\log_{10}$ scale. 
\item For each of the values of $n$, and a particular choice of $(C, c)$, we approximate the expectation in \eqref{eq:lemconj} by taking average over $500$ monte-carlo iterations. In other words, we generate a multinomial random variable from $\mult(n; n^{-1}, n^{-1}, \dots, n^{-1})$ $500$ times, calculate the term inside the expectation of \eqref{eq:lemconj} and take the average over these $500$ terms. 
\item We set $c = 20$ and vary $C \in \{1, 2, 5, 10, 100, 200, 500, 1000\}$. 
\end{enumerate}
We present the eight plots corresponding to eight different values of $C$ in Figure \ref{fig:fixed_c}. As the figure shows, the monte carlo approximation of the expectation remains bounded away from $0$. 
The main takeaway from our simulation study is that it seems plausible that the conjecture is true if one can choose $c$ judiciously based on $C$. However, 
a proper mathematical justification is necessary to reach a conclusion, which we leave as a future research problem.  

\section{Conclusion and future work}
\label{sec:conclusion}
This manuscript investigates the relationship between the monotone shuffled regression model and the monotone unlinked regression model, explicitly considering the vanishing variance condition where the noise variance approaches $0$ as $n \rightarrow \infty$. 
Our findings reveal a phase transition; if $\sigma_n$ falls below a certain threshold close to $n^{-1/2}$, then the minmax risk of the shuffled regression model (in terms of both $n$ and $\sigma_n$) decays faster than that of the unlinked regression model. 
However, when $\sigma_n$ is larger than the threshold, both problems exhibit similar minimax risk, which aligns with the minimax risk of the Wasserstein deconvolution problem. The main findings for Gaussian noise are summarized in Table \ref{tab: minimaxdeconvolutionsuffled}.

Nevertheless, several unresolved questions remain, serving as avenues for future investigations: 
\begin{enumerate}
\item Our analysis relies on supersmooth error, as defined in Definition \ref{def:supersmooth}. It would be intriguing to explore whether similar phase transitions also occur with ordinary smooth errors -- i.e., error distributions characterized by heavier-tailed characteristic functions.

\item As the primary objective of this paper is to compare the minimax risk between unlinked and shuffled regression, we adopt a random design model, where the unobserved signal values are treated as random -- a prerequisite for the unlinked regression model. In contrast, shuffled regression does not necessitate random signals. Notably, recent research by \cite{rigollet2019uncoupled,slawski2022permuted} explored monotone shuffled regression within a fixed design framework, assuming signals to be fixed unobserved points, albeit with $\sigma_n \asymp 1$. 
Therefore, it is intriguing to investigate the minimax risk of shuffled regression in a fixed design scenario with a diminishing noise variance and possibly unbounded regression function.

\item In this paper, we demonstrate that the minimax risk associated with shuffled regression, unlinked regression, and the deconvolution problem is of the same order when $\sigma_n$ is not too small (of the order much larger than $n^{-1/2}$). This implies a comparable difficulty in estimating the signal distribution for all three problems. However, achieving a more \emph{accurate} result would involve establishing asymptotic equivalence among these three experiments in the sense of Le-Cam (\cite{le2012asymptotic}). This endeavor calls for an entirely new set of tools and is deferred as an intriguing avenue for future exploration. 

\item Finally, our findings are contingent on the univariate signal. Extending these results for the multivariate signals would be a compelling direction for future exploration. 
\end{enumerate} 
{\bf Acknowledgement: }We thank Prof. Moulinath Banerjee, Prof. Ya'acov Ritov, Prof. J\'er\^ome Dedecker, and Prof. Phillip Rigollet for their insightful comments. This research has been conducted within the FP2M federation (CNRS FR 2036).

\bibliographystyle{plain}
\bibliography{shuffreg2}
 \begin{appendix}
\section{Proof of the main results}
\label{app:main}

\subsection{Notations and basic facts}\label{sec: notation}
Let $X$ be a random variable supported on $[0,1]$ with a continuous distribution. 
For all $m\in\cal M$, we denote by $m^{-1}$  the generalized inverse defined for $x\in\R$ by 
\begin{equation}\label{eq: inv}
m^{-1}(x)=\sup\{t\in[0,1]\mbox{ such that }m(t)\leq x\}
\end{equation}
with the convention that the supremum of an empty set is 0. By left-continuity of $m$ we then have $m\circ m^{-1}(x)\leq x$ for all $x\geq m(0)$, from which we derive that for all $t\neq0$ and all $x\in\R$, we have the equivalence
\begin{equation*}
m(t)\leq x \Longleftrightarrow t\leq m^{-1}(x).
\end{equation*}
Hence, if we 
define $\mu_m$ as the distribution of the random variable $m(X)$,  for all $m\in{\mathcal M}$, its distribution function is $F_X\circ m^{-1}$ where $F_X$ is the distribution function of $X$. (We used here that $X\neq 0$ almost surely.)
Moreover,
\begin{eqnarray}\label{eq: mum}
\int \vert x\vert ^{a+2}d\mu_m(x)
&=&
\int \vert m(x)\vert ^{a+2} d\mu_X(x)
\end{eqnarray}
so $m$ belongs to $\mathcal M(M,a)$ defined in Section \ref{sec:unlinked}  if and only if the corresponding $\mu_m$ belongs to $\mathcal D(M,a)$ defined in Section \ref{sec: deconvolutionSM}. 

We denote by $\hat\mu_X$ the empirical distribution of the $X$-sample $(X_1,\dots,X_{n})$ and by $\widehat F_X$ and $\widehat F_X^{-1}$ the corresponding distribution and quantile functions.

For all $m\in\cal M$, $\hat\mu_{m}$ denotes the distribution that puts mass $n^{-1}$ at each $m(X_i$), $i\in\{1,\dots,n\}$. Using the above properties of the generalized inverse, one obtains that  the corresponding distribution function is $\widehat F_X\circ m^{-1}$ and the quantile function is $m\circ \widehat F_X^{-1}$.

We recall that from \cite[Theorem 2.10 and Theorem 2.9]{bobkov2019one},  for all probability measures $\mu$ and $\mu'$ on $\R$ with finite first moment and corresponding distribution functions $F$ and $G$, and quantile functions $F^{-1}$ and $G^{-1}$ one has
\begin{equation}\label{eq: W1inv}
W_1(\mu,\mu')=\int_{0}^1\vert F^{-1}(x)-G^{-1}(x)\vert dx
\end{equation}
and
\begin{equation}\label{eq: W1dir}
W_1(\mu,\mu')=\int_{-\infty}^{+\infty}\vert F(x)-G(x)\vert dx.
\end{equation}

Now, for arbitrary $m _0,m_1\in\cal M$ and $U$ a uniform variable on $[0,1]$ we have
\begin{eqnarray}\label{eq: L1toW1}\notag
\int_0^1\vert m_0(x)-m_1(x)\vert d\mu_X(x)
&=&
E\vert m_0(X)-m_1(X)\vert\\ \notag
&=&E\vert m_0\circ F_X^{-1}(U)-m_1\circ F_X^{-1}(U)\vert\\
&=& W_1(\mu_{m_0},\mu_{m_1})
\end{eqnarray}
provided that $\mu_{m_0}$ and $\mu_{m_1}$ have a finite first moment, where we used \eqref{eq: W1inv} for the last equality. Similarly,
\begin{eqnarray}\label{eq: L1toW1emp}
\int_0^1\vert m_0(x)-m_1(x)\vert d\hat\mu_X(x)
&=& W_1(\hat\mu_{m_0},\hat\mu_{m_1}).
\end{eqnarray}
Note that $\hat \mu_m$ has finite first moment for any $m\in\cal M$ as a discrete probability measure on a finite support.

\subsection{Entropy results}
We recall that the $\delta$-entropy with bracketing of a class of functions $\mathcal M$ with respect to a  norm $\|\,.\,\|$ is denoted by
\begin{equation*}
\log N_{[]}(\delta,\mathcal M,\|\,.\,\|)
\end{equation*}
and is defined as the logarithm of the minimum number of brackets $[l,u]$ with $\|l-u\|\leq\delta$ needed to cover $\mathcal M$. Here, $u$ and $l$ are functions that do not need to belong to $\mathcal M$ and the bracket $[l,u]$ is the set of all functions $f$ with $l\leq f\leq u$.

\begin{lemma}\label{lem: entropy}
Suppose that $X$ has a continuous distribution.
Then, for all positive $a$ and $M$, there exists $L>0$  such that for all $\delta>0$, one has
\begin{equation*}
\log N_{[]}(\delta,\mathcal M(M,a),\|\,.\,\|_2)\leq L\delta^{-1}
\end{equation*}
where for all functions $f$, $\|f\|_2^2=\int |f|^2d\mu_X.$ 
\end{lemma}

\begin{lemma}\label{lem: emptopop} 
Let $\mu_X$ be a continuous distribution on $\R$ and $\hat\mu_X$ be the empirical distribution of a $n$-sample from $\mu_X$. Then, for all positive $a$ and $M$, there exists $C>0$  such that for arbitrary random function  $\widehat m_n\in\mathcal M(M,a)$ and  $m\in\mathcal M(M,a)$ one has
\begin{eqnarray*}
E\left|\int\vert \widehat m_n-m\vert (d\mu_X-d\hat \mu_X)\right|&\leq&
E\left|\int\sup_{f\in \mathcal M(M,a)}\vert f-m\vert (d\mu_X-d\hat \mu_X)\right|\\
&\leq& Cn^{-1/2}.
\end{eqnarray*}
\end{lemma}

\subsection{Definition of a finite familly in $\mathcal D(M,a)$}\label{sec: family}

In this section, we build a finite familly in $\mathcal D(M,a)$ that will be usefull to compute lower bounds for the minimax risk, using the fact that the supremum over $\mathcal D(M,a)$ in the definition of the minimax risk is not smaller than the supremum over the finite familly.
If $E|\delta|^{2+a}$ is finite, it follows from \cite[Lemma 1]{dedecker2013minimax} (where we take $p=d=A=1$ and $b=1+a$) that one can find $\kappa_1\in(0,1)$ and $\kappa_2>1$ such that
\begin{equation}\label{assm:DM}
P(\vert\delta-t\vert\leq\vert t\vert^{\kappa_1})=O(\vert t\vert^{-\kappa_2}) \mbox{ as }\vert t\vert\to\infty
\end{equation}
and 
\begin{equation}\label{eq: DMkappa}
\max\left(2+\frac a2,\frac{\kappa_2}{2\kappa_1}+\frac 12\right)<\kappa_2.
\end{equation} 
Hence, we can consider  $r>0$ such that
\begin{equation}\label{eq: condr}
\max\left(\frac 32+\frac a2,\frac{\kappa_2}{2\kappa_1}\right)<r<\kappa_2-\frac 12
\end{equation}
and define the density function
\begin{equation*}
f_{0,r}(t)=C_r(1+t^2)^{-r}
\end{equation*}
where $C_r$ is a normalizing constant so that $f_{0,r}$ integrates to one. Inspecting the proof of \cite[Lemma 1]{dedecker2013minimax}, one can see that we can choose $\kappa_2=2+a$ and $\kappa_1$ close enough to one so that one can choose $r>0$ such that 
\begin{equation}\label{eq: r}
r\in\left(\frac 32+\frac a2,\frac 32+a\right).
\end{equation}

 For an arbitrary sequence of positive integers $b_n$ (to be fixed latter) 
and $\theta\in\{0,1\}^{b_n}$, let $\mu_\theta$ be the probability measure with density
\begin{equation}\label{eq: ftheta}
f_\theta(t)=\sigma_n^{-1}f_{0,r}(\sigma_n^{-1}t)+\sum_{s=1}^{b_n}\theta_s\sigma_n^{-1}H(b_n(\sigma_n^{-1}t-t_{s,n})),\ t\in\R,
\end{equation}
where $(\sigma_n)_
n$ is a bounded sequence of  positive numbers, $t_{n,s}=(s-1)/b_n$ and $H$ is a bounded function whose integral on the line is 0 and which satisfies conditions (A1), (A2) and (A3) in the proof of Theorem 3 of \cite{dedecker2013minimax}. 
As noticed in \cite{dedecker2013minimax}, following the definition of these conditions, 
there indeed exists $H$ that satisfies the conditions and such that $f_\theta$ is a density with respect to the Lebesgue measure on $\R$; and $\mu_\theta\in\mathcal D(M,a)$ for all $\theta$ provided that $M$ is sufficiently large. (To see this, note that if $U_\theta$ is a random variable with density $f_\theta$ as above, then $U_\theta$ has the same distribution as $\sigma_n V_\theta$ where $V_\theta$ has the density in Formula (9) in \cite{dedecker2013minimax}, so $E|U_\theta|^{a+2}=\sigma_n^{a+2}E|V_\theta|^{a+2}$ where according to Formula (10)  in \cite{dedecker2013minimax},  $E|V_\theta|^{a+2}$ is bounded uniformly in $\theta$.) Note that we can choose $H$ in such a way that $f_\theta(t)>0$ for all $t$ and $\theta$.

\begin{lemma}\label{lem: mutheta} Let $(\sigma_n)_n$ be a bounded sequence of positive numbers, $a>0$ and $M>0$.
Let $r$ satisfies \eqref{eq: r} and $H$ satisfies the conditions above.
 Then, there exists $b_n\sim(\log n)^{1/\beta}$  such that $f_\theta$ in \eqref{eq: ftheta} is the density of a probability measure $\mu_\theta\in\mathcal D(M,a)$ for all $\theta$ provided that $M$ is sufficiently large. If moreover,  Assumption  {\bf (L)} holds and $E|\delta|^{2+a}<\infty,$ then $b_n$ can be choosen in such a way that for some $C>0$,
\begin{eqnarray}\label{eq: mutheta}
\mathbb{E}_{\tilde \theta_n} E_{(\mu_{\tilde \theta_n}\star\mu_\epsilon)^{\otimes n}}W_1(\tilde\mu_n,\mu_{\tilde\theta_n})\geq C\sigma_n(\log n)^{-1/\beta} \,,
\end{eqnarray}
 for all estimators $\tilde\mu_n$,
where $\mathbb{E}$ denotes the expectation with respect to the distribution of $\tilde\theta_n$,  a random vector taking values in $\{0,1\}^{b_n}$ whose components are i.i.d. Bernoulli random variables with parameter of success $1/2$. 
\end{lemma}

\begin{lemma}\label{lem: Ftheta}
Let $\sigma_n$ be a bounded  sequence of positive numbers, $r>1/2$ and $b_n$ a sequence of positive numbers.
Assume that for all $\theta\in\{0,1\}^{b_n}$ the function in \eqref{eq: ftheta} is the density of a probability measure $\mu_\theta$.
Then there exists $C>0$ such that for all $\theta$,
\begin{eqnarray*}
\int W_1(\hat\mu_Z,\mu_\theta)d\mu_\theta^{\otimes n}(Z_1,\dots,Z_n)\leq Cn^{-1/2}
\end{eqnarray*}
where $\hat\mu_Z$ denotes the empirical distribution of a sample $Z_1,\dots,Z_n$ from $\mu_\theta$.
\end{lemma}

\subsection{Proof of Theorem \ref{theo: DMlower}}\label{sec: proofDMlower}
By Lemma \ref{lem: mutheta}, $\mu_\theta\in\mathcal D(M,a)$ for all $\theta\in\{0,1\}^{b_n}$. Hence,  Lemma \ref{lem: mutheta} also proves that for all estimators $\tilde\mu_n$ we have
\begin{equation*}
\sup_{\mu\in\mathcal D(M,a)}E_{(\mu\star\mu_\epsilon)^{\otimes n}}W_1(\mu,\tilde \mu_n)\geq \mathbb{E}_{\tilde \theta_n} E_{(\mu_{\tilde \theta_n}\star\mu_\epsilon)^{\otimes n}}W_1(\tilde\mu_n,\mu_{\tilde\theta_n})\geq C\sigma_n(\log n)^{-1/\beta}.
\end{equation*}
We shall prove below that for arbitrary estimator $\tilde\mu_n$ we have
\begin{equation}\label{eq: paramlower}
\sup_{\mu\in\mathcal D(M,a)}E_{(\mu\star\mu_\epsilon)^{\otimes n}}W_1(\mu,\tilde\mu_n)\geq Cn^{-1/2}.
\end{equation}
Theorem \ref{theo: DMlower} immediately follows from the previous two displays.

We turn to the proof of \eqref{eq: paramlower}. 
From the definition of Wasserstein distance $W_1$, for arbitrary distributions $\mu_1$ and $\mu_2$ we have
\begin{eqnarray*}
W_1(\mu_1\star\mu_\epsilon,\mu_2\star\mu_\epsilon)=\inf_{Z_1,Z_2,\epsilon_1,\epsilon_2}E\vert (Z_1+\epsilon_1)-(Z_2+\epsilon_2)\vert
\end{eqnarray*}
where the infimum is taken over all $\{Z_1,Z_2,\epsilon_1,\epsilon_2\}$ such that $Z_i$ has distribution $\mu_i$ and is independent of $\epsilon_i$ that has distribution $\mu_\epsilon$, for $i=1$ and $i=2$. Restricting the above infimum to cases where $\epsilon_1=\epsilon_2$ yields
\begin{eqnarray*}
W_1(\mu_1\star\mu_\epsilon,\mu_2\star\mu_\epsilon)
&\leq&
 \inf_{Z_1,Z_2,\epsilon}E\vert (Z_1+\epsilon)-(Z_2+\epsilon)\vert\\
&=&
 \inf_{Z_1,Z_2}E\vert Z_1-Z_2\vert\\
&=&
W_1(\mu_1,\mu_2) \,.
\end{eqnarray*}
Therefore,
\begin{eqnarray*}
W_1(\mu_1,\mu_2)\geq W_1(\mu_1\star\mu_\epsilon,\mu_2\star\mu_\epsilon).
\end{eqnarray*}
Hence, 
\begin{eqnarray*}
\inf_{\tilde \mu_n}\sup_{\mu\in\mathcal D(M,a)}E_{(\mu\star\mu_\epsilon)^{\otimes n}}W_1(\mu,\tilde\mu_n)
&\geq&
\inf_{\tilde \mu_n}\sup_{\mu\in\mathcal D(M,a)}E_{(\mu\star\mu_\epsilon)^{\otimes n}}W_1(\mu\star\mu_\epsilon,\tilde\mu_n\star\mu_\epsilon).
\end{eqnarray*}
By convexity of the function $x\mapsto x^b$ where $b=a+2$, for $Y=Z+\epsilon$ we have
\begin{eqnarray*}
E\vert Z\vert ^b
&\leq& 2^{b-1}\left(E\vert Y\vert ^b+E\vert \epsilon\vert ^b\right) \le M \,.
\end{eqnarray*}
 provided that $E\vert 2Y\vert ^b\leq M$ and $E\vert2\epsilon\vert^{b}\leq M$. Hence, $\mu$ belongs to $\mathcal D(M,a)$ whenever $\mu\star\mu_\epsilon$ belongs to $\mathcal D(2^{-b}M,a)$ and therefore, 
\begin{eqnarray*}
\inf_{\tilde \mu_n}\sup_{\mu\in\mathcal D(M,a)}E_{(\mu\star\mu_\epsilon)^{\otimes n}}W_1(\mu,\tilde\mu_n)\geq
\inf_{\tilde \mu_n}\sup_{\mu\star\mu_\epsilon\in\mathcal D(2^{-b}M,a)}E_{(\mu\star\mu_\epsilon)^{\otimes n}}W_1(\mu\star\mu_\epsilon,\tilde\mu_n\star\mu_\epsilon).
\end{eqnarray*}
This implies that
\begin{eqnarray*}
\inf_{\tilde \mu_n}\sup_{\mu\in\mathcal D(M,a)}E_{(\mu\star\mu_\epsilon)^{\otimes n}}W_1(\mu,\tilde\mu_n)\geq
\inf_{\tilde \mu_n}\sup_{\mu\in\mathcal D(2^{-b}M,a)}E_{\mu^{\otimes n}}W_1(\mu,\tilde\mu_n)
\end{eqnarray*} 
where the infimum on the right hand side is now taken over all estimators $\tilde\mu_n$ of the common  distribution of the observations $Y_1,\dots,Y_n$.
Denote by $\mu_0$  the standard Gaussian distribution and by $\mu_1$  the Gaussian distribution with variance 1 and mean $n^{-1/2}$. For simplicity, we assume $M$ large enough  so that both $\mu_0$ and $\mu_1$ belong to $\mathcal D(2^{-b}M,a)$. (For smaller $M$, it suffices to normalize the above Gaussian laws to ensure that they both belong to  $\mathcal D(2^{-b}M,a)$). We then have
\begin{eqnarray*}
\inf_{\tilde \mu_n}\sup_{\mu\in\mathcal D(M,a)}E_{(\mu\star\mu_\epsilon)^{\otimes n}}W_1(\mu,\tilde\mu_n)\geq
\inf_{\tilde \mu_n}\sup_{\mu\in\{\mu_0,\mu_1\}}E_{\mu^{\otimes n}}W_1(\mu,\tilde\mu_n).
\end{eqnarray*}
We will show below that for the  right hand side we have
\begin{eqnarray}\label{eq: Lsmallsigma}
\inf_{\tilde \mu_n}\sup_{\mu\in\{\mu_0,\mu_1\}}E_{\mu^{\otimes n}}W_1(\mu,\tilde\mu_n)\geq Cn^{-1/2}
\end{eqnarray}
for some $C>0$. Equation \eqref{eq: paramlower} then follows from the previous two displays.

For an arbitrary estimator $\tilde\mu_n$, we denote by $\hat \mu_n$ a distribution in $\{\mu_0,\mu_1\}$ that is closest to $\tilde\mu_n$ in the $W_1$-distance (it does not matter how we handle ties).
For any $\mu\in\{\mu_0,\mu_1\}$ we then have
\begin{equation*}
W_1(\hat \mu_n,\mu)\leq W_1(\tilde \mu_n,\mu)+W_1(\tilde\mu_n,\hat\mu_n) \leq 2W_1(\tilde \mu_n,\mu). 
\end{equation*}
Hence,
\begin{eqnarray*}
E_{\mu^{\otimes n}} W_1(\tilde\mu_n,\mu)&\geq& \frac12 E_{\mu^{\otimes n}} W_1(\hat\mu_n,\mu) \\
&\geq&\frac 12 P_{\mu^{\otimes n}} (\hat\mu_n\neq \mu)W_1(\mu_0,\mu_1).
\end{eqnarray*}
Now, with $F_0$ and $F_1$ the respective distribution functions of $\mu_0$ and $\mu_1$ we have 
\begin{equation*}
W_1(\mu_0,\mu_1)=\int_\R\vert F_0(t)-F_1(t)\vert dt=n^{-1/2}
\end{equation*}
where the first equality follows from \cite[Theorem 2.9]{bobkov2019one} whereas the second one follows from simple algebra combined to the  definitions of $\mu_0$ and $\mu_1$. Therefore,
\begin{equation*}
\inf_{\tilde \mu_n}\sup_{\mu\in\{\mu_0,\mu_1\}}E_{\mu^{\otimes n}}W_1(\mu,\tilde\mu_n)\geq 
\frac {n^{-1/2}}2\inf_{\hat\mu_n}\left(1-\min_{\mu\in\{\mu_0,\mu_1\}}P_{\mu^{\otimes n}} (\hat\mu_n= \mu)\right).
\end{equation*}
Hence, it follows from Birg\'e's inequality \cite[Theorem B.10]{giraud2021introduction} that with $K$ the Kullback divergence we have
\begin{eqnarray*}
\inf_{\tilde \mu_n}\sup_{\mu\in\{\mu_0,\mu_1\}}E_{(\mu\star\mu_\epsilon)^{\otimes n}}W_1(\mu,\tilde\mu_n)&\geq& 
\frac {n^{-1/2}}2\left(1-\frac{2e}{2e+1}\vee \frac{nK(\mu_0,\mu_1)}{\log2}\right).
\end{eqnarray*}
We have $K(\mu_0,\mu_1)=n^{-1}/2$ so the right-hand side of the above display is larger than $n^{-1/2}/13$. This completes the proof of \eqref{eq: Lsmallsigma}, whence the proof ov \eqref{eq: paramlower} and the proof of Theorem \ref{theo: DMlower}. \hfill{$\Box$}

\subsection{Proof of Theorem \ref{theo: DFMupper}}
Let $Y\sim\mu\star\mu_\epsilon$ for some $\mu\in\mathcal D(M,a)$. Since $E\vert \delta\vert^{a+2}<M$, it follows from the Markov inequality together with convexity of the function $x\mapsto \vert x\vert ^{a+2}$ that
\begin{eqnarray*}
\int_0^\infty\sqrt{P(\vert Y\vert \geq x)}dx &\leq& 1+\int_1^\infty\sqrt{x^{-(a+2)}E(\vert Y\vert^{a+2})}dx\\
&\leq& 1+\sqrt{2^{a+1}(1+ A)M}\int_1^\infty x^{-(a+2)/2}dx\\
&\leq& 1+\sqrt{2^{a+1}(1+ A)M}\times \frac{2}{a}
\end{eqnarray*}
where $A>0$ is such that $\sigma_n^{a+2}\leq A$ for all $n$. Hence, there exists $L>0$ that depends only on $a,M,A$ such that
\begin{eqnarray}\label{eq: intsqrt}
\int_0^\infty\sqrt{P(\vert Y\vert \geq x)}dx &\leq& L.
\end{eqnarray}
Moreover, possibly enlarging $L$ we also have
\begin{equation*}
E\vert Y\vert ^{3/2} \leq L.
\end{equation*}
Furthermore, the function $r^{(\ell)}$ is bounded on the interval $[-2,2]$ for every $\ell=0,1,2$ under the assumptions of Theorem \ref{theo: DFMupper}, and with  $\mu^*_\epsilon$ the Fourier transform of the distribution of $\epsilon:=\sigma_n\delta$,  and $r_\epsilon=1/\mu^*_\epsilon$, we have $r_\epsilon(x)=r_\delta(\sigma_n x)$ for all $x\in\R$. Combining all of this with  \cite[Proposition 3.1] {dedecker2015improved} (with $p=1$ and $\rho$ sufficiently small) proves that for all $h\in(0,1)$, we can find an estimator $\tilde\mu_{n,h}$ of $\mu$, and a constant $L>0$ that depends only on $a,M,A$ and the distribution of $\delta$ such that for all $n$,
\begin{eqnarray*}
&&E_{(\mu\star\mu_\epsilon)^{\otimes n}}W_1(\tilde\mu_{n,h},\mu)\\
&&\qquad \leq L\left(h+n^{-1/2}+n^{-1/2}\left(\sum_{\ell=0}^1\int_{-1/h}^{1/h}\frac{\sigma_n^{2\ell}\vert r_\delta^{(\ell)}(\sigma_nx)\vert ^2}{\vert x\vert ^2}1_{[-1,1]^c}(x)dx\right)^{1/2}\right).
\end{eqnarray*}
If Assumption  \ref{assm:error_char_deriv}
 holds, then we conclude that there exists $L>0$ that depends only on $a,M,A$ and the distribution of $\delta$ such that for all $n$,
\begin{eqnarray}\label{eq: gub}\notag
&&E_{(\mu\star\mu_\epsilon)^{\otimes n}}W_1(\tilde\mu_{n,h},\mu)\\ \notag
&&\qquad\leq
 L\left(h+n^{-1/2}+n^{-1/2}\exp\left(\vert\sigma_n /h\vert^{\beta}/\gamma_2\right)\times\left(\int_1^{1/h}( x ^{-2}+\sigma_n^{2\tilde\beta} x^{2\tilde\beta-2})dx\right)^{1/2}\right)\\
&&\qquad\leq
 2L\left(h+n^{-1/2}\exp\left(\vert\sigma_n /h\vert^{\beta}/\gamma_2\right)\times\left(1+A(n,h,\tilde\beta)\right)\right)
\end{eqnarray}
where
\begin{eqnarray*}
A(n,h,\tilde\beta)=\begin{cases}
\sigma_n^{1/2}(1/h)^{\tilde\beta-1/2} &\mbox{ if }\tilde\beta>1/2\\
\sigma_n^{1/2}\log(1/h)&\mbox{ if }\tilde\beta=1/2\\
1&\mbox{ if }\tilde\beta<1/2\\
\end{cases}
\end{eqnarray*}

Assume that $\sigma_n\geq n^{-1/2}$ and choose
\begin{equation*}
h=\sigma_n \left(C\gamma_2\log \left(n\sigma_n^2\log n\right)\right)^{-1/\beta}
\end{equation*}
for some constant $C>0$ small enough so that
\begin{equation*}
\eta>\frac{C}{1-2C}.
\end{equation*}
Note that $h$ belongs to $(0,1)$ as required. Possibly enlarging $L$, it follows from \eqref{eq: gub} that with $\hat\mu_Z=\tilde\mu_{n,h}$ with $h$ taken from the previous display, we have
\begin{equation*}
E_{(\mu\star\mu_\epsilon)^{\otimes n}}W_1(\hat\mu_Z,\mu)\leq 
L\left(\sigma_n\left(\log \left(n\sigma_n^2\log n\right)\right)^{-1/\beta} +n^{-1/2}(n\sigma_n^2\log n)^C\left(1+A(n,h,\tilde\beta)\right)\right).
\end{equation*}
The first term on the right-hand side is of larger order than the second one if $\sigma_n\geq n^{-1/2}(\log n)^{\eta}$ so in that case, possibly enlarging $L$ we obtain the upper bound
\begin{equation*}
E_{(\mu\star\mu_\epsilon)^{\otimes n}}W_1(\hat\mu_Z,\mu)\leq 
L\sigma_n\left(\log \left(n\sigma_n^2\log n\right)\right)^{-1/\beta}.
\end{equation*}
If $\sigma_n\in(n^{-1/2}, n^{-1/2}(\log n)^{\eta})$
we obtain the upper bound
\begin{equation*}
E_{(\mu\star\mu_\epsilon)^{\otimes n}}W_1(\hat\mu_Z,\mu)\leq 
Ln^{-1/2}\left(\log \log n\right)^{-1/\beta}(\log n)^{\eta}.
\end{equation*}
Now, assume that $\sigma_n\leq n^{-1/2}$. Choosing $h=n^{-1/2}$ yields an upper bound of the order $n^{-1/2}$, which completes the proof of Theorem \ref{theo: DFMupper}. \hfill{$\Box$}

\subsection{Proof of Theorem \ref{theo: DMlowerUL}}
\subsubsection{Lower bound for $\mu_X$}
It follows from the definition of minimax risk for the unlinked regression model that
\begin{equation}
\label{eq: lowerindep}
\mathcal R_U(M,a,n, \mu_X)\geq \inf_{\widehat m}\sup_{m\in\mathcal M(M,a)}E_{m}\left[\int_0^1\vert\widehat m-m\vert d \mu_X\right]
\end{equation}
where for fixed $m\in\mathcal M(M,a)$, $E_{m}$ denotes the expectation under the distribution of observations $\cX_n \cup \cY_n$, where $\cX_n $ is independent of $\cY_n$. Note that we can restrict the infimum on the right-hand side to those $\widehat m$ for which $E_m\int|\widehat m|d\mu_X$ is finite for all $m$ since otherwise the supremum is infinite. 
Let $\widehat m$ be an arbitrary such estimator. 
If the $X$-sample is independent of the $Y$-sample, then the conditional expectation $\widehat{m}^Y:=E_m(\widehat m|Y_1,\dots, Y_{n})$ does not depend on $m$ and can be viewed as an estimator of $m_0$ based solely on the $Y$-sample. By Jensen's inequality for conditional expectations,
\begin{eqnarray*}
E_{m}\vert\widehat m-m\vert \geq 
E_{m}\vert\widehat {m}^Y-m\vert 
\end{eqnarray*}
and therefore, \eqref{eq: lowerindep} still holds with the infimum restricted to estimators $\widehat m$ that are solely based on observations $Y_1,\dots,Y_{n}$.

It is immediate from \eqref{eq: mum}, that the measure $\mu_m$ of $m(X)$ has finite first moment for all $m\in\mathcal M(M,a)$. This implies that in the definition \eqref{eq: defRU} of the minimax risk with $\mu=\mu_X$, we can restrict the infimum to estimators $\widehat m$ such that the corresponding measure $\mu_{\widehat m}$ has finite first  moment since otherwise, the supremum is infinite. 
Hence, it follows from  \eqref{eq: L1toW1} that
\begin{equation*}
\mathcal R_U(M,a,n, \mu_X) \geq \inf_{\widehat m}\sup_{m\in\mathcal M(M,a)}E_{m}W_1(\mu_{\widehat m},\mu_m).
\end{equation*}
Since for all $m$, $\mu_m$ is the distribution of $m(X)$ (with distribution function $F_X\circ m^{-1}$), $\mu_{\widehat m}$ is a pseudo-estimator of the distribution of $m_0(X)$. Here, we call pseudo-estimator an estimator that is allowed to depend on $F_X$. Hence,
\begin{equation*}
\mathcal R_U(M,a,n, \mu_X) \geq \inf_{\hat \mu}\sup_{m\in\mathcal M(M,a)}E_{m}W_1(\hat\mu,\mu_m)
\end{equation*}
where the infimum is extended to all pseudo-estimators of the distribution of $m_0(X)$ that relies only on the $Y$-sample. As the $Y$-sample has distribution $(\mu_m\star\mu_\epsilon)^{\otimes n_y}$ under $E_m$ (for any  $m\in\mathcal M(M,a)$), we conclude: 
\begin{equation*}
\mathcal R_U(M,a,n, \mu_X) \geq\inf_{\hat \mu}\sup_{m\in\mathcal M(M,a)}E_{(\mu_m\star\mu_\epsilon)^{\otimes n_y}}W_1(\hat\mu,\mu_m).
\end{equation*}

Recall that $m$ belongs to $\mathcal M(M,a)$ if and only if the corresponding $\mu_m$ belongs to $\mathcal D(M,a)$ (see \eqref{eq: mum}). Moreover, any $\mu\in \mathcal D(M,a)$ can be written in the form $\mu=\mu_m$ with $m\in\mathcal M(M,a)$ such that $m^{-1}=F_X^{-1}\circ L$ with $L$ the distribution function of $\mu$. Note that we use here the fact that $F_X\circ F_X^{-1}$ is the identity function by continuity of  the distribution of $X$. Hence, 
\begin{eqnarray*}
\mathcal R_U(M,a,n, \mu_X) & \geq &\inf_{\hat\mu}\sup_{\mu\in\mathcal D(M,a)}E_{(\mu\star\mu_\epsilon)^{\otimes n_y}}W_1(\hat\mu,\mu)
\end{eqnarray*}
where the infimum is taken over all pseudo-estimators that are based on the $Y$-sample. The conclusion of Theorem \ref{theo: DMlowerUL} now follows from Theorem \ref{theo: DMlower}.

\subsubsection{Lower bound for $\hat \mu_X$}\label{sec: prooflbemp}
For arbitrary estimator $\widehat m$, let $\widetilde m$  be an estimator in $\mathcal M(M,a)$ such that
\begin{equation*}
\int|\widehat m-\widetilde m|d\hat\mu_X \leq \inf_{m\in\mathcal M(M,a)}\int|\widehat m-m|d\hat\mu_X +\eta
\end{equation*}
where $\eta>0$ is arbitrarily small.
It then follows from the triangle inequality that for all $m\in\mathcal M(M,a)$
\begin{eqnarray*}
\int|\widetilde m-m|d\hat\mu_X &\leq& \int|\widehat m-m|d\hat\mu_X + \int|\widehat m-\widetilde m|d\hat\mu_X\\
&\leq&2 \int|\widehat m-m|d\hat\mu_X +\eta
\end{eqnarray*}
and therefore,
\begin{equation*}
\mathcal R_U(M,a,n, \hat \mu_X) \geq \frac 12\left(\inf_{\widetilde m \in\mathcal M(M,a)}\sup_{\pi}E_{\pi}\int|\widetilde m-m|d\hat\mu_X -\eta\right).
\end{equation*}
Combining with  Lemma \ref{lem: emptopop} yields that there exists $C>0$ such that
\begin{equation*}
\mathcal R_U(M,a,n, \hat \mu_X) \geq \frac 12\left(\inf_{\widetilde m \in\mathcal M(M,a)}\sup_{\pi}E_{\pi}\int|\widetilde m-m|d\mu_X-Cn^{-1/2} -\eta\right).
\end{equation*}
We choose $\eta=Cn^{-1/2}$. If $\sigma_n\geq C' n^{-1/2}(\log n)^{1/\beta}\vee n^{-1/2}$ for some sufficiently large $C'>0$, then  it follows from the lower bound already obtained for $\mu_X$ that there exists $c>0$ such that the infimum in the right-hand side above is larger than 
$
c\sigma_n\left(\log n\right)^{-1/\beta}
$
which is larger than $3Cn^{-1/2}$. Therefore,
\begin{equation*}
\mathcal R_U(M,a,n, \hat \mu_X) \geq \frac {c\sigma_n}6\left(\log n\right)^{-1/\beta} \,.
\end{equation*}
Now we consider the case when $\sigma_n \le C' n^{-1/2}(\log n)^{1/\beta}\vee n^{-1/2}$. As before, consider only the estimators of the form $\hat m$ that solely relies on $\cY_n$. 
For such estimators,
\begin{eqnarray*}
E_{m}\left[\int\vert\widehat m-m\vert d\hat \mu_X\right]&=& \frac 1n\sum_{i=1}^n E_m\vert\widehat m(X_i)-m(X_i)\vert\\
&=& E_m\int |\widehat m-m|d\mu_X.
\end{eqnarray*}
Hence, similar as above it follows from  \eqref{eq: L1toW1} that
\begin{eqnarray*}
\inf _{\widehat m} \sup_{\pi}E_{\pi}\left[\int\vert\widehat m-m\vert d\hat \mu_X\right]
&\geq & \inf_{\hat \mu}\sup_{m\in\mathcal M(M,a)}E_{(\mu_m\star\mu_\epsilon)^{\otimes n}}W_1(\hat\mu,\mu_m)\\
&\geq & \inf_{\hat \mu}\sup_{\mu\in\mathcal D(M,a)}E_{(\mu\star\mu_\epsilon)^{\otimes n}}W_1(\hat\mu,\mu)
\end{eqnarray*}
where  the infimum in the right-hand side is extended to all pseudo-estimators of the distribution of $m_0(X)$ based on the $Y$-sample. It follows from Theorem \ref{theo: DMlower} that there exists $C>0$ such that the right-hand side is larger than $Cn^{-1/2}$, which completes the proof of 
Theorem \ref{theo: DMlowerUL}. \hfill{$\Box$}

\subsection{Proof of  Corollary \ref{cor: DMlowerUL}}
Assume that $X$ has a density that is bounded from above and bounded away from zero.
Consider an arbitrary estimator $\widehat m$ and let $\widehat m_0\in\cal M$ be an estimator such that 
\begin{equation*}
\int_0^1|\widehat m(x)-\widehat m_0(x)|dx \leq\inf_{m\in\cal M}\int_0^1|\widehat m(x)-m(x)|dx+\eta
\end{equation*}
where $\eta>0$ is arbitrarily small. We consider here integrals with respect to the Lebesgue measure because by assumption, it is known and equivalent to the unknown distribution of $X$. It follows from the triangle inequality and the definition of $\widehat m_0$ that for all $m\in\cal M$
\begin{eqnarray*}
\int_0^1|\widehat m_0(x)-m(x)|dx
&\leq& \int_0^1|\widehat m_0(x)-\widehat m(x)|dx+\int_0^1|\widehat m(x)-m(x)|dx\\
&\leq&2\int_0^1|\widehat m(x)-m(x)|dx+\eta.
\end{eqnarray*}
Denoting by $c_X$ and $C_X$ two positive numbers such that $c_X\leq f_X\leq C_X$
one obtains
\begin{eqnarray*}
\int_0^1|\widehat m_0(x)-m(x)|d\mu_X(x)
&\leq& C_X \int_0^1|\widehat m_0(x)-m(x)|dx\\
&\leq&\frac{2C_X}{c_X}\int_0^1|\widehat m(x)-m(x)|d\mu_X(x)+C_X\eta.
\end{eqnarray*}
Hence,
\begin{equation*}
\inf_{\widehat m}\sup_{\pi}E_{\pi}\left[\int\vert\widehat m-m\vert d\mu_X\right]
\geq
\frac{c_X} {2C_X}\left(\inf_{\widehat m_0\in\cal M}\sup_{\pi}E_{\pi}\left[\int\vert\widehat m_0-m\vert d\mu_X\right]-C_X\eta\right).
\end{equation*}
Choosing $\eta$ small enough so that $C_X\eta$ is smaller than half the infimum in the right-hand side then yields
\begin{equation*}
\inf_{\widehat m}\sup_{\pi}E_{\pi}\left[\int\vert\widehat m-m\vert d\mu_X\right]
\geq
\frac{c_X} {4C_X}\inf_{\widehat m\in\cal M}\sup_{\pi}E_{\pi}\left[\int\vert\widehat m-m\vert d\mu_X\right]
\end{equation*}
so same the lower bound on $\mu_X$ than in  Theorem \ref{theo: DMlowerUL}, with now the infimum extended to all possible estimators, follows from Theorem \ref{theo: DMlowerUL}.

To show that the lower bound on $\hat\mu_X$ in Theorem \ref{theo: DMlowerUL} still holds if the infimum is extended to all possible estimators, for  arbitrary estimator $\widehat m$ we defined $\widehat m_0\in\cal M$ as an estimator such that 
\begin{equation*}
\int_0^1|\widehat m-\widehat m_0|d\hat\mu_X \leq\inf_{m\in\cal M}\int_0^1|\widehat m-m|d\hat\mu_X+\eta
\end{equation*}
where $\eta>0$ is arbitrarily small. Similar as above we then have
\begin{eqnarray*}
\int_0^1|\widehat m_0-m|d\hat\mu_X
&\leq&2\int_0^1|\widehat m-m|d\hat\mu_X+\eta
\end{eqnarray*}
  for all $m\in\cal M$, so repeating similar arguments as above completes the proof of
Corollary \ref{cor: DMlowerUL}
 \hfill{$\Box$}

\subsection{Proof of Theorem \ref{theo: DMUpperUL}}
\label{sec: proofDMUpperUL}

We consider a slight modification of the set $\cM(M, a)$, namely $\cM^0(M, a)$, which is the set of all  functions $m\in\mathcal M$ with: 
\begin{equation*}
\int \vert m(x)\vert^{a+2}dx\leq M/c_X.
\end{equation*}
It is immediate that $\mathcal M(M,a)\subset \mathcal M^0(M,a)$. 
Note that, the set $\mathcal M^0(M, a)$ does not depend on $\mu_X$, except on $c_X$ that is assumed to be known. 

We will build an estimator which maximal risk is smaller than $v_n$ in order. The  estimator is built  in two steps. First, note that from Theorem \ref{theo: DFMupper}, it follows that there exists  an estimator $\hat\mu_{Z}$  of the distribution of $Z:=m_0(X)$ based on the $Y$-sample that satisfies
\begin{equation*}
\sup_{\mu\in\mathcal D(M,a)}E_{(\mu\star\mu_\epsilon)^{\otimes n}}W_1(\mu,\hat\mu_Z)\leq Lv_n.
\end{equation*}
Second,  we define our estimator for $m_0$ as the minimum constrast estimator with respect to $W_1$-distance:  the estimator $\widehat m_n$ belongs to $\mathcal M^0(M,a)$ and satisfies
\begin{equation*}
W_1\left(\hat\mu_{Z},\hat\mu_{ \widehat m_n}\right)\leq \inf_{m\in\mathcal M^0(M,a)}W_1\left(\hat\mu_{Z},\hat\mu_{ m}\right)+\eta
\end{equation*}
where $\eta>0$ is arbitrarily small. This is indeed an estimator since it does not depend on unknown parameters. Here for each $m \in \mathcal M^0(M,a)$,  $\hat \mu_m$ is the empirical measure    that puts mass $n^{-1}$ at each $m(X_i$), $i\in\{1,\dots,n\}$. Hence, $\hat\mu_m$ estimates the distribution of $m(X)$ and our estimator minimizes the $W_1$-distance between this estimator and an estimator of the distribution of $m_0(X)$.

 Since $X$ has a bounded density, $\mathcal M^0(M,a)$ is included in $\mathcal M(M',a)$ where $M'=C_XM/c_X$ with $c_X > 0$ and $C_X < \infty$ are lower bound and the upper bound on the density of $X$. Hence, $\widehat m_n$ belongs to $\mathcal M(M',a)$ so
it follows from Lemma \ref{lem: emptopop} (with $M$ replaced by $M'$) that there exists $L>0$ such that for all probability measures $\pi$ that belong to $\Pi(m,n)$ for some $m\in\mathcal M(M,a)$ we have

\begin{eqnarray*}
E_\pi\int\vert \widehat m_n-m\vert d\mu_X&\leq&E_\pi\int\vert \widehat m_n-m\vert d\hat \mu_X+Ln^{-1/2}.
\end{eqnarray*}
For the first term on the right-hand side, it follows from \eqref{eq: L1toW1emp} that
\begin{eqnarray}\label{eq: empW1}
\int\vert \widehat m_n-m\vert d\hat \mu_X&=&W_1(\hat\mu_{\widehat m_n},\hat \mu_{m}).
\end{eqnarray}
Since $\mathcal M(M,a)\subset \mathcal M^0(M,a)$, it follows from the triangle inequality combined to the definition of $\widehat m_n$ that for all $m\in\mathcal M(M,a)$,
\begin{eqnarray*}
\int\vert \widehat m_n-m\vert d\hat \mu_X&\leq&W_1(\hat\mu_{\widehat m_n},\hat \mu_Z)+W_1(\hat\mu_Z,\hat \mu_{m})\\&\leq&2W_1(\hat\mu_Z,\hat \mu_{m})+\eta\\
&\leq&2W_1(\hat\mu_Z,\mu_{m})+2W_1(\hat\mu_m,\mu_{m})+\eta
\end{eqnarray*}
where $\eta>0$ is arbitrarily small. In the sequel we choose $\eta=n^{-1/2}$.
Due to \eqref{eq: mum}, we have $\mu_m\in\mathcal D(M,a)$ if and only if $m\in\mathcal M(M,a)$ so  we obtain that there exists $L>0$ such that
\begin{eqnarray}
\label{eq: ent2}
E_\pi\int\vert \widehat m_n-m\vert d\mu_X&\leq&Lv_{n}+Ln^{-1/2}+2E_\pi W_1(\hat\mu_m,\mu_{m}).
\end{eqnarray}
Now, $\hat \mu_m$ is the empirical distribution of the sample $m(X_1),\dots,m(X_n)$ whereas $\mu_m$ is the common distribution of those variables so by
\cite[Theorem 3.2]{bobkov2019one}, with $F_m$ the distribution function of $\mu_m$ we have
\begin{eqnarray*}
E_\pi W_1(\hat\mu_m,\mu_{m})&\leq& n^{-1/2}\int_{-\infty}^\infty \sqrt{F_m(t)(1-F_m(t))}dt\\
&\leq& 2n^{-1/2}\int_0^\infty\sqrt{P\left(|m(X)|>t\right)}dt\\
&\leq& 2n^{-1/2}\left(1+\frac{2\sqrt{M}}{a}
\right)
\end{eqnarray*}
where the last inequality is obtained with similar arguments as for the proof of \eqref{eq: intsqrt}. 
Combining the two last displays yields that there exists $L>0$ such that
\begin{eqnarray*}
\sup_{m\in\mathcal M(M,a)}\sup_{\pi\in \Pi(m,n)} E_\pi\int\vert \widehat m_n-m\vert d\mu_X&\leq&Lv_{n}+Ln^{-1/2}.
\end{eqnarray*}
Hence,
\begin{eqnarray*}
\mathcal R_U(M,a,n,\mu)&\leq&
\sup_{m\in\mathcal M(M,a)}\sup_{\pi\in \Pi(m,n)} E_\pi\int\vert \widehat m_n-m\vert d\mu_X\\&\leq&Lv_{n}+Ln^{-1/2}.
\end{eqnarray*}
The inequality in Theorem \ref{theo: DMUpperUL} with $\mu=\mu_X$ follows since $v_n\geq n^{-1/2}$.

Now, since (as proved above)  $\widehat m_n$ belongs to $\mathcal M(M',a)$ where $M'=C_XM/c_X$,  the inequality with $\mu=\hat\mu_X$  follows from  Lemma \ref{lem: emptopop} together with the the inequality with $\mu=\mu_X$. This completes the proof of Theorem \ref{theo: DMUpperUL}.

\subsection{Proof of Theorem \ref{thm:shuffled_random_upper}}
We use the same notation $\mathcal M^0(M,a)$ as in Section \ref{sec: proofDMUpperUL}.
Let $\hat\mu_Y$ denote the empirical distribution generated by $\cY_n$, i.e. $\hat\mu_Y$ puts mass $n^{-1}$ at each $Y_i$, for $i\in\{1,\dots,n\}$. For any $m \in \cM(M, a)$, we further recall that  $\hat\mu_{m}$ denotes the distribution  that puts mass $n^{-1}$ at each $m(X_i$), $i\in\{1,\dots,n\}$. With these notations at our disposal, we define the minimum contrast estimator $\widetilde m_n$ as an element of $\mathcal M^0(M,a)$ that satisfies : 
 \begin{equation}\label{eq: estimatorp}
W_2\left(\hat\mu_{Y},\hat\mu_{ \widetilde m_n}\right)\leq \inf_{m\in\mathcal M^0(M,a)}W_2\left(\hat\mu_{Y},\hat\mu_{ m}\right)+\eta
\end{equation}
where $\eta>0$ is arbitrarily small (we will choose $\eta=\sigma_n^2$) and the $W_2$-Wasserstein distance beetween probability measures $\mu$ and $\mu'$ on $\R$ is defined by
\begin{equation}\label{def: Wassersteinp}
W_2^2(\mu,\mu')=\inf_{\pi\in\Pi(\mu,\mu')}\int_{\R\times\R}\vert x-y\vert^2\pi(dx,dy) \,,
\end{equation}
where $\Pi(\mu,\mu')$ is the set of all probability measures on $\R\times\R$ with marginals $\mu$ and $\mu'$. We will also consider a piece-wise constant version of the estimator $\widetilde m_n$: let $\widehat m_n$ denote the estimator that coincides with $\widetilde m_n$ at each order statistic $X_{(i)}$, for $i\in\{1,\dots,n\}$ and that is constant on the intervals $[0,X_{(1)}]$, $[X_{(n)},1]$ and on all intervals $(X_{(i-1)},X_{(i)}]$ for $i\in\{2,\dots,n\}$.
Let  $m\in\mathcal M(M,a)$, and $U$ be a random variable independent of the $(X_i, \eps_i)$'s, with uniform distribution on the finite set $\{1,\dots,n\}$. Define
\begin{equation*}
T_m=\sum_{i=1}^nm(X_i)\mathds{1}_{U=i}, \ \qquad \ T_Y=\sum_{i=1}^nY_{i}\mathds{1}_{U=i}  \,.
\end{equation*}
Conditionally on the $X_i$'s and $\delta_i$'s, $T_m$ has distribution $\hat \mu_m$ and the distribution of $T_Y$ is $\hat \mu_Y$. By definition \eqref{def: Wassersteinp} of the $W_2$-Wasserstein distance, with $\bbE_{U \mid (X, \delta)}$ the conditional expectation given the $X_i$'s and $\delta_i$'s, and $m_0=m$, one concludes that 
\begin{eqnarray*}
W^2_2\left(\hat \mu_Y, \hat \mu_m\right)  \le \bbE_{U \mid (X, \delta)}\left| T_m-T_Y\right|^2 =\sigma_n^2\bbE_{U \mid (X, \delta)}\left|\sum_{i=1}^n\delta_{i} \mathds{1}_{U=i} \right|^2 = \frac{\sigma_n^2}n \sum_{i=1}^n \delta_i^2.
\end{eqnarray*}
Let $m_n$ denote either $\widehat m_n$ or $\widetilde m_n$. 
By \cite[Theorem 2.10]{bobkov2019one} (that extends \eqref{eq: L1toW1emp} for $W_1$ to $W_p$ for arbitrary $p\geq 1$), 
\begin{eqnarray*}
\int\vert m_n-m\vert^2 d\hat \mu_X& =&\frac1n \sum_{i = 1}^n \left(\widetilde m_n(X_{(i)}) - m(X_{(i)})\right)^2 \\ \notag 
& = &W_2^2(\hat \mu_{\widetilde m_n}, \hat \mu_m) 
\end{eqnarray*}
so it follows from the triangle inequality  that
\begin{eqnarray*}
\int\vert m_n-m\vert^2 d\hat \mu_X&\leq&2\left(W_2^2(\hat\mu_{\widetilde m_n},\hat \mu_Y)+W_2^2(\hat\mu_Y,\hat \mu_{m})\right).
\end{eqnarray*}
Since $\mathcal M(M,a)\subset \mathcal M^0(M,a)$, it follows from the previous inequality combined to  the definition of $\widetilde m_n$ that for all $m\in\mathcal M (M,a)$, if $m_0=m$ then we have
\begin{eqnarray*}
\int\vert m_n-m\vert^2 d\hat \mu_X
\leq 4 W_2^2(\hat\mu_Y,\hat \mu_{m})+2\eta.
\end{eqnarray*}
Hence we arrive at
\begin{eqnarray}\label{eq:emp_W2}\notag
\int\vert m_n-m\vert^2 d\hat \mu_X
&=&\frac1n \sum_{i = 1}^n \left(\widetilde m_n(X_{(i)}) - m(X_{(i)})\right)^2 \\
&\leq& \frac{4\sigma_n^2}n \sum_{i=1}^n \delta_i^2+2\eta\,.
\end{eqnarray}
Taking expectation  on both sides and combining with Holder inequality proves that
$$
\sup_{m \in \cM(M, a)} E_{m}\left[\int\left| m_n- m\right| \ d\hat\mu_X\right] \le K\left(\sigma_n+\eta^{1/2}\right).
$$
Choosing $\eta=\sigma_n^2$ yields 
\begin{equation}\label{eq: B1}
\mathcal R_S(M,a,n,\hat\mu_X)\le K\sigma_n.
\end{equation}
This  proves the first claim of Theorem  \ref{thm:shuffled_random_upper}.
Now, the shuffled regression model is a special case of the unlinked regression model in the sense that the distribution corresponding to $E_m$ in the shuffled regression model belongs to $\Pi(m,n)$ for all $m$, so using the notation of Section \ref{sec:unlinked}, for all estimators $\widehat m_n$  we have
\begin{equation}\label{eq:SvsU}
\sup_{m \in \cM(M, a)} E_{m}\left[\int\left| \widehat m_n- m\right| \ d\hat\mu_X\right] \le \sup_{m\in\mathcal M(M,a)}\sup_{\pi\in \Pi(m,n)}E_{\pi}\left[\int|\widehat m_n-m|d\hat\mu_X \right]
\end{equation}
and therefore,
\begin{eqnarray*}
\mathcal R_S(M,a,n,\hat\mu_X)&\le& \mathcal R(M,a,n,\hat\mu_X)\\
&\leq& Lv_n
\end{eqnarray*}
provided that the assumptions of the Theorem \ref{theo: DMUpperUL}  are satisfied. This completes the proof of Theorem  \ref{thm:shuffled_random_upper} for the empirical minimax risk since for the upper bound, we can consider the smallest of the two  bounds above.

Next, we turn to the population minimax risk. For this task, we calculate the distance between $\widehat m_n$ and $m$ with respect to $\mu_X$, where $m\in\mathcal M(M,a)$. With $X_{(0)}:=0$ we have
\begin{align*}
& \int |\widehat m_n - m| \ d\mu_X \\
& = \sum_{j = 1}^{n} \int_{X_{(j-1)}}^{X_{(j)}} |\widehat m_n(x) - m(x) |  \ d\mu_X(x) + \int_{X_{(n)}}^1 |\widehat m_n(x) - m(x) |  \ d\mu_X(x) \\
& =  \sum_{j = 1}^n \int_{X_{(j-1)}}^{X_{(j)}} |\widetilde m_n(X_{(j)}) - m(x) |  \ d\mu_X(x) + \int_{X_{(n)}}^1 |\widetilde m_n(X_{(n)}) - m(x) |  \ d\mu_X(x) \\
& \le \sum_{j = 2}^n \max\left\{|\widetilde m_n(X_{(j)}) - m(X_{(j-1)}) |, |\widetilde m_n(X_{(j)}) - m(X_{(j)}) | \right\}(F_X(X_{(j)}) - F_X(X_{(j-1)})) \\
& \qquad \qquad + \int_{0}^{X_{(1)}} |\widetilde m_n(X_{(1)}) - m(x) |  \ d\mu_X(x) +  \int_{X_{(n)}}^1 |\widetilde m_n(X_{(n)}) - m(x) |  \ d\mu_X(x)
\end{align*}
where the sum in the last term is bounded above by
\begin{align*}
&\sum_{j = 2}^n |\widetilde m_n(X_{(j)}) - m(X_{(j)}) |(F_X(X_{(j)}) - F_X(X_{(j-1)})) \\
& \qquad + \sum_{j = 2}^n (m(X_{(j)}) - m(X_{(j-1)})) (F_X(X_{(j)}) - F_X(X_{(j-1)})) \,.
\end{align*}
Hence,
\begin{align*}
 \int |\widehat m_n - m| \ d\mu_X
 \le T_1+T_2+T_3+T_4
\end{align*}
where
\begin{align*}
T_1& =  \sum_{j = 1}^n |\widetilde m_n(X_{(j)}) - m(X_{(j)}) |(F_X(X_{(j)}) - F_X(X_{(j-1)}))\\
T_2& =\sum_{j = 2}^n (m(X_{(j)}) - m(X_{(j-1)})) (F_X(X_{(j)}) - F_X(X_{(j-1)}) )\\
T_3&= |\widetilde m_n(X_{(n)}) - m(X_{(n)}) | (1 - F_X(X_{(n)}))\\
T_4&= \int_{0}^{X_{(1)}} |m(X_{(1)}) - m(x) |  \ d\mu_X(x) +  \int_{X_{(n)}}^1 | m(X_{(n)}) - m(x) |  \ d\mu_X(x) \,.
\end{align*}
We now analyze each $T_i$ separately. The rest of the proof relies on some properties of order statistics, which we state below: 
\begin{lemma}
\label{lem:order_stat}
If $X$ has a continuous distribution function $F$, then $F(X) \sim \unif([0, 1])$. As a consequence, if $X_1, \dots, X_n \stackrel{i.i.d.}{\sim} F$ with the order statistics $X_{(1)}, \dots, X_{(n)}$, then $F(X_{(1)}), \dots, F(X_{(n)})$ have same distribution as $U_{(1)}, \dots, U_{(n)}$, n order statistics from uniform distribution on $[0,1]$. Furthermore, if we define $\delta_{n, i} = U_{(i)} - U_{(i-1)}$ for $1 \le i \le n+1$ (with the convention $U_{(0)} = 0$ and $U_{(n+1)} = 1$), then: 
$$
\left(\delta_{n, 1}, \dots, \delta_{n, n+1}\right) \overset{d}{=} \left(\frac{\alpha_1}{\sum_{i = 1}^{n+1} \alpha_i}, \dots, \frac{\alpha_{n+1}}{\sum_{i = 1}^{n+1} \alpha_i}\right)
$$
where $\alpha_1, \dots, \alpha_{n+1} \overset{i.i.d.}{\sim} \Exp(1)$.
\end{lemma}
For $T_1$ we have
\begin{align*}
T_1 & \le \sqrt{\sum_{j = 1}^n (\widetilde m_n(X_{(j)}) - m(X_{(j)}))^2} \sqrt{\sum_{j = 1}^n (F_X(X_{(j)}) - F_X(X_{(j-1)}))^2}\,.
\end{align*}
Hence, \eqref{eq:emp_W2} with $\eta=\sigma^2$ yields: 
\begin{align*}
E_m[T_1] & \leq E_m\left[\sqrt{\frac1n \sum_{j = 1}^n (\widetilde m_n(X_{(j)}) - m(X_{(j)}))^2} \sqrt{n\sum_{j = 1}^n (F_X(X_{(j)}) - F_X(X_{(j-1)}))^2}\right] \\
& \le 2\sigma_n E_m\left[\sqrt{\frac1n \sum_{i=1}^n \delta_i^2 }\right]E_m\left[\sqrt{n\sum_{j = 1}^n (F_X(X_{(j)}) - F_X(X_{(j-1)}))^2}\right] \\
& \le C\sigma_n \|\delta\|_2 
\end{align*}
where $C>0$ is a universal constant.
Here the last line follows from an application of Jensen's inequality and using the fact that $F_X(X) \sim \unif([0, 1])$, see Lemma \ref{lem:order_stat}. The later fact implies that (omitting the subscript $m$ since the expectation below does not depend on $m$)
$$
E[(F_X(X_{(j)}) - F_X(X_{(j-1)}))^2] = E[(U_{(j)} - U_{(j-1)})^2]  \le Kn^{-2}  
$$
for some constant $K$. Hence, fot $T_1$ we have
$$E_m(T_1)\leq C\sigma_n$$
for some universal constant $C>0$.

Now for $T_2$ we use the following: 
\begin{align*}
T_2& \le \left[\sup_j \left(F_X(X_{(j)}) - F_X(X_{(j-1)})\right) \right]\left( m\left(X_{(n)}\right) -  m\left(X_{(1)}\right)\right)\end{align*}
where by monotonicity of $m$,
\begin{align*}
m\left(X_{(n)}\right) -  m\left(X_{(1)}\right)
& \le\left(\frac{1}{1 - F_X(X_{(n)})} \int_{X_{(n)}}^1 m\ d\mu_X - \frac{1}{F_X(X_{(1)})} \int_0^{X_{(1)}} m\ d\mu_X\right) \\
& \le M^{\frac{1}{a+2}} \left(\left(1 - F_X(X_{(n)})\right)^{-\frac{1}{a+2}} + F_X(X_{(1)})^{-\frac{1}{a+2}}\right)\,. 
\end{align*}
Now it follows from Lemma \ref{lem:order_stat}
 that (omitting again the unnecessary subscript $m$)
\begin{align*}
E^{1/2}\left[\sup_j \left(U_{(j)} - U_{(j-1)}\right)^2\right] & = O\left(\frac{\log{n}}{n}\right) \,, \\
E^{1/2}[U_{(1)}^2] & = O(n^{-1}) \,, \\
E^{1/2}[(1 - U_{(n)})^2] & = O(n^{-1}) \,
\end{align*}
where $U_1,\dots,U_n$ are i.i.d. variables with a uniform distribution on $[0,1]$ and $U_{(1)}\leq\dots\dots \le U_{(n)}$ are the corresponding order statistics.
Hence 
\begin{align*}
& E\left[ \left[\sup_j \left(F_X(X_{(j)}) - F_X( X_{(j-1)})\right) \right]  \left(\left(1 - F(X_{(n)})\right)^{-\frac{1}{a+2}} + F(X_{(1)})^{-\frac{1}{a+2}}\right) \right] \\
& = E\left[ \left[\sup_j \left(U_{(j)} - U_{(j-1)})\right) \right]  \left(\left(1 - U_{(n)}\right)^{-\frac{1}{a+2}} + U_{(1)}^{-\frac{1}{a+2}}\right) \right] \\
& \le E^{1/2}\left[\sup_j \left(U_{(j)} - U_{(j-1)})\right)^2 \right]E^{1/2}\left[  \left(1 - U_{(n)}\right)^{-\frac{2}{a+2}} + U_{(1)}^{-\frac{2}{a+2}} \right] \\
& \le Cn^{-1 + \frac{1}{a+2}}\log{n}
\end{align*}
where $C>0$ is a universal constant.
Here in the last line we used the following fact: 
$$
E\left[U_{(1)}^{-\frac{2}{a+2}}\right] \le Cn^{\frac{2}{a+2}}, \ \ \ \ E\left[\left(1 - U_{(n)}\right)^{-\frac{2}{a+2}}\right]  \le Cn^{\frac{2}{a+2}} \,.
$$
We prove below the first inequality, as the proof the second inequality is similar. Note that the distribution function and the density function of $U_{(1)}$ are respectively given by
$$
F_{U_{(1)}}(x) = 1 - (1- x)^n, \ \ \  f_{U_{(1)}}(x) = n(1-x)^{n-1} \,.
$$ 
Therefore we have: 
\begin{align*}
E[U_{(1)}^{-\frac{2}{a+2}}] & = n\int_0^1 x^{-\frac{2}{a+2}}(1-x)^{n-1} \ dx \\
& =  n\int_0^1 x^{\frac{a}{a+2} - 1}(1-x)^{n-1} \ dx \\
& = nB\left(\frac{a}{a+2}, n\right) \\
& = n \frac{\Gamma\left(\frac{a}{a+2}\right)\Gamma(n)}{\Gamma\left(\frac{a}{a+2} + n\right)} \,.
\end{align*}
It follows from Wendall's limit theorem (see 1.1.6 of \cite{qi2010bounds}) for Gamma functions that
$$
\lim_{n \uparrow \infty} \frac{\Gamma(n + \alpha)}{\Gamma(n) n^{\alpha}} = 1\,, \ \ \forall \ \ \alpha > 0 \,.
$$
Therefore for all large $n$ (with $\alpha = a/(a+2))$ we have $\Gamma(n + a/(a+2)) \ge (1/2)\Gamma(n) n^{a/(a+2)}$, whence
\begin{align*}
E[U_{(1)}^{-\frac{2}{a+2}}]  \le 2n\Gamma\left(\frac{a}{a+2}\right)n^{-\frac{a}{a+2}} = 2\Gamma\left(\frac{a}{a+2}\right) n^{\frac{2}{a+2}} \,.
\end{align*}
We have elaborate on the calculation for $U_{(1)}$ and the calculation of $1 - U_{(n)}$ is similar. 
Hence
\begin{equation*}
E_m(T_2)\le Cn^{-1 + \frac{1}{a+2}}\log{n}
\end{equation*}
where $C>0$ depends only on $M$ and $a$.

We next bound $T_3$.  The bound on $T_3$ follows from the similar line of argument as of $T_1$. Note that: 
\begin{align*}
E_m[T_3] & \le E_m\left[\left(\sum_j |\widetilde m_n(X_{(j)}) - m(X_{(j)}) |\right) (1 - F_X(X_{(n)})) \right] \\
& \le \sqrt{E_m\left[\left(\frac1n \sum_j |\widetilde m_n(X_{(j)}) - m(X_{(j)}) |\right)^2\right]}\sqrt{ E\left[\left(n(1 - F_X(X_{(n)}))\right)^2 \right]} \\
& \le \sqrt{E_m\left[\frac1n \sum_j \left(\widetilde m_n(X_{(j)}) - m(X_{(j)}) \right)^2\right]} \sqrt{E\left[\left(n(1 - F_X(X_{(n)}))\right)^2 \right]} \\
& \le 2\sigma_n (\|\delta\|_2 +1) \sqrt{E\left[\left(n(1 - U_{(n)})\right)^2 \right]} \le C \sigma_n 
\end{align*}
where $C>0$ is a universal constant.

Finally for $T_4$, note that: 
\begin{align*}
E_m[T_4] & = E_m\left[\int_{0}^{1} (m(X_{(1)}) - m(x) ) \mathds{1}_{x \le X_{(1)}} \ d\mu_X(x) \right] \\
& \qquad \qquad+ E_m\left[\int_{0}^1 ( m(x) - m(X_{(n)}))\mathds{1}_{x \ge X_{(n)}}  \ d\mu_X(x)\right] \\
& \le E_m\left[m(X_{(1)})F_X(X_{(1)})\right]  + M^{\frac{1}{a+2}}E\left[F_X(X_{(1)})^{1- \frac{1}{a+2}}\right] \\
& \qquad \qquad -E_m\left[m(X_{(n)})\left(1 -F(X_{(n)})\right)\right] +  M^{\frac{1}{a+2}}E\left[\left(1 -F_X(X_{(n)})\right)^{1- \frac{1}{a+2}}\right]  \\
& \le E_m\left[m(X_{(1)})F_X(X_{(1)})\right]  -E_m\left[m(X_{(n)})\left(1 -F(X_{(n)})\right)\right] + Cn^{-1 + \frac{1}{a+2}}
\end{align*}
where $C>0$ depends only on $M$ and $a$.
For the first term, observe that: 
\begin{align*}
E_m[m(X_{(1)})F_X(X_{(1)})] & \le E_m\left[\frac{F_X(X_{(1)})}{F_X(X_{(2)}) - F_X(X_{(1)})} \int_{X_{(1)}}^{X_{(2)}} m(x) \ d\mu_X(x)\right] \\
& \le  M^{\frac{1}{a+2}} E\left[F_X(X_{(1)}) \left(F_X(X_{(2)}) - F_X(X_{(1)})\right)^{-\frac{1}{a+2}}\right] \\
& \le M^{\frac{1}{a+2}}  \sqrt{E\left[F_X^2(X_{(1)})\right]} \sqrt{E\left[\left(F_X(X_{(2)}) - F_X(X_{(1)})\right)^{-\frac{2}{a+2}}\right]} \\
& \le Cn^{-1 + \frac{1}{a+2}}
\end{align*}
where $C>0$ depends only on $M$ and $a$.
. The calculation for the second term is similar and skipped for the brevity. Combining, we get
\begin{equation*}
E_m(T_1+T_2+T_3+T_4)\le C\left(\sigma_n+n^{-1 + \frac{1}{a+2}}\log{n}\right)
\end{equation*}
where $C>0$ depends only on $M$ and $a$. Hence,
\begin{align*}
\sup_{m\in\mathcal M(M,a)}E_m \int |\widehat m_n - m| \ d\mu_X
 \le C\left(\sigma_n+n^{-1 + \frac{1}{a+2}}\log{n}\right)
\end{align*}
and the second claim in Theorem \ref{thm:shuffled_random_upper} follows.  To prove the last claim for the population risk, note that if $\sigma_n\geq n^{-1/2}$ then the above upper bound simplifies to
\begin{align*}
\sup_{m\in\mathcal M(M,a)}E_m \int |\widehat m_n - m| \ d\mu_X
 \le 2C\sigma_n
\end{align*}
since $-1+1/(a+2)<-1/2$.
Then, we again use \eqref{eq:SvsU} and Theorem \ref{theo: DMUpperUL} to prove that the population minimax risk is bounded above by $Lv_n$, and we consider the smallest of the two obtained upper bounds. 
This completes the proof of Theorem \ref{thm:shuffled_random_upper}.  \hfill{$\Box$}

\subsection{Proof of Theorem \ref{thm:mlb_shuffled_random_noiseless}}
  To prove the lower bound, we use Le-Cam's two point method. Towards that end, set $X \sim U([0, 1])$, define $m_0 \equiv 0$ and $\tilde m$ as a function which is $0$ on $(Cn^{-1}, 1]$ (for some constant $C$ to be chosen later) and  $\left(x\log^{1+\eps}{(1/x)}\right)^{-\frac{1}{a+2}}$ on $[0, Cn^{-1}]$. Note that the function $\left(x\log^{1 + \eps}{(1/x)}\right)^{-\frac{1}{a+2}}$ is decreasing in a neighborhood of $0$ and consequently $\tilde m$ is decreasing function. Furthermore we have: 
\begin{align*}
\bbE[|\tilde m(X)|^{a+2}] \le \int_0^{C/n} \left(x\log^{1 + \eps}{(1/x)}\right)^{-1} \ dx \le M \,.
\end{align*}
Hence both $m_0, \tilde m \in \cM(M, a)$. Moreover we have the following lower bound on the $L_1$ distance between $\tilde m$ and $m_0$:
  \begin{align*}
    d(\tilde m, m_0) = \int \left|\tilde m - m_0\right| \ d\mu_X & = \int_0^{Cn^{-1}} \left(x\log^{1 + \eps}{(1/x)}\right)^{-\frac{1}{a+2}}\ dx \\
    & \ge (n/C)^{\frac{1}{a+2}}(\log{(n/C)})^{-\frac{1+\eps}{a+2}} \int_0^{Cn^{-1}} \ dx \\
    & \ge (n/C)^{-1 + \frac{1}{a+2}}(\log{(n/C)})^{-\frac{1+\eps}{a+2}} \triangleq \tau_n \,.
    \end{align*}
    Recall that, for any two measure $P$ and $Q$ the TV (total variation) distance between them is defined as: 
    $$
    \|P -Q\|_\TV = \sup_A |P(A) - Q(A)|
    $$
    for all measurable set $A$. Let  $P_{m_0}$ (resp. $P_{\tilde m}$) denote the distribution of $(X, m_0(X))$ (resp. $(X, \tilde m(X))$).  Using Le-Cam's method we have: 
    \begin{align*}
    & \inf_{\hat m} \sup_{m \in \cM(M, a)} E_m\left[d(\hat m, m)\right] \\
    & \ge \frac{\tau_n}{2}    \inf_{\hat m} \sup_{m \in \cM(M, a)} P^{\otimes n}_m(d(\hat m, m) \ge \tau_n/2)  \\
    & \ge  \frac{\tau_n}{2}    \inf_{\hat m} \max_{m \in \{m_0, \tilde m\}} P^{\otimes n}_m(d(\hat m, m) \ge \tau_n/2) \\
     & \ge \frac{\tau_n}{4} \inf_{\hat m} \left\{P^{\otimes n}_{m_0}(d(\hat m, m_0) \ge \tau_n/2) + P^{\otimes n}_{\tilde m}(d(\hat m, \tilde m) \ge \tau_n/2)\right\} \\
     & \ge \frac{\tau_n}{4} \inf_{\hat m} \left\{P^{\otimes n}_{m_0}(d(\hat m, m_0) \ge \tau_n/2) + P^{\otimes n}_{\tilde m}(d(\hat m,  m_0) < \tau_n/2)\right\}\,,
  \end{align*}
using that $d(m_0, \tilde m) \ge \tau_n$ for the last inequality. Hence,
 \begin{align*}
& \inf_{\hat m} \sup_{m \in \cM(M, a)} E_m\left[d(\hat m, m)\right] \\
        & \ge \frac{\tau_n}{4} \inf_{\hat m} \left\{1 - \left(P^{\otimes n}_{\tilde m}(d(\hat m,  m_0) \ge \tau_n/2) - P^{\otimes n}_{m_0}(d(\hat m, m_0) \ge \tau_n/2) \right)\right\} \\
          & \ge \frac{\tau_n}{4} \left\{1 - \left\|P^{\otimes n}_{m_0} - P^{\otimes n}_{\tilde m}\right\|_\TV \right\} \\
           & = \frac{\tau_n}{4} \left\{1 - n\left\|P_{m_0} - P_{\tilde m}\right\|_\TV \right\} \,.
    \end{align*}
  Now we bound $\left\|P_{m_0} - P_{\tilde m}\right\|_\TV$ for which we use the following fact: 
    \begin{lemma}
        For any two measure $(\mu, \nu)$ we have
$ \|\mu - \nu\|_\TV \le \bbP(X \neq Y)$ for any $(X, Y)$ with $X \sim \mu$ and $Y \sim \nu$. 
    \end{lemma}
\noindent
We now generate a coupling as follows: i) generate $X \sim U([0, 1])$; ii) set $Z_1 = (X, m_0(X)) \equiv (X, 0)$ and iii) set $Z_2 = (X, \tilde m(X))$. Then it is immediate $Z_1 \sim P_{m_0}$ and $Z_2 \sim P_{ \tilde m}$. Consequently: 
\begin{align*}
    \|P_{m_0} - P_{\tilde m}\|_\TV \le \bbP(Z_1 \neq Z_2) = \bbP(X \le Cn^{-1}) =Cn^{-1} \,.
\end{align*}
Therefore, choosing $C = 1/2$, we have: 
$$
\inf_{\hat m} \sup_{m \in \cM(M, a)} E_m\left[\int \left|\hat m - m\right| \ d\mu_X\right] \ge \frac{\tau_n}{8} = \frac{1}{8}(2n)^{-1 + \frac{1}{a+2}}(\log{(2n)})^{-\frac{1+\eps}{a+2}} \,.
$$
This completes the proof of Theorem \ref{thm:mlb_shuffled_random_noiseless}.   \hfill{$\Box$}

\subsection{Proof of Theorem \ref{thm:mlb_shuffled_random_noisy}}
Let $r$ satisfies \eqref{eq: r} and $H$ satisfies the conditions from Subsection \ref{sec: family}.
Let  $a_0= 2a+2$  and $b_n\sim(\log n)^{1/\beta}$ and $\tilde\theta_n$ be taken from Lemma \ref{lem: mutheta} where $a$ is replaced by $a_0$ and $M$ is assumed to be sufficiently large. We choose $\tilde\theta_n$ independent of the $\epsilon_i$'s and $X_i$'s.  Let $Z_1,\dots,Z_n$ be i.i.d. random variables taking values in $\mathbb{R}$ such that for all fixed $\theta\in\{0,1\}^{b_n}$, the conditionnal distribution of $Z_i$ given $\tilde\theta_n=\theta$ is $\mu_\theta$ with density $f_\theta$ taken from \eqref{eq: ftheta}. 
The $Z$-sample can be assumed to be independent of  the $X$- and $\epsilon$-samples conditionally on $\tilde\theta_n$, which we will assume in the sequel. This implies that the $Z$-sample is independent of the $X$- and $\epsilon$-samples unconditionally. Without loss of generality (by continuity of the distributions in the $X$- and $Z$-samples) we assume that the $Z_i$'s are all different from each other, and that the $X_i$'s are all different from each other.
Denoting by $Z_{(i)}$  the order statistics of the $Z$-sample, and by $\hat F^{-1}_Z$ the corresponding empirical quantile function, we consider the non-decreasing function $m_Z=\hat F_Z^{-1}\circ F_X.$ Hence,
\begin{equation*}
m_Z(X_i)=Z_{(j)}\quad if \quad F_X(X_i)\in\left(\frac{j-1}n,\frac jn\right].
\end{equation*}
In other words, 
\begin{equation*}
m_Z(t)=Z_{(\lceil nF_X(t)\rceil)}\quad for\ all\quad t\in(0,1].
\end{equation*}
For all $j\in\{1,\dots,n\}$ let $n_j$ be the number of indexes $i$ for which $m_Z(X_i)=Z_{(j)}$:
\begin{equation}\label{eq: nj}
n_j=\sum_{i=1}^n1\left(F_X(X_i)\in\left(\frac{j-1}n,\frac jn\right]\right).
\end{equation}
Define the random variables
\begin{equation}\label{eq: ZhatN}
(\hat Z_{(1)},\dots,\hat Z_{(n)})=(Z_{(1)},\dots,Z_{(1)},Z_{(2)},\dots,Z_{(2)},\dots,Z_{(n)},\dots,Z_{(n)})
\end{equation}
where for all $j$, $Z_{(j)}$ is repeated $n_j$ times, see \eqref{eq: nj}. Let $ \rho$ be the permutation on $\{1,\dots,n\}$ such that $Z_{(i)}=Z_{\rho(i)}$ (note that the $Z_i$'s can be assumed to be all different from each other since they are independent with a continuous distribution, whence $\rho$ is uniquely defined). We reconstruct the sample $\hat Z_1,\dots, \hat Z_n$ by defining $\hat Z_{\rho(i)}=\hat Z_{(i)}$ and then define $Y_i=\hat Z_i+\epsilon_i$ for all $i$.

In the sequel, $E_\theta$ denotes expectation under observation of $\cY_n$ and $\cX_{ord}$ as above  conditionally on $\tilde\theta_n=\theta$, and $P_{\tilde\theta_n}$ denotes the distribution of $\tilde\theta_n$.  The minimax risk is bounded from below as follows (see Section \ref{sec: minrandomdeconv} for a proof). Recall  that we call pseudo-estimator an estimator that is allowed to depend on $\mu_X$.
\begin{lemma}\label{lem: minrandomdeconv}
Under the assumptions of Theorem \ref{thm:mlb_shuffled_random_noisy} we have
\begin{eqnarray}\label{eq: minrandomdeconv}
\inf_{\widehat m}\sup_{m\in{\mathcal M}(M,a)}E_{m}\left[\int\vert\widehat m-m\vert d\mu_X\right]
&\geq &
 \inf_{\hat \mu}\int E_\theta W_1(\hat\mu,\mu_\theta)dP_{\tilde\theta_n}(\theta)-Cn^{-1/2}
\end{eqnarray}
where the infimum on the right-hand side is extended to all pseudo-estimators based on observations $N=(n_1,\dots,n_n)$ and $Y_1,\dots,Y_n$. 
\end{lemma}

Since the estimators in the right-hand side no longer depend on the $X$-sample, we can reformulate the model in a simpler way. We still have $Z_1,\dots,Z_n$  with the same distribution $\mu_\theta$, and
$\epsilon_1,\dots,\epsilon_n$  with distribution $\mu_\epsilon$ and we now have $N$ with multinomial distribution with parameter $n$ and probabilities $(n^{-1},\dots,n^{-1}).$ All the variables are mutually independent. Let
\begin{equation*}
(\hat Z_{1},\dots,\hat Z_{n})=(Z_{1},\dots,Z_{1},Z_{2},\dots,Z_{2},\dots,Z_{n},\dots,Z_{n})
\end{equation*}
where for all $j$, $Z_{j}$ is repeated $n_j$ times. On observes $N$ and $Y_i=\hat Z_i+\epsilon_i$ for all $i$.

Henceforth, our goal is to obtain a lower bound on
\begin{equation}\label{eq: W1shuf}
\int E_\theta  W_1(\hat\mu,\mu_\theta) dP_{\tilde\theta_n}(\theta) = \int_{\mathbb{R}}\int E_\theta|\hat F(t)-F_\theta(t)|dP_{\tilde\theta_n}(\theta) dx
\end{equation}
where $\hat F$ and $F_\theta$ denote the distribution functions associated to $\hat\mu$ and $\mu_\theta$ respectively, see \eqref{eq: W1dir}.

Let  $c_k=\sum_{j\leq k} n_j$ for all $k\in\{1,\dots,n\}$ and $c_0=0$. Then, conditionally on $N$ and on $\tilde\theta_n$, the vectors $(Y_{c_{j-1}+1},\dots,Y_{c_{j}})$, $j\in\{1,\dots,n\}$,  are mutually independent
with density
\begin{eqnarray*}
h_{\theta,n_j}(y_{1},\dots,y_{n_{j}})= \int_{\mathbb{R}}f_\theta(z)f_\epsilon(y_1-z)\times\dots\times f_\epsilon(y_{n_j}-z)dz
\end{eqnarray*}
if $n_j>0$. If $n_j=0$ we set $h_{\theta,n_j}(y_{1},\dots,y_{n_{j}}) =1$.
Similar as  \cite[Equation (14)]{dedecker2013minimax} we have for all 
estimators $\hat\mu$ and $s\in\{1,\dots,b_n\}$: 
\begin{align}
\label{eq:chi_square_bound_1}
&\int E_\theta|\hat F(t)-F_\theta(t)|dP_{\tilde\theta_n}(\theta)  \notag \\
& \ge \frac{|H^{(-1)}(b_n(t-t_{s,n})|}{2 b_n} \bbE\int_{\mathbb{R}^n}\prod_{j=1}^n\min\left( h_{\tilde\theta_n,n_j,s,0}, h_{\tilde\theta_n,n_j,s,1}\right)(y_{c_{j-1}+1},\dots,y_{c_{j}})dy_1\dots dy_n \notag \\
& \geq \frac{|H^{(-1)}(b_n(t-t_{s,n})|}{2 b_n} \bbE\prod_{j=1}^n \int_{\mathbb{R}^{n_j}}\min\left( h_{\tilde\theta_n,n_j,s,0}, h_{\tilde\theta_n,n_j,s,1}\right)(y_{1},\dots,y_{n_{j}})dy_1\dots dy_{n_j} \notag \\
& \ge \frac{|H^{(-1)}(b_n(t-t_{s,n})|}{4 b_n} \bbE\left[\prod_{j\ s.t.\ n_j>0} \left(1-\frac 12 \chi^2( h_{\tilde\theta_n,n_j,s,0}, h_{\tilde\theta_n,n_j,s,1})\right)^2\right] \,,
\end{align}
where, we use \cite[Equations (2.12) and (2.21)]{tsybakov2009intro} for the last inequality. Therefore, we need to bound (from above) the $\chi^2$-divergence. We now present a lemma, whose proof is deferred to Appendix \ref{sec:chi_square}. 
\begin{lemma}
\label{lem:chi_square}
Suppose that $\delta$ is standard Gaussian and $b_n\geq 1$. Then for each $j$ such that $n_j > 0$, the $\chi^2$-divergence in the above equation can be upper bounded as: 
$$
\int \frac{(h_{\theta, n_j, s, 1} - h_{\theta, n_j, s, 0})^2}{ h_{\theta, n_j, s, 0}} (y_1, \dots, y_{n_j}) \ dy_1 \dots dy_{n_j} \le C e^{-c\frac{\sigma_n^2b_n^2}{n_j}} \,.
$$
for some positive constants $C, c$. 
\end{lemma}
\noindent
Using the bound obtained in the lemma above in equation \eqref{eq:chi_square_bound_1} we get: 
\begin{align}\label{eq:conj}\notag
&\int E_\theta|\hat F(t)-F_\theta(t)|dP_{\tilde\theta_n}(\theta) \\
& \ge  \frac{|H^{(-1)}(b_n(t-t_{s,n})|}{4 b_n} \bbE\left[\prod_{j\ s.t.\ n_j>0} \left(1-\frac 12 C e^{-c\frac{\sigma_n^2b_n^2}{n_j}} \right)^2\right] \\
\notag
& \ge \frac{|H^{(-1)}(b_n(t-t_{s,n})|}{4 b_n} \bbE\left[\prod_{j\ s.t.\ n_j>0} \left(1-\frac 12 C e^{-c\frac{\sigma_n^2b_n^2}{n_j}} \right)^2\mathds{1}_{\max_j n_j \le K\frac{\log{n}}{\log{\log{n}}}}\right] \\
\notag
& \ge \frac{|H^{(-1)}(b_n(t-t_{s,n})|}{4 b_n} \bbE\left[\prod_{j\ s.t.\ n_j>0} \left(1-\frac 12 C e^{-\frac{c}{K}\sigma_n^2b_n^2 \frac{\log{\log{n}}}{\log{n}}} \right)^2\mathds{1}_{\max_j n_j \le K\frac{\log{n}}{\log{\log{n}}}}\right]\,.
\end{align}
Now taking $b_n = \sigma_n^{-1}\sqrt{K/c_1}(\log{n}/\sqrt{\log{\log{(n)}}})$ we have: 
\begin{align*}
&\int E_\theta|\hat F(t)-F_\theta(t)|dP_{\tilde\theta_n}(\theta) \\
& \ge C_2 \sigma_n |H^{(-1)}(b_n(t-t_{s,n})| \frac{\sqrt{\log{\log{n}}}}{\log{n}} \left(1- \frac{C_1}{2n}\right)^{2n} \bbP\left(\max_j n_j \le K\frac{\log{n}}{\log{\log{n}}}\right)\\
& \ge C_3\sigma_n |H^{(-1)}(b_n(t-t_{s,n})| \frac{\sqrt{\log{\log{n}}}}{\log{n}}\,,
\end{align*}
where we use Lemma \ref{lem: multinomial} below for the last inequality. Although the lemma is a standard tail bound since $N$ has a multinomial distribution, we present a simple proof for the convenience of the reader in Section \ref{sec: multinomial}.
\begin{lemma}\label{lem: multinomial}
There exists $K>0$ such that
\begin{equation*}
\lim_{n\to\infty}\bbP\left(\max_j n_j > K\frac{\log{n}}{\log{\log{n}}}\right)=0.
\end{equation*}
\end{lemma}
Integrating with respect to $t$ then yields
\begin{align*}
\int_{\mathbb{R}}\int E_\theta|\hat F(t)-F_\theta(t)|dP_{\tilde\theta_n}(\theta)dt & \ge C_3\sigma_n\frac{\sqrt{\log{\log{n}}}}{\log{n}}\sum_{s=1}^{b_n}\int_{t_{s,n}}^{t_{s+1,n}} |H^{(-1)}(b_n(t-t_{s,n})| dt\\
& \ge C_3\sigma_n\frac{\sqrt{\log{\log{n}}}}{\log{n}}\int_0^1 |H^{(-1)}(u)| du.
\end{align*}
The integral on the last line is strictly positive since we have assumed that condition (A1) in the proof of Theorem 3 of \cite{dedecker2013minimax} was satisfied. Combining with \eqref{eq: W1shuf} and  Lemma \ref{lem: minrandomdeconv} proves that
\begin{eqnarray*}
\inf_{\widehat m}\sup_{m\in{\mathcal M}(M,a)}E_{m}\left[\int\vert\widehat m-m\vert d\mu_X\right]
&\geq &
 \inf_{\hat \mu}\int_{\mathbb{R}}\int E_\theta|\hat F(t)-F_\theta(t)|dP_{\tilde\theta_n}(\theta) dt-Cn^{-1/2}\\
& \ge& c\sigma_n\frac{\sqrt{\log{\log{n}}}}{\log{n}}-Cn^{-1/2}.
\end{eqnarray*}
The term on the right-hand side is larger than $(c/2)\sigma_n\frac{\sqrt{\log{\log{n}}}}{\log{n}}$ provided that $c\sigma_n\frac{\sqrt{\log{\log{n}}}}{\log{n}}\ge 2Cn^{-1/2}$, which 
completes the proof of Theorem \ref{thm:mlb_shuffled_random_noisy}.  \hfill{$\Box$} 

\subsection{Proof of Theorem \ref{thm:mlb_shuffled_random_noisy2}}
\subsubsection{Lower bound for $\mu_X$}
If conjecture \ref{conj} is true, then \eqref{eq:conj}
yields that for $b_n=c\sigma_n^{-1}(\log n)^{1/2}$ with some sufficiently large $c>0$ one has
\begin{align*}
\int E_\theta|\hat F(t)-F_\theta(t)|dP_{\tilde\theta_n}(\theta)  \ge C \frac{|H^{(-1)}(b_n(t-t_{s,n})|}{ b_n} \,.
\end{align*} 
Integrating with respect to $t$ then yields
\begin{align*}
\int_{\mathbb{R}}\int E_\theta|\hat F(x)-F_\theta(x)|dP_{\tilde\theta_n}(\theta)dt & \ge C_3\sigma_n(\log n)^{-1/2}\sum_{s=1}^{b_n}\int_{t_{s,n}}^{t_{s+1,n}} |H^{(-1)}(b_n(t-t_{s,n})| dt.
\end{align*}
Similar to the  proof of Theorem \ref{thm:mlb_shuffled_random_noisy}, the sum on the last line is strictly positive and  Theorem \ref{thm:mlb_shuffled_random_noisy2} for the case where $\mu=\mu_X$ follows from similar arguments as for the end of the  proof of Theorem \ref{thm:mlb_shuffled_random_noisy}.

\subsubsection{Lower bound for $\hat \mu_X$}
For arbitrary estimator $\widehat m$, let $\widetilde m$  be an estimator in $\mathcal M(M,a)$ be defined as in Section \ref{sec: prooflbemp} for some arbitrarily small $\eta>0$. As in Section \ref{sec: prooflbemp}  we get 
\begin{eqnarray*}
\int|\widetilde m-m|d\hat\mu_X 
&\leq&2 \int|\widehat m-m|d\hat\mu_X +\eta
\end{eqnarray*}
and therefore,
\begin{equation*}
\mathcal R_S(M,a,n, \hat \mu_X) \geq \frac 12\left(\inf_{\widetilde m \in\mathcal M(M,a)}E_{m}\int|\widetilde m-m|d\hat\mu_X -\eta\right).
\end{equation*}
Combining with  Lemma \ref{lem: emptopop} yields that there exists $C>0$ such that
\begin{equation*}
\mathcal R_S(M,a,n, \hat \mu_X) \geq \frac 12\left(\inf_{\widetilde m \in\mathcal M(M,a)}E_{m}\int|\widetilde m-m|d\mu_X-Cn^{-1/2} -\eta\right).
\end{equation*}
We choose $\eta=Cn^{-1/2}$. Then  it follows from the lower bound already obtained for $\mu_X$ that there exists $c>0$ such that 
\begin{equation*}
\mathcal R_S(M,a,n, \hat \mu_X) \geq 2c\sigma_n\left(\log n\right)^{-1/2}-Cn^{-1/2}.
\end{equation*}
The right-hand side above is larger than $c\sigma_n\left(\log n\right)^{-1/2}$ provided that
$
c\sigma_n\geq Cn^{-1/2}(\log n)^{1/2}$. This completes the proof of Theorem \ref{thm:mlb_shuffled_random_noisy2} for the case where $\mu=\hat\mu_X$.
 \hfill{$\Box$} 

\section{Proof of the lemmas}
\label{app:lemmas}

\subsection{Proof of Lemma \ref{lem: entropy}}

Let $m\in\mathcal M(M,a)$ and $m_-=(-m)\vee 0$. Then $m_-$ is non-negative and non-increasing whence for all $x_0$ in the support of $\mu_X$,
\begin{eqnarray*}
M&\geq&\int_0^{x_0}(m_-(x))^{a+2}d\mu_X(x)\\
&\geq&(m_-(x_0))^{a+2}\int_0^{x_0}d\mu_X(x)\\
&=&(m_-(x_0))^{a+2}F_X(x_0)
\end{eqnarray*}
with $F_X$ being the distribution function of $X$. Therefore, 
\begin{equation}\label{eq: env+}
m_-(x_0)\leq\left(\frac{M}{F_X(x_0)}\right)^{1/(a+2)}
\end{equation}
for all $x_0$ in the support of $\mu_X$. The right-hand side is larger than 1 for all $x_0$ provided that $M\geq 1$, and it defines a non-increasing function. Moreover, for $\alpha\in(0,a)$ we have
\begin{eqnarray*}
\int_0^1 \left(\frac{M}{F_X(x)}\right)^{(2+\alpha)/(a+2)}d\mu_X(x)&=& \int_0^1 \left(\frac{M}{x}\right)^{(2+\alpha)/(a+2)}dx\\
&<&\infty,
\end{eqnarray*}
using that $F_X$ is continuous for the first equality.
Hence, it follows from \cite[Lemma 7.11]{van2000empirical} that for some $A$ (that depends only on $a$ and $M$) and all $\delta>0$,
\begin{equation*}
\log N_{[]}(\delta,\left\{m_-\mbox{ s.t. }m\in\mathcal M(M,a)\right\},\|\,.\,\|_2)\leq A\delta^{-1}.
\end{equation*}
With $m_+=m\vee0$  one obtains similarly that 
\begin{equation}\label{eq: env-}
m_+(x_0)\leq\left(\frac{M}{1-F_X(x_0)}\right)^{1/(a+2)}
\end{equation}
for all $x_0$ in the support of $\mu_X$ and therefore,
\begin{equation*}
\log N_{[]}(\delta,\{m_+\mbox{ s.t. }m\in\mathcal M(M,a)\},\|\,.\,\|_2)\leq A\delta^{-1}.
\end{equation*}
for all $\delta>0$. Lemma \ref{lem: entropy} easily follows from the two entropy bounds above.  \hfill{$\Box$}

\subsection{Proof of Lemma \ref{lem: emptopop} }

It follows from \eqref{eq: env+} and  \eqref{eq: env-} that the set of functions $\mathcal M(M,a)$ has enveloppe
\begin{equation}\label{eq: envelope}
e (x)= \left(\frac{M}{F_X(x)(1-F_X(x))}\right)^{1/(a+2)}
\end{equation}
that satisfies
\begin{eqnarray*}
\int e^2d\mu_X&=& M^{2/(a+2)}\int_0^1\frac{1}{(x(1-x))^{2/(a+2)}}dx\\
&<&\infty,
\end{eqnarray*}
using that $F_X$ is continuous for the equality and that $a>0$ for the inequality.
Therefore, it follows from \cite[Theorem 2.14.2]{wellner2013weak} that
\begin{eqnarray*}
E\left|\int\sup_{f\in \mathcal M(M,a)}\vert f-m\vert (d\mu_X-d\hat \mu_X)\right|\lesssim n^{-1/2} \int_0^1\sqrt{1+\log N_{[]}(\delta C,\mathcal M_0,\|\,.\,\|_2)}d\delta
\end{eqnarray*}
where $C=\int e^2d\mu_X$ and $\mathcal M_0$ is the set of all functions $|f-m|$ with $f\in\mathcal M(M,a)$. It follows from Lemma \ref{lem: entropy} that the above integral is finite,
which completes the proof of Lemma \ref{lem: emptopop}.
\hfill{$\Box$}

\subsection{Proof of Lemma \ref{lem: mutheta}}
That $f_\theta$ in \eqref{eq: ftheta} is the density of a probability measure $\mu_\theta\in\mathcal D(M,a)$ for all $\theta$ provided that $M$ is sufficiently large and $b_n\sim(\log n)^{1/\beta}$ has been proved before the statement of the lemma. Hence, we turn to the proof of the inequality.

For $\theta\in\{0,1\}^{b_n}$, $s\in\{1,\dots,b_n\}$ and $j\in\{0,1\}$, let us define the measure $\mu_{\theta,s,j}$ with density $f_{\theta,s,j}:=f_{\theta_1,\dots,\theta_{s-1},j,\theta_{s+1},\dots,\theta_{b_n}}$. Let $ h_{\theta,s,j}= f_{\theta,s,j}\star \mu_\epsilon$ where $\star$ denotes convolution. With exactly the same arguments as in \cite[Section 2.2.1]{dedecker2013minimax} (with $p=d=A=1$), and using their notations, one obtains that if for some $b_n>0$ and $c>0$ one has
\begin{equation}\label{eq: DM16}
\chi^2(h_{\theta,s,0},h_{\theta,s,1}):=\int_\R\frac{(h_{\theta,s,0}-h_{\theta,s,1})^2}{h_{\theta,s,0}}\leq cn^{-1}
\end{equation}
for all $\theta$, 
then for sufficiently large $M$, there exists a constant $C>0$ such that for all  $n$,
\begin{eqnarray}\notag \label{eq: DM16'}
\inf_{\tilde \mu_n}\mathbb{E}_{\tilde\theta_n}E_{(\mu_{\tilde\theta_n}\star\mu_\epsilon)^{\otimes n}}W_1(\mu_{\tilde\theta_n},\tilde\mu_n)
&\geq& \frac C{b_n}\sum_{s=1}^{b_n}\int_{\sigma_nt_{s,n}}^{\sigma_nt_{s+1;n}}|H^{(-1)}(b_n(\sigma_n^{-1}t-t_{s,n}))|dt\\
&\geq& \frac {C\sigma_n}{b_n}\int_0^1 |H^{(-1)}(t)|dt
\end{eqnarray}
where  the infimum is over all estimators $\tilde \mu_n$.

Now, we choose $b_n$ as small  as possible  while satisfying \eqref{eq: DM16}. For this task, we follow the line of proof of \cite[Equation (16)]{dedecker2013minimax}. Similar to \cite[Equation (17)]{dedecker2013minimax}, we obtain that there exists $C,C'>0$ such that for arbitrary $\theta$

\begin{eqnarray*}
\chi^2(h_{\theta,s,0},h_{\theta,s,1})
&\leq& 
\int_\R\frac{\{\sigma_n^{-1}\int_\R H(b_n(\sigma_n^{-1}(t-x)-t_{s,n}))f_\epsilon(x)dx\}^2}{f_{\theta,s,0}\star f_\epsilon(t)}dt\\
&\leq&C\int_\R\frac{\{\sigma_n^{-1}\int_\R H(b_n(\sigma_n^{-1}(t-x)-t_{s,n}))f_\epsilon(x)dx\}^2}{f_{0}\star f_\epsilon(t)}dt\\
&\leq&C'\int_\R\frac{\{\sigma_n^{-1}\int_\R H(b_n\sigma_n^{-1}(t-x))f_\epsilon(x)dx\}^2}{f_{0}\star f_\epsilon(t)}dt
\end{eqnarray*}
where $f_\epsilon$ denotes the density of $\epsilon$ and $f_0(t)=\sigma_n^{-1}f_{0,r}(\sigma_n^{-1}t)$ for all $t$. Hence, With $f_\delta$ the density of $\delta=\epsilon/\sigma_n$ given by $f_\delta(t)=\sigma_nf_\epsilon(\sigma_n t)$  we conclude that
\begin{eqnarray*}
\chi^2(h_{\theta,s,0},h_{\theta,s,1})
&\leq& 
 C'\sigma_nb_n^{-1}\int_\R\frac{\{\sigma_n^{-1}\int_\R H(v-y)f_\delta(y/b_n)dy/b_n\}^2}{f_{0}\star f_\epsilon(v\sigma_n/b_n)}dv.
\end{eqnarray*}
For the denominator we have
\begin{eqnarray*}
f_{0}\star f_\epsilon(v\sigma_n/b_n)&=&\int_\mathbb{R}f_{0}\left(\frac{v\sigma_n}{b_n}-x\right)f_\epsilon(x)dx\\
&=&\sigma_n^{-2}\int_\mathbb{R}f_{0,r}\left(\frac{v}{b_n}-\frac{x}{\sigma_n}\right)f_\delta\left(\frac{x}{\sigma_n}\right)dx\\
&=&\sigma_n^{-1}f_{0,r}\star f_\delta(v/b_n).
\end{eqnarray*}
Hence,
\begin{eqnarray*}
\chi^2(h_{\theta,s,0},h_{\theta,s,1})
&\leq& C' b_n^{-1}\int_\R\frac{\{\int_\R H(v-y)f_\delta(y/b_n)dy/b_n\}^2}{f_{0,r}\star f_\delta(v/ b_n)}dv.
\end{eqnarray*}
Similar to the two last displays on page 284 of \cite{dedecker2013minimax}, we conclude that we can choose $b_n$ of the order $(\log n)^{1/\beta}$ so that $\chi^2(h_{\theta,s,0},h_{\theta,s,1})=O(n^{-1})$. 
Hence, \eqref{eq: DM16} holds  so  \eqref{eq: DM16'} holds as well. since we have assumed that $H$ satisfies Assumption (A1)  taken from  \cite[page 282]{dedecker2013minimax}, the integral in \eqref{eq: DM16'} is strictly positive,  so the lemma follows.
\hfill{$\Box$}

\subsection{Proof of Lemma \ref{lem: Ftheta}}\label{sec: proofFtheta}
Denoting by $F_\theta$ the distribution function associated to $\mu_\theta$, we first show that there exists $C>0$ such that for all $n$,
\begin{equation}\label{eq: empZtopopZ}
\sup_{\theta\in\{0,1\}^{b_n}}\int_{-\infty}^\infty\sqrt{F_\theta(t)(1-F_\theta(t))}dt\leq C.
\end{equation} 

With $F_\theta^1$ the distribution function of $\mu_\theta$ in the particular case where $\sigma_n=1$ we have $F_\theta(t)=F_\theta^1(\sigma_n^{-1}t)$ for all $t\in\mathbb{R}$, so using the change of variable $x=\sigma_n^{-1}t$, it suffices to prove the result in the particular case where $\sigma_n=1$ for all $n$, which we assume in the sequel.

 It follows from Lemma 3 in Appendix A of \cite{dedecker2013minimax} that there exists $C>0$ such that 
\begin{equation*}
|f_\theta(t)|\leq C(1+t^2)^{-r}
\end{equation*}
for all $t\in\mathbb{R}$ and $\theta \in\{0,1\}^{b_n}$. Using that $F_\theta(t)(1-F_\theta(t))$ is less than or equal to $F_\theta(t)$ for all $t<-1$ and to $1-F_\theta(t)$ for all $t>1$ we conclude that  for all $t$ with $|t|>1$ and $\theta \in\{0,1\}^{b_n}$, one has
\begin{eqnarray*}
F_\theta(t)(1-F_\theta(t))
&\leq& \int_{|t|}^\infty C(1+x^2)^{-r}dx\\
&\leq& \int_{|t|}^\infty Cx^{-2r}dx\\
&\leq& \frac{C|t|^{1-2r}}{2r-1}
\end{eqnarray*}
where we recall that $r>1/2$. Hence,
\begin{equation*}
\int_{-\infty}^\infty\sqrt{F_\theta(t)(1-F_\theta(t))}dt \leq 2+2\int_1^\infty \sqrt{\frac{Ct^{1-2r}}{2r-1}}dt.
\end{equation*}
The integral on the right-hand side is finite since we have chosen $r>3/2$. This proves the inequality in \eqref {eq: empZtopopZ}. 

Now, combining \eqref{eq: empZtopopZ} with \cite[Theorem 3.2]{bobkov2019one} proves that there exists $C>0$ that does not depend on $\theta$ nor $n$ such that
\begin{eqnarray*}
\int W_1(\hat\mu_Z,\mu_\theta)d\mu_\theta^{\otimes n}(Z_1,\dots,Z_n)\leq Cn^{-1/2}
\end{eqnarray*}
which completes the proof of Lemma \ref{lem: Ftheta}.  \hfill{$\Box$}

\subsection{Proof of Lemma \ref{lem:order_stat}}

See e.g. \cite[Chapter 21]{convergence1996empirical}.

\subsection{Proof of Lemma \ref{lem: minrandomdeconv}}
\label{sec: minrandomdeconv}

First, we prove that 
\begin{equation}\label{eq: devZ}
P\left(\frac 1n\sum_{i=1}^n|Z_i|^{a+2}>M\right)\leq 4 n^{-1} M^{-1}.
\end{equation}
Recall that from   Lemma \ref{lem: mutheta} where $a$ is replaced by $a_0= 2a+2$, it follows that $\mu_\theta\in\mathcal D(M,a_0)$ for all $\theta\in\{0,1\}^{b_n}$. The conditional distribution of $Z$ given $\tilde\theta_n=\theta$ is $\mu_\theta$ so denoting by $\mathbb{E}_{\tilde\theta_n}$
 the expectation with respect to $\tilde\theta_n$ we have
\begin{eqnarray}\label{eq: espZ}\notag
E|Z|^{2a+4}&=&\mathbb{E}_{\tilde\theta_n}\int |z|^{2a+4}d\mu_{\tilde\theta_n}(z)\\
\notag
&=&\mathbb{E}_{\tilde\theta_n}\int |z|^{a_0+2}d\mu_{\tilde\theta_n}(z)\\
&\leq& M.
\end{eqnarray}
It follows from Jensen's inequality that $E|Z|^{a+2}\leq M^{1/2}$ so combining the previous two displays with Markov's inequality yields that for $M\geq 4$ (which implies that $M\geq 2E|Z|^{a+2}$), one has
\begin{eqnarray*}
P\left(\frac 1n\sum_{i=1}^n|Z_i|^{a+2}>M\right)
&\leq &P\left(\frac 1n\sum_{i=1}^n|Z_i|^{a+2}-E|Z|^{a+2}>\frac M2\right)\\
&\leq& \frac{n^{-1} Var(|Z|^{a+2})}{(M/2)^2}\\
&\leq& \frac{n^{-1} E|Z|^{2a+4}}{(M/2)^2}\\
&\leq& 4 n^{-1} M^{-1}.
\end{eqnarray*}
This proves \eqref{eq: devZ}.

Now, for all estimators $\widehat m$ of $m_0$ we have
\begin{eqnarray*}
&&\sup_{m\in{\mathcal M}(M,a)}E_{m}\left[\int\vert\widehat m-m\vert d\mu_X\right]\\
&&\qquad\qquad \geq 
\int 1(m_Z\in{\mathcal M}(M,a))E^Z_{m_Z}\left[\int\vert\widehat m-m_Z\vert d\mu_X\right]d\mu_Z^{\otimes n}(Z_1,\dots,Z_n)
\end{eqnarray*}
where we denote by 
$\mu_Z$ the distribution of $Z_i$ and $E^Z$  the expectation conditionally on the $Z$-sample. 
The next step is to remove the characteristic function in the integral. To do that we first note that $F_X(X_i)$ is uniformly distributed on $[0,1]$ by continuity of $F_X$ so
by definition of $m_Z$, one has
\begin{equation*}
\int |m_Z|^{a+2}d\mu_X=E|\hat F_Z^{-1}(U)|^{a+2}
\end{equation*}
where $Z_1,\dots,Z_n$ are considered as fixed and $U$ is uniformly distributed on $[0,1]$. Since $\hat F_Z^{-1}(U)\sim\hat\mu_Z$, the empirical distribution of the $Z$-sample, this implies that
\begin{eqnarray*}
\int |m_Z|^{a+2}d\mu_X=\int|z|^{a+2}d\hat\mu_Z(z)=\frac 1n\sum_{i=1}^n|Z_i|^{a+2}.
\end{eqnarray*}
Hence, it follows from \eqref{eq: devZ} that
\begin{eqnarray}\label{eq: mZnotinM}
P\left(m_Z\not\in{\mathcal M}(M,a)\right)
&=&P\left(\frac 1n\sum_{i=1}^n|Z_i|^{a+2}>M\right)\leq 4n^{-1}M^{-1}.
\end{eqnarray}
Now, for arbitrary  $\widehat m \in {\mathcal M}(M,a)$ one has
\begin{eqnarray*}
\int|\widehat m-m_Z|d\mu_X
&\leq&\int|\widehat m|d\mu_X+\int_0^1| m_Z|d\mu_X\\
&\leq& M^{1/(a+2)}+\frac 1n\sum_{i=1}^n|Z_i|,
\end{eqnarray*}
using that $\widehat m \in {\mathcal M}(M,a)$ combined with Jensen's inequality for the last line.
It follows from \eqref{eq: espZ} combined with Jensen's inequality that $E|Z|\leq M^{1/(2a+4)}$ so we can find $C>0$ sufficiently large so that
\begin{eqnarray*}
\int|\widehat m-m_Z|d\mu_X&\leq& C+\left|\frac 1n\sum_{i=1}^n|Z_i|-E|Z|\right|.
\end{eqnarray*}
Hence, for all $\widehat m\in{\mathcal M}(M,a)$ we have
\begin{eqnarray*}
&&
\int 1(m_Z\not\in{\mathcal M}(M,a))E^Z_{m_Z}\left[\int\vert\widehat m-m_Z\vert d\mu_X\right]d\mu_Z^{\otimes n}(Z_1,\dots,Z_n)\\
&&\qquad\qquad \leq 
CP\left(m_Z\not\in{\mathcal M}(M,a)\right)+P^{1/2}\left(m_Z\not\in{\mathcal M}(M,a)\right)
E^{1/2}\left(\frac 1n\sum_{i=1}^n|Z_i|-E|Z|\right)^2 
\end{eqnarray*}
using H\"older's inequality for the last line. For the last expectation we have
\begin{eqnarray*}
E^{1/2}\left(\frac 1n\sum_{i=1}^n|Z_i|-E|Z|\right)^2 = 
\sqrt{n^{-1}Var(Z)}
\end{eqnarray*}
where $Var(Z)\leq E|Z|^2\leq M^{1/(a+2)}$. Combining with \eqref{eq: mZnotinM} proves that
\begin{eqnarray*}
\int 1(m_Z\not\in{\mathcal M}(M,a))E^Z_{m_Z}\left[\int_0^1\vert\widehat m-m_Z\vert d\mu_X\right]d\mu_Z^{\otimes n}(Z_1,\dots,Z_n)
 \leq 
Cn^{-1}
\end{eqnarray*}
and therefore, one can remove the characteristic function as required: for all $\widehat m\in\mathcal M(M,a)$,
\begin{eqnarray*}
\sup_{m\in{\mathcal M}(M,a)}E_{m}\left[\int\vert\widehat m-m\vert d\mu_X\right]
&\geq &
\int E^Z_{m_Z}\left[\int\vert\widehat m-m_Z\vert d\mu_X\right]d\mu_Z^{\otimes n}(Z_1,\dots,Z_n)-Cn^{-1}\\ \notag
&=&
\int E^Z_{m_Z}W_1(\mu_{\widehat m},\mu_{m_Z}) d\mu_Z^{\otimes n}(Z_1,\dots,Z_n)-Cn^{-1}
\end{eqnarray*}
where  we used \eqref{eq: L1toW1} for the last equality. By definition of $m_Z$, $\mu_{m_Z}$ is the empirical distribution $\hat\mu_Z$ of the $Z$-sample whence
\begin{eqnarray*}
\sup_{m\in{\mathcal M}(M,a)}E_{m}\left[\int\vert\widehat m-m\vert d\mu_X\right]
&\ge &
\int E^Z_{m_Z}W_1(\mu_{\widehat m},\hat\mu_{Z}) d\mu_Z^{\otimes n}(Z_1,\dots,Z_n)-Cn^{-1}.
\end{eqnarray*}
The conditional distribution of $(Z_1,\dots,Z_n)$ given $\tilde\theta_n=\theta$ is $\mu_\theta^{\otimes n}$ so integrating the inequality in Lemma \ref{lem: Ftheta} with respect to $\theta$ yields that we can find $C>0$ such that  for all $\widehat m\in\mathcal M(M,a)$,
\begin{align}
\label{eq: minrandom}
\sup_{m\in{\mathcal M}(M,a)}E_{m}\left[\int\vert\widehat m-m\vert d\mu_X\right]
\ge
\int E_\theta W_1(\mu_{\widehat m},\mu_\theta) dP_{\tilde\theta_n}(\theta)-Cn^{-1/2}.
\end{align} 

 Similar to the proof of Theorem \ref{theo: DMlowerUL},  we can restrict  to estimators $\widehat m$ for which $E_\theta\int|\widehat m|d\mu_X$ is finite for  almost all $\theta$. Since conditionally on $N:=(n_1,\dots,n_n)$, the $Y$-sample is independent of the $X$-sample,  for all such estimators $\widehat m$, the conditional expectation
\begin{equation*}
\widehat m^{NY}:=E_\theta(\widehat m|N,Y_1,\dots,Y_n)
\end{equation*}
is a function of $N,Y_1,\dots,Y_n$ that depends only on $\mu_X$. Moreover, the Jensen inequality for conditional expectations implies that
\begin{equation*}
E_\theta|\widehat m-m_\theta|\geq E_\theta|\widehat m^{NY}-m_\theta|
\end{equation*}
where for all $\theta$ we set $m_\theta=F_\theta^{-1}\circ F_X$. Since $\mu_\theta=\mu_{m_\theta}$, combining with \eqref{eq: L1toW1} yields
\begin{eqnarray*}
E_\theta W_1(\mu_{\widehat m},\mu_\theta) = \int E_\theta|\widehat m-m_\theta|d\mu_X & \ge \int E_\theta|\widehat m^{NY}-m_\theta|d\mu_X\\
&= E_\theta W_1(\mu_{\widehat m^{NY}},\mu_\theta).
\end{eqnarray*}
Combining with \eqref{eq: minrandom} yields
\begin{eqnarray*}
\sup_{m\in{\mathcal M}(M,a)}E_{m}\left[\int\vert\widehat m-m\vert d\mu_X\right]
&\ge&
\int E_\theta W_1(\mu_{\widehat m^{NY}},\mu_\theta) dP_{\tilde\theta_n}(\theta)-Cn^{-1/2}\\
&\ge&
\inf_{\hat\mu}\int E_\theta W_1(\hat\mu,\mu_\theta) dP_{\tilde\theta_n}(\theta)-Cn^{-1/2}
\end{eqnarray*}
where the infimum is extended to all psudo-estimators.
Taking the infimum on the left hand side of the inequality completes the proof of Lemma \ref{lem: minrandomdeconv}. \hfill{$\Box$}

\subsection{Proof of Lemma  \ref{lem:chi_square}}
\label{sec:chi_square}
Recall the definition of $h_0$: 
\begin{align*}
&h_{0, n_j}(y_1, \dots, y_{n_j})\\
&=  \int_{\mathbb{R}}\frac{C_r}{(1+z^2)^r}\times f_\epsilon(y_1-z)\times\dots\times f_\epsilon(y_{n_j}-z)dz\\
&= \int_{\mathbb{R}}\frac{C_r}{(1+z^2)^r}\times (2\pi\sigma_n^2)^{-n_j/2}\exp\left(-\frac 1{2\sigma_n^2} \sum_j y_j^2 -\frac{n_j}{2\sigma_n^2}z^2+\frac{z}{\sigma_n^2}\sum_{i=1}^{n_j} y_i \right)
dz \\
& = (2\pi\sigma_n^2)^{-n_j/2} e^{-\frac{1}{2\sigma_n^2}\sum_j (y_j - \bar y)^2}\int_{\mathbb{R}}\frac{C_r}{(1+z^2)^r}\times \exp\left(-\frac{n_j}{2\sigma_n^2} (z - \bar y)^2\right)
dz \\
& = (2\pi\sigma_n^2)^{-n_j/2} e^{-\frac{1}{2\sigma_n^2}\sum_j (y_j - \bar y)^2}\int_{\mathbb{R}}\frac{C_r}{(1+(z + \bar y)^2)^r}\times \exp\left(-\frac{n_j}{2\sigma_n^2} z^2\right) 
dz .
\end{align*}
 For the numeration of the chi-square divergence note that: 
\begin{align*}
h_{\theta, n_j, s, 1}(\by) - h_{\theta, n_j, s, 0}(\by) & = \int H(b_n(z - t_{s, n})) f_\eps(y_1 - z) \dots f_\eps(y_{n_j} - z) \ dz \\
& = (2\pi\sigma_n^2)^{-\frac{n_j}{2}}e^{-\frac1{2\sigma_n^2}\sum_j (y_j - \bar y)^2} \int H(b_n(z - t_{s, n})) e^{-\frac{n_j}{2\sigma_n^2}(z - \bar y)^2} \ dz .
\end{align*}
 Therefore, we have: 
 \begin{align*}
 & \int \frac{\left(h_{\theta, n_j, s, 1}(\by) - h_{\theta, n_j, s, 0}(\by)\right)^2}{h_{0, n_j}(\by)} \ d \by\\
 & = \int (2\pi\sigma_n^2)^{-n_j/2} e^{-\frac{1}{2\sigma_n^2}\sum_j (y_j - \bar y)^2} \frac{ \left(\int H(b_n(z - t_{s, n})) e^{-\frac{n_j}{2\sigma_n^2}(z - \bar y)^2} \ dz \right)^2}{\int_{\mathbb{R}}\frac{C_r}{(1+(z + \bar y)^2)^r}\times \exp\left(-\frac{n_j}{2\sigma_n^2} z^2\right) \ dz} \ d\by
 \end{align*}
 Now we do the change of variable $\bs = A \by$ where $A \in \reals^{n_j \times n_j}$ is an orthonormal matrix. The first row of A is $(n_j^{-1/2}, \dots, n_j^{-1/2})$ and the rest of the rows are obtained by applying Gram-Schmidt algorithm. Therefore, we have $s_1 = \sqrt{n_j} \bar y$. Furthermore, $\|\bs\|^2 = \|\by\|^2$ and consequently, $\sum_j (y_j - \bar y)^2 = \sum_j y_j^2 - n_j \bar y^2 = \|\bs\|^2 - s_1^2 = \sum_{j = 2}^{n_j} s_j^2 $.  Using this change of variable we have: 
 \begin{align*}
 & \int \frac{\left(h_{\theta, n_j, s, 1}(\by) - h_{\theta, n_j, s, 0}(\by)\right)^2}{h_{0, n_j}(\by)} \ d \by\\
& = \int (2\pi\sigma_n^2)^{-n_j/2} e^{-\frac{1}{2\sigma_n^2}\sum_{j = 2}^{n_j} s_j^2 } \frac{ \left(\int H(b_n(z - t_{s, n})) e^{-\frac{n_j}{2\sigma_n^2}(z -s_1/\sqrt{n_j})^2} \ dz \right)^2}{\int_{\mathbb{R}}\frac{C_r}{(1+(z + s_1/\sqrt{n_j})^2)^r}\times \exp\left(-\frac{n_j}{2\sigma_n^2} z^2\right) \ dz} \ d\bs \\
& = \frac{1}{\sqrt{2\pi\sigma_n^2}}\int \frac{ \left(\int H(b_n(z - t_{s, n})) e^{-\frac{n_j}{2\sigma_n^2}(z -s_1/\sqrt{n_j})^2} \ dz \right)^2}{\int_{\mathbb{R}}\frac{C_r}{(1+(z + s_1/\sqrt{n_j})^2)^r}\times \exp\left(-\frac{n_j}{2\sigma_n^2} z^2\right) \ dz} \ ds_1 \\
& = \sqrt{\frac{n_j}{2\pi\sigma_n^2}} \int \frac{ \left(\int H(b_n(z - t_{s, n})) e^{-\frac{n_j}{2\sigma_n^2}(z -s_1)^2} \ dz \right)^2}{\int_{\mathbb{R}}\frac{C_r}{(1+(z + s_1)^2)^r}\times \exp\left(-\frac{n_j}{2\sigma_n^2} z^2\right) \ dz} \ ds_1.
 \end{align*} 
For notational simplicity let $\phi_n$ be the density of $\cN(0, n_j^{-1}\sigma_n^2)$, that is
\begin{equation*}
\phi_n(t)=\sqrt{\frac{n_j}{2\pi\sigma_n ^2}} e^{-\frac{n_j}{2\sigma_n ^2}t^2}.
\end{equation*}
We then have
 \begin{align*}
 & \int \frac{\left(h_{\theta, n_j, s, 1}(\by) - h_{\theta, n_j, s, 0}(\by)\right)^2}{h_{0, n_j}(\by)} \ d \by\\
&\qquad\qquad= \int \frac{ \left(\int H(b_n(z - t_{s, n})) \phi_n(z -s_1) \ dz \right)^2}{\int_{\mathbb{R}}\frac{C_r}{(1+(z + s_1)^2)^r} \phi_n(z) \ dz} \ ds_1
 \\
& \qquad\qquad=  \int \frac{ \left(\int H(b_n(z + s_1 - t_{s, n}))  \phi_{n}(z) \ dz \right)^2}{\int_{\mathbb{R}}\frac{C_r}{(1+(z + s_1)^2)^r}  \phi_{n}(z) \ dz} \ ds_1 \\
& \qquad\qquad= \int_0^{\infty}  \frac{ \left(\int H(b_n(z + s_1 - t_{s, n}))  \phi_{n}(z) \ dz \right)^2}{\int_{\mathbb{R}}\frac{C_r}{(1+(z + s_1)^2)^r}  \phi_{n}(z) \ dz} \ ds_1 \\
&\qquad\qquad\qquad\qquad+ \int_{-\infty}^0  \frac{ \left(\int H(b_n(z + s_1 - t_{s, n}))  \phi_{n}(z) \ dz \right)^2}{\int_{\mathbb{R}}\frac{C_r}{(1+(z + s_1)^2)^r}  \phi_{n}(z) \ dz} \ ds_1 \\
& \qquad\qquad\triangleq T_1 + T_2 \,.
 \end{align*} 
Here $T_1$ denotes the integral over the positive real line and $T_2$ denotes the integral over the negative real line. We start with $T_1$ and the calculation for $T_2$ is similar. Before delving deep into the analysis of $T_1$, we give a generic lower bound on the denominator. Considering the two cases on whether $0\leq s\leq 1$ or $s>1$ we obtain:

\begin{align}
\label{eq:lb_s_ge_1}
 \int_{\mathbb{R}}\frac{C_r}{(1+(z + s_1)^2)^r}  \phi_{n}(z) \ dz \notag & \ge \int_{|z + s_1| \le s_1\vee 1}\frac{C_r}{(1+(z + s_1)^2)^r}  \phi_{n}(z) \ dz \notag \\
& \ge (s_1\vee 1)^{-2r} \frac{C_r}{2^r}  \int_{-s_1\vee 1-s_1}^{s_1\vee 1-s_1} \phi_{n}(z) \ dz \notag \\
 & \ge  (s_1\vee 1)^{-2r} \frac{C_r}{2^r} \int_{-1}^0  \phi_{n}(z) \ dz \notag \\
& = (s_1\vee 1)^{-2r} \frac{C_r}{2^r} \left\{\frac12 - \Phi(-\sigma_n^{-1}\sqrt{n_j})\right\} \notag \\
& \ge  (s_1\vee 1)^{-2r} \frac{C_r}{2^r} \left\{\frac12 - \Phi(-1)\right\} 
\end{align}
where $\Phi$ denotes the distribution function of the standard Gaussian distribution. For the last inequality we use that $\sigma_n^{-1}\sqrt{n_j}\geq 1$.

We now establish a bound on $T_1$. Let $\{M_n\}$ be sequence $\ge 1$ and diverges to $\infty$ as $n \uparrow \infty$. We further divide it into two parts: 
$$
T_1 = \int_0^{M_n} + \int_{M_n}^{\infty} \triangleq T_{11} + T_{12}  \,.
$$
For $T_{11}$ observe that: 
\begin{align*}
T_{11} & = \int_0^{M_n}  \frac{ \left(\int H(b_n(z + s_1 - t_{s, n}))  \phi_{n}(z) \ dz \right)^2}{\int_{\mathbb{R}}\frac{C_r}{(1+(z + s_1)^2)^r}  \phi_{n}(z) \ dz} \ ds_1 \notag \\
& \le C\int_0^1 \left(\int H(b_n(z + s_1 - t_{s, n}))  \phi_{n}(z) \ dz \right)^2 \ ds_1 \\
&\qquad+
C\int_1^{M_n} s_1^{2r} \left(\int H(b_n(z + s_1 - t_{s, n}))  \phi_{n}(z) \ dz \right)^2 \ ds_1 \notag \\
& \le C(1+M_n^{2r})\int \left(\int H(b_n(z + s_1 - t_{s, n}))  \phi_{n}(z) \ dz \right)^2 \ ds_1 \notag \\
& \le 2C M_n^{2r} \int \left(H_n \star \phi_n(s_1)\right)^2 \ ds_1 \notag 
\end{align*}
where $H_n(z)=H(b_n(z-t_{s,n}))$ and $\star$ stands for the convolution product. With $H_n^*$ and $\phi_n^*$ the fourier transforms of $H_n$ and $\phi_n$ respectively, the fourier transform of $H_n \star \phi_n$ is the product $H_n^*\times\phi_n^*$ so it follows from the Parceval theorem that
\begin{align*}
T_{11} & \leq 2C M_n^{2r}   \int |H_n^*(t)|^2 |\phi_n^*(t)|^2 \ dt  \notag \\
& \le \frac{2C M_n^{2r} }{b^2_n} \int \left|H^*\left(\frac{t}{b_n}\right)\right|^2 e^{-\frac{\sigma_n^2t^2}{n_j}} \ dt \notag \\
& = \frac{2C M_n^{2r} }{b_n} \int \left|H^*(t)\right|^2 e^{-\frac{\sigma_n^2t^2b_n^2}{n_j}} \ dt \notag \\
& = \frac{2C M_n^{2r} }{b_n} \int_1^2 \left|H^*(t)\right|^2 e^{-\frac{\sigma_n^2t^2b_n^2}{n_j}} \ dt  
\end{align*}
since $H$ satisfies the conditions in Section \ref{sec: family} and in particular, $H^*(t)$ is zero outside $[1,2]$. Hence,
\begin{align*}
T_{11} & \lesssim\frac{M_n^{2r}}{b_n} e^{-\frac{\sigma_n^{2}b_n^2}{n_j}} \,.
\end{align*}
Now for $T_{12}$ we have
\begin{align*}
T_{12} & = \int_{M_n}^\infty  \frac{ \left(\int H(b_n(z + s_1 - t_{s, n}))  \phi_{n}(z) \ dz \right)^2}{\int_{\mathbb{R}}\frac{C_r}{(1+(z + s_1)^2)^r}  \phi_{n}(z) \ dz} \ ds_1 \notag \\
& \le C\int_{M_n}^\infty s_1^{2r} \left(\int H(b_n(z + s_1 - t_{s, n}))  \phi_{n}(z) \ dz \right)^2 \ ds_1 \notag \\
& \le C\int_{M_n}^\infty s_1^{2r} \left(\int_{b_n |s_1 + z - t_{s, n}| \ge s_1^{\alpha_0}}  |H(b_n(s_1 + z - t_{s, n}))| \phi_n(z) \ dz \right. \notag \\
& \qquad \qquad \qquad \left. + \int_{b_n |s_1 + z - t_{s, n}| \le s_1^{\alpha_0}}  |H(b_n(s_1 + z - t_{s, n}))| \phi_n(z) \ dz\right)^2 \ ds_1\notag 
\end{align*}
where $\alpha_0\in(1/2,1)$ will be chosen later. Recall that by assumption, $|H(t)|\leq c(1+t^2)^{-r}$ for all $t\in\mathbb{R}$. Hence, for all $t>1$ we have $|H(t)|\leq c2^{-r}t^{-2r}$ so defining $H_+$ to be the supremum of $|H(t)|$ we get
\begin{align*}
T_{12}& \le C\int_{M_n}^\infty s_1^{2r} \left(c2^{-r} \int_{b_n |s_1 + z - t_{s, n}| \ge s_1^{\alpha_0}} b_n^{-2r} |s_1 + z - t_{s, n}|^{-2r}\phi_n(z) \ dz \right. \notag \\
& \qquad \qquad \qquad \left. + H_+ \int_{b_n |s_1 + z - t_{s, n}| \le s_1^{\alpha_0}} \phi_n(z) \ dz\right)^2 \ ds_1\notag \\
& \le C\int_{M_n}^\infty s_1^{2r} \left(c2^{-2r} s_1^{-2r\alpha_0} + H_+ \int_{b_n |s_1 + z - t_{s, n}| \le s_1^{\alpha_0}} \phi_n(z) \ dz\right)^2 \ ds_1\notag \\
& \le 2Cc^22^{-4r} \int_{M_n}^\infty s_1^{-2r(2\alpha_0 - 1)} \ ds_1 + 2CH_+^2 \int_{M_n}^\infty s_1^{2r} \left(\bbP_{s_1}\right)^2 \ ds_1 
\end{align*}
where
\begin{equation*}
\bbP_{s_1}=\bbP\left(\left|\frac{\sigma_nZ}{\sqrt{n_j}} + s_1 - t_{s,. n}\right| \le \frac{s_1^{\alpha_0}}{b_n}\right).
\end{equation*}
Here, $Z$ denotes a standard Gaussian random variable.
Now, we bound $\bbP_{s_1}$. Before going into the details, note that as long as $b_n > 1$ (which will be immediate from our choice later), $\alpha_0 < 1$  and $s_1 \ge M_n \ge 4$ (for all large $n$ as $M_n$ diverges to infinity), we have $s_1 > (4s_1^{\alpha_0})/b_n$ and by definition $t_{s, n} \le 1 \le s_1/4$. This implies: 
$$
s_1 - \frac{s_1^{\alpha_0}}{b_n} - t_{s, n} \ge \frac{s_1}{2} \,.
$$   
Since $Z$ has the same distribution as $-Z$ this implies that
\begin{align*}
\bbP_{s_1}
&\leq  \bbP\left(-\frac{\sigma_nZ}{\sqrt{n_j}} + s_1 - t_{s,. n} \le \frac{s_1^{\alpha_0}}{b_n}\right)  \\
& \le  \bbP\left(\frac{\sigma_nZ}{\sqrt{n_j}} \ge \frac{s_1}{2}\right) \\
& \le  \bbP\left(Z \ge \frac{s_1}{2}\right) 
\end{align*}
since  $\sqrt {n_j}/\sigma_n\geq 1$. Hence,
\begin{align*}
\bbP_{s_1}
& \le e^{- \frac{s_1^2}{8}}  \le C_1 s_1^{-2r\alpha_0} .
\end{align*}
Here the last line follows from the fact that exponetial decay is faster than polynomial decay and $C_1$ is some constant. Therefore, 
\begin{align*}
T_{12} & \le 2C(c^22^{-4r} +C_1^2H_+^2)\int_{M_n}^\infty s_1^{-2r(2\alpha_0 - 1)} \ ds_1 \\
& \le C_2 M_n^{1 - 2r(2\alpha_0 - 1)}
\end{align*}
provided that $1<2r(2\alpha_0 - 1)$. Since by assumption, $r>3/2$, a sufficient condition is that $2\alpha_0-1>1/3$, that is $\alpha_0>2/3$.
The calculation for $T_2$ (i.e. when $s_1 < 0$) is similar, hence skipped for the brevity. Finally combining the bounds obtained above, we obtain: 
\begin{equation*}
\int \frac{\left(h_{\theta, n_j, s, 1}(\by) - h_{\theta, n_j, s, 0}(\by)\right)^2}{h_{0, n_j}(\by)} \ d \by  \le C_3\left(\frac{M_n^{2r}}{b_n} e^{-\frac{\sigma
_n^2b_n^2}{n_j}} + M_n^{1 - 2r(2\alpha_0 - 1)} \right) \,.
\end{equation*}

If we choose $\alpha_0\in (2/3,1)$ and $M_n = e^{c \frac{\sigma_n^2b_n^2}{n_j}}$ for some constant $c$ such that $2cr < 1$ then we get that: 
$$
\int \frac{\left(h_{\theta, n_j, s, 1}(\by) - h_{\theta, n_j, s, 0}(\by)\right)^2}{h_{0, n_j}(\by)} \ d \by  \le C e^{-c_1\frac{\sigma
_n^2b_n^2}{n_j}} \,.
$$
for some constants $C, c$. This completes the proof of Lemma \ref{sec:chi_square}. \hfill{$\Box$}

\subsection{Proof of Lemma \ref{lem: multinomial}}\label{sec: multinomial}
The proof is based on Chernoff's bound on Binomial random variable (e.g., see Lemma 1.3 of \cite{mitzenmacher2001power}), which states that if $X_1, \dots, X_n$ are i.i.d. Bernoulli random variable with $P(X = 1) = p = 1 - P(X = 0)$, then for any $t \ge np$, we have: 
$$
\bbP\left(\sum_{i = 1}^n X_i \ge t\right) \le \left(\frac{np}{t}\right)^t e^{t - np} \,.
$$
We will now use the fact that $(n_1, \dots, n_n) \sim \mult(n: n^{-1}, \dots, n^{-1})$. Hence, each coordinate $n_j$ is marginally distributed as Binomial random variable with parameter $n$ and probability $n^{-1}$.  Therefore we have: 
\begin{align*}
\bbP\left(\max_j n_j > t\right) & \le \sum_{j = 1}^n \bbP\left(n_j >  t\right) \\
& \le  \sum_{j = 1}^n \left(\frac{1}{t}\right)^t e^{t - 1} = \frac{1}{e} e^{\log{n} + t - t\log{t}}\hspace{0.1in} 
\end{align*}
Now using $t = K\frac{\log{n}}{\log{\log{n}}}$ we have: 
\begin{align*}
\log{n} + t  - t\log{t} & = \log{n} +  K\frac{\log{n}}{\log{\log{n}}} -  K\frac{\log{n}}{\log{\log{n}}}\log{\left( K\frac{\log{n}}{\log{\log{n}}}\right)} \\
& =  \log{n} +  K\frac{\log{n}}{\log{\log{n}}} - K\log{n} - K\frac{\log{n}}{\log{\log{n}}}\log{\left(\frac{K}{{\log{\log{n}}}}\right)} \\
& = -\log{n}\left[K-1 + \frac{K}{{\log{\log{n}}}}\left(1 - \log{K} + \log{\log{\log{n}}}\right)\right]
\end{align*}
Therefore, it is immediate that taking $K = 2$ suffices as for large $n$, the term in square bracket of the above equation is larger than $1/2$. 
This completes the proof of Lemma \ref{lem: multinomial}. \hfill{$\Box$}

 \end{appendix}

\end{document}